\documentstyle{article}
\newtheorem{theorem}{Theorem}[section]
\newtheorem{pro}{Proposition}[section]
\newtheorem{lemma}{Lemma}[section]
\newtheorem{definition}{Definition}[section]
\newtheorem{remark}{Remark}[section]
\newtheorem{cor}{Corollary}[section]

\newcommand{\proof}[1]{\noindent{\it\bf Proof:#1\ }}
\newcommand{\QED}{\hfill$\Box$\medskip}

\newcommand{\Fc}{${\cal F}$-curve }
\newcommand{\Gc}{${\cal G}$-curve }
\newcommand{\ptl}{\partial}
\newcommand{\la}{\langle}
\newcommand{\ra}{\rangle}
\title{Weinstein Conjecture and GW Invariants} 
\author{Gang Liu and Gang Tian\\
Department of Mathematics MIT\\
Cambridge, MA 02139}
\date {Preliminary version\\ August 1997\\ Revised November 1997}
\begin{document}
\maketitle

\section{Introduction}

The purpose of this paper is to provide a new and general method to find
closed orbits for the characteristic foliation on a compact hypersurface of
contact type by using Gromov-Witten invariants. The question on whether such closed
orbits exist has been known as the Weinstein conjecture, proposed in [W]. 
As one of the applications of our method, we completely  solve a stabilized version of this
conjecture in this paper. To describe the conjecture, we need to introduce some basic notations first.

Let $V$ be a connected,
symplectic manifold with a symplectic form $\omega$. A hypersurface
$S$ is said to be of contact type if there exists a vector field $X$
defined on some neighborhood $U$ of $S$ such that (i) $X$ is transversal
to $S$ and (ii) $L_{X}\omega=\omega.$

Now for any hypersurface $S$ in a symplectic manifold $V$, there exists
a 1-dimensional characteristic foliation $\xi$ of $S$ defined by:
$$\xi_{x}=\{v_x\,|\,v_x\in T_xS,\,\omega(v_x, u_x)=0, for\,\, all\,\,
u_x\in T_xS\}$$ for $x\in S.$

The Weinstein conjecture claims that if $S$ is of contact type 
and compact, then $S$ carries at least one closed orbit of $\xi$.
The stabilized version of this conjecture claims the same 
conclusion as above
under the assumption that $S$ is contained in
$(V\times {\bf C}^n, \omega \oplus \omega_{0}),$
the stabilization of $V.$

Before we state our result, we
recall that  given $A\in H_{2}(V, {\bf Q}),$
the  $(n+2)$-pointed
GW invariant  is a homomorphism
$$ \Psi^V_{A, g, n+2}:  H_{*}({\bar {\cal M}}_{g, n+2},{\bf Q})
\times H_*(V, {\bf Q})^{n+2} 
\longrightarrow {\bf Q}, $$
(see [FO] and [LT]).
Here for convenience, we use homology instead of cohomology
as in [FO] and [ LT].
We will omit the upscript $V$ if no confusion arises.  

 Throughout this paper,   we  will assume that $ S$
separates $V$, i.e.
 there
 exist two
submanifolds $V_-$ and $V_+$ of $V$
 with common boundary $S$ such that
 $V_-\cup V_+=V,V_-\cap V_+=S.$
  This can be achieved by imposing, for example, that $H^1(V,{\bf Z}_2)=0$.

The main theorems of this paper are stated as follows. 

\begin{theorem}\label{first}
 
 If there exist $A\in H_2(V,{\bf Z})$ and
$\alpha_+,\alpha_-,
 \in H_*(V,\bf Q)$,          such
that 

\noindent
(i) $supp (\alpha_+)\hookrightarrow \stackrel{\circ}{V_-},$ and
$supp (\alpha_-)\hookrightarrow \stackrel{\circ}{V_+}$;

\noindent
(ii) the GW-invariant 
$\Psi_{A, g, n+2}(C;\alpha_-,\alpha_+, \beta _1, \cdots, \beta_n)$
$\not =0$,

then $S$ carries at least one closed orbit of $\xi.$

\end{theorem}

In particular, we have

\begin{theorem}\label{second}

Let $S$ be as above. If there exist $A\in H_2(V,{\bf Z})$  such that
the invariant $\Psi_{A,g, n+2}(\cdot ; e, e, \cdots)\not = 0,$ 
where $e$ denotes the generator of $H_0(V, {\bf Z})$ represented by a point,
then $S$ carries at least
one closed orbit of $\xi.$ 
\end{theorem}

Among various potential applications of these two theorems,
we only mention the following corollaries.

As a corollary to Theorem \ref{first}, we have completely solved the stabilized Weinstein
conjecture in the following theorem.

\begin{theorem}\label{third}
The Weinstein conjecture holds
for $(V\oplus {\bf C}^l, \omega\oplus\omega_0).$
That is, after $V$ is stabilized
by ${\bf C}^l$, the  Weinstein conjecture holds.

\end{theorem}

As a corollary to Theorem \ref{second}, we have

\begin{theorem}

The Weinstein conjecture holds for $\prod _{i=1}^k{\bf C}P^{n_i}$
with the symplectic form $\omega=
\oplus_{i=1}^k\omega_i$, where $\omega_i$ is the standard symplectic
from of ${\bf C}P^{n_i}.$ Moreover, the Weinstein conjecture holds
for any rational algebraic manifolds $(V,\omega)$, provided there is a surjective
morphism $\pi:$ $ V \rightarrow CP^n$  such that $\pi$ is one to one
over $ V\setminus S$ for some subvariety $S$ of $V$ with $codim_{\bf C}\pi(S) 
\geq 2.$
 In particular,the  Weinstein
conjecture holds for any blow-ups of ${\bf C}P^n$ along its subvarieties.

\end{theorem}

A special case of this theorem, where $V={\bf C}P^n$,
was proved by Hofer and Viterbo in [HV].

Closely related to this conjecture is the existence of closed orbits of some
Hamiltonian function, which can be described as follows.

Let $\Psi :S\times (-\epsilon, \epsilon)\rightarrow U\hookrightarrow
V$ be the flow of the vector field $X$. Since $X$ is transversal to $S$,
$\Psi$ is a diffeomorphism from $S\times (-\epsilon, \epsilon)$
to some neighborhood $W$ of $S$ in $V$. Let $S_t=\Psi(S
\times \{t\}),$  $W_{-}=\cup_{t<0}S_t$ and $W_+=\cup_{t>0}S_t.$
Note that $S=S_0.$ Then $W=S\cup W_+\cup W_-.$

Because of our assumption that $S$ separates $V,$
we may
assume further that 
  
\[
(*)\left\{
\begin{array}{cl}
 & \mbox{ there exist two submanifolds $V_-$ and $V_+$ of $V$}\\
 & \mbox{ with common boundary $S$ such that}\\
& \mbox{(i) $V_-\cup V_+=V,V_-\cap V_+=S;$}\\
 &\mbox{(ii) $W_-\hookrightarrow V_-, W_+\hookrightarrow V_+.$}
\end{array}
\right.
\]
This condition implies that $S$ can be realized as a zero set of some
Hamiltonian function.


A particular defining  Hamiltonian function ${\tilde H}={\tilde H}_{S, X}$ of $S$ can
be defined as follows.

\[
{\tilde H}(x)=\left\{\begin{array}{cl}
\epsilon, & {x\in V_+ \setminus W_+}\\
\phi(t) & {x\in S_t}\\
-\epsilon & {x\in V_- \setminus W_-},
\end{array}
\right.
\]
where $\phi :[-\epsilon,\epsilon]\rightarrow [-\epsilon,\epsilon] $ is 
a smooth function defined by
\[
\phi(t)=\left\{\begin{array}{cl}
t  & -\epsilon+2\delta< t< \epsilon-2\delta\\
-\epsilon &{ t< -\epsilon+\delta}\\
\epsilon &{t>\epsilon-\delta},
\end{array}
\right.
\]
for some $0<\delta<<\epsilon.$
The Hamiltonian vector field $X_{{\tilde H}}$ is defined by 
$$\omega(X_{{\tilde H}},\cdot)=d{\tilde H}.$$  
Consider the Hamiltonian equation 
\begin{equation}
\frac{dx}{dt}=X_{{\tilde H}}(x(t)).
\label {1}
\end{equation} 
Any non-trivial closed orbit
$x$ of (1) will lie on some level hypersurface $S_t={\tilde H}^{-1}(t),$
$ -\epsilon +\delta< t<\epsilon-\delta.$ Now the condition ${\cal L}_X\omega
=\omega$ implies that the characteristic foliation $\xi_t$ on $S_t$
is conjugate to $\xi=\xi_0$ on $S$ under the flow $\Psi.$ 
It follows that 
$S$ will also carry a closed orbit of $\xi$ given by $\Psi^{-1}_t(x).$

Therefore, the Weinstein conjecture for those compact hypersurfaces of contact
type satisfying (*) can be proved as long as the existence of some non-trivial
closed orbits of (1) can be established.


The Weinstein conjecture was first proved for a convex 
or star-shaped hypersurface in
$({\bf R}^{2n}, \omega_0)$ by Weinstein and Rabinowitz in [W] and [R] respectively. In 1986, a substantial progress was made by
Viterbo in [V]. He proved the conjecture for any compact 
hypersurface of contact type of $({\bf R}^{2n},\omega_0)$. 
A simpler proof of this was 
given by Hofer and Zehnder in [HZ]. We notice here that for any hypersurface
of ${\bf R}^{2n}$, the above condition (*) always holds. Due to the  
 work of Gromov and 
Floer,  it is  possible to generalize this result to hypersurfaces
in certain general symplectic manifolds. In [FHV], Floer, Hofer and Viterbo 
proved the  Weinstein conjecture for $M\times {\bf R}^{2n}$ with symplectic
form $\omega\oplus\omega_0$ under the assumption that $\pi_2(M)=0.$ 
Note that any compact hypersurface $S$ of $M\times {\bf R}^{2n}$ 
can be embedded into $M\times \{S^2_r\}^n$ for some large $r$, 
where $S^2_r$ is the 2-dimensional sphere
of radius $r$ with the standard symplectic form $\omega_0$ given by the
area form.
In [HV], Hofer and Viterbo proved the same statement 
under the weaker but rather technical
assumption that $\min \omega (A)>\omega_0([S^2_r])$ for all effective classes  $A\in H_2(M).$
Here a second homology class A is said to be effective if there exists an 
$\omega$-compatible
almost complex structure $J$ and a non-trivial $J$-holomorphic sphere $f :S^2\rightarrow V$ such
that $[f]=A.$ We note that in [FHV] and [HV], the condition (*) was never
 stated explicitly, although
such a restriction seems to be  necessary for their method of
using Hamiltonian functions as we remarked above.
We may view the main
results obtained in [FHV] and [HV] as a stabilized version of
the  Weinstein conjecture.
In this aspect, as we mentioned before,
we are able to solve such a stabilized
Weinstein conjecture completely
without any restriction on $(V, \omega)$. (See Theorem \ref{third} above.)
For three dimensional contact manifolds, many
deep results have been proved on
the Weinstein conjecture and
related problems
by Eliashberg, Hofer, Zehnder and others (c.f.   [EH], [HWZ]).
For example, Hofer solved the Weinstein conjecture for overtwisted
contact 3-manifolds.

The main focus of this paper, however, is not only to prove the
stabilized Weinstein conjecture, but to establish
the full relationship between the existence of
$J$-holomorphic curves of any genus and the existence of non-trivial closed orbits of 
${\widetilde H}$.
Such a relationship 
obtained by using Gromov-Witten invariants of any genus  did
not appear in previous literature even for the semi-positive case.

The general idea of proving the existence of closed orbits for
the Hamiltonian equation (1) by using genus zero $J$-holomorphic curves or
perturbed $J$-holomorphic curves was already realized by Floer, 
Hofer and Viterbo. One quantitative form of such an idea was developed
in [FHV] and [HV] as their theory of d-index. The desired
existence results were then obtained by exploiting the 
deformation invariance of the d-index.  
However, the results obtained by this theory are quite limited.
It may be partly because the well-known difficulty of the 
transversality of multiple covered $J$-holomorphic spheres of negative
first Chern class and partly because the pathological nature of Hamiltonian function
used in d-index. 
In fact, most of the results obtained in [HV]  were not proved  even
for semi-positive symplectic manifolds before.

The recent progress on the Floer homology theory and GW-invariants
(cf. [FO], [LT], [LiuT1]) enables us to
overcome the difficulty of transversality.
Furthermore, in this paper, we will describe the full relationship of the 
existence of
 $J$-holomorphic curves of any genus 
with the existence of closed orbits of the Hamiltonian equation (1)
in its general form. We believe that this new finding will throw light
on solving the Weinstein conjecture completely. Hopefully, this new 
finding also gives clues to understanding  the mystery of Gromov-Witten invariants
on symplectic manifolds.

To prove Theorem \ref{first}, we choose an $\omega$-compatible almost complex
structure $J $ 
and consider
a family of Hamiltonian functions
${\widetilde H}_\lambda=\lambda\cdot {\widetilde H},$ $\lambda\in[0,\infty).$ 
Let $A\in H_2(V)$
and $\alpha_+,$  $\alpha_-,$ 
$ \beta_j\in H_*(V)$, $ j=1,\cdots, n.$
As before, we assume that $\alpha_{-}\in V_{+}$ and $\alpha_{+}\in V_{-}.$
The key step now is to obtain a Morse function $H$ which is a
small perturbation  of ${\tilde H}$ such that $H_{\lambda}=\lambda \cdot H$ has no non-trivial
closed orbits in $V\setminus W$ and has the  same  closed orbits as ${\tilde H}_{\lambda}$
 has
in $W$, if $0<\lambda< 1+\omega(A)/2\epsilon.$  The reason for choosing the quantity
$1+\omega(A)/2\epsilon$ will be explained in Section \ref{gw} of
this paper (see [HV] also).
By the usual Morse theory, $\alpha_{+}$ $(\alpha_{-})$ can be represented by
a linear combination of some critical points of $H$, denoted it by $c_{+}$ 
$(c_{-}),$ together with the associated
descending (ascending) manifold $M(c_{+})$ $(M(c_{-}))$ in $V_{-}$ $(V_{+}).$

We  will define a perturbed GW-invariant
$\Phi_{A,J_\lambda,H_\lambda, g}(c_-,c_{+},\beta_{1}, \cdots,\beta_{n})$ in 
Section \ref{gw} and Section \ref{gw1}. In genus zero case, it counts
algebraically the  $\nu$-perturbed $(J,H_{\lambda})$-maps $$u: ({\bf R}^1 \times S^1; x_{1},
\cdots, x_{n})  \rightarrow (V; \beta_{1}, \cdots, \beta_{n}),$$
satisfying the conditions 

\noindent  (i) ${\bar{\partial}}_{J, F_{\lambda},\nu} u=0,$
(ii) $\lim _{s \rightarrow {+} \infty}u(s,t) =c_{+},$
$\lim _{s\rightarrow {-} \infty}u(s,t )=c_{-},$ (iii) $[u]=A$ 
(see Section \ref{gw1} for higher genus case).
  
\begin{theorem}\label{fifth}

When $\lambda$ is small enough, 
 $$\Phi_{A, J_\lambda, H_\lambda,g}(c_-,c_{+},
\beta_{1}, \cdots,\beta_{n})=
 \Psi_{A,g}(\alpha_-,
\alpha_{+},\beta_{1}, \cdots, \beta_n).$$
\end{theorem}

This theorem were claimed in [PSS] and [RT] for semi-positive case.
The  third different method were described in [L2].
However, these methods are not sufficient for general symplectic manifolds.
Using the techniques devoloped in [LiuT1], we will 
  prove  this theorem  in  [LiuT3].

 The following  theorem is the main technique part of this paper.  
  
\begin{theorem}\label{six} 

If $H$ has no non-trivial closed orbits, then the perturbed $GW$-invariant
$\Phi_{A,J_\lambda,H_\lambda,g}(c_-, c_{+}, \beta_{1},\cdots, \beta_n)$
is well-defined, independent of the choice of
$\lambda\in (0,1+\omega (A)/2\epsilon).$ Moreover,
$\Phi_{A,J_\lambda,H_\lambda,g}(c_-,c{+},\beta_{1},\cdots,\beta_n)=0$ when
$\lambda> \omega (A)/2\epsilon+{1\over 2}.$
\end{theorem}

 The proof of this theorem requires a $T^{N_P}$-equivariant version of
the new technique developed in [LT] and [LiuT1]. 
The simplest case of  such a
theory ,  the $S^1$-equivariant case,  was already used in our computation 
of  Floer homology in [LiuT1]. 

Now Theorem \ref{first} follows from Theorem \ref{fifth} and Theorem 
\ref{six} easily.

The  main body of this paper  ( from Section 3 to Section \ref{gw} ) 
is devoted to establish
 a Morse theoretic version of $GW$ invariants  of genus zero case for
 general symplectic 
 manifolds under the assumption   that $H$ has no non-trivial closed orbits. 
 We then prove Theorem \ref{first} for genus zero case  in Section \ref{gw}. In the last section,
 we generalize the theory of genus zero case to higher genus
 case and prove Theorem 1.1 for any genus. 

 We note that if  the symplectic manifold $V$ is semi-positive, 
 our theory  in this paper can be developed  in a  much simpler
 manner. This includes the case where dimension of $V$ is four or six.

 This paper is the detailed version of our announcement [LiuT2].
 During the preparation of this paper, we learned that W. Chen proved
 some relevant results for the $4$-dimensional case in [C] 
 by a different method.

 From now on until the end of Section \ref{gw}, we will only deal with the
 case of genus zero. The discussions for higher genus cases are identical and  will be outlined in last section.  

\section{Compactness}
In this section, we will set up our main assumption, which will be used 
throughout the rest of this paper. We then explore the two simple consequences of the
assumption, the existence of the Morse function $H$ mentioned in Section 1 
and the compactness of the moduli spaces of cuspidal $(J,H)$-maps.

\medskip

\noindent $\bullet$ {\bf Main Assumption (I):}
${\tilde H}$ has no non-trivial closed orbits.

\medskip

The first consequence of this assumption is the following lemma,
 \begin{lemma}
There exists a Morse function $H$ such that 
(i) $H$ has same level sets as ${\tilde H}$ has in $W$;

(ii) $H$ is $C^0$-close to ${\tilde H}$ so that for any critical 
points $c_-$ in $V_-$ and $c_+$ in $V_+$,
$$ 0< \frac {\omega (A)}{ H(c_{-})-H_(c_{+})} < \lambda_{0},$$
where $\lambda_{0}= \frac {1} {2} + \frac {\omega (A)}{ 2\epsilon}$
and $A$ is an effective second homology class in the sense that it can be represented by some $J$-holomorphic sphere;

(iii) for  any  $0<\lambda < \lambda_{0}+\frac{1}{2},$ $H_{\lambda}=\lambda\cdot H$
has no non-trivial closed orbits of period one.

\end{lemma}
\proof

Choose  $r>0$ such that $$ 0< \frac {\omega (A)}{2\epsilon -2r -4 \delta}<
\lambda_{0}.$$
Here $\delta$ is the same as the one appeared in the definition of ${\tilde H}.$
As before we assume that  $\delta << \epsilon .$
Recall that ${\tilde H }(x)= t   $ if $x\in S_{t},$  $-\epsilon+2\delta  < t< \epsilon 
-2\delta .$  Let $V_{s}$ be the manifold $V_{+}\setminus S\times [0, s]$ with boundary
$S_{s}.$ Set ${\tilde H}_{+} ={\tilde H}|_{V_{\epsilon -2\delta}}$. 
Then $\nabla {\tilde H}_{+}\not=
0$ along  the boundary of $ V_{\epsilon-2\delta}.$ 

It is well-known that there is a $C^{2}$-small $C^{\infty}$-function   
${\tilde G}_{+}: V_{\epsilon-2\delta}\rightarrow {\bf R}$ such that  
${\tilde F}_{+} = {\tilde H}_{+}+{\tilde G}_{+}$ is a Morse function on  
$V_{\epsilon-2\delta}.$  Since ${\tilde H}_{+}$ is  regular along  
$S_{\epsilon-2\delta },$ we may arrange that ${\tilde G}_{+}$ vanishes  
near $S_{\epsilon-2\delta}.$ Now  decompose ${\tilde F}_{+}$ 
as: ${\tilde F}_{+}= \epsilon -2 \delta + {\bar F}_{+}$ and define 
${\tilde F}^{\lambda}_+ =\epsilon -2\delta + \lambda {\bar F}_+.$ 
Let ${\bar F}_+^{\lambda}=\lambda {\bar F}_+.$  Then $\nabla 
{\tilde  F}_+^{\lambda}=
\nabla{\bar F}_+^{\lambda}.$ Hence  ${\tilde F}^{\lambda}_{+}$ has  
non-trivial closed orbits if and only if ${\bar F}^{\lambda}_{+}$ has. 
Now $\| {\bar F}^{\lambda}_{+}\|_{C^2}= \lambda \|{\bar F}_{+}\|_{C^2},$ 
which implies that the $C^2$-norm  of  ${\bar F}^{\lambda}_{+}$ is small 
when $\lambda$ is small enough. Therefore, there exists a $\lambda_1>0$ 
such that ${\overline F}^\lambda_+$ ( hence ${\tilde F}^\lambda_+$) has 
no non-trivial closed orbits of period $1$ for 
$0<\lambda<\lambda_1.$

Similarly we can also define ${\tilde F}_-, {\tilde F}^\lambda_-$, etc on 
$V_{-\epsilon+2\delta}.$

Fix a $\lambda>0$ satisfying the following two conditions:

\noindent (a) $\lambda(\lambda_0+1/2)<\lambda_1;$

\noindent (b) for any critical point $c_+$ of ${\overline F}_+$ and 
$c_-$ of ${\overline F}_-$,
$$\lambda | {\overline F}_+(c_+)-{\overline F}_-(c_-)|<2r.$$

Define $H_+={\tilde F}_+^\lambda$ and $H_-={\tilde F}_-^\lambda.$ 
We will extend $H_+\cup H_-$ to $V$ to get an $H$ with the same 
level set as ${\tilde H}$ has in the ``middle part" $S\times 
[-\epsilon+2\delta, \epsilon-2\delta].$ If this is done, 
then it follows from (a) 
and (b) that $H$ has the required properties of the lemma.

We define $H$ on  $S \times (-\epsilon+2\delta, \epsilon-2\delta)$ 
to extend $H_+\cup H_-$ as follows,
\[
H(x)=
\left\{
\begin{array}{cl}
{\tilde H}(x), & x\in S\times (-\epsilon+3\delta, \epsilon-3\delta)\\
\psi (t(x)), & x\in S\times \{(-\epsilon+2\delta, -\epsilon+3\delta)\cup (\epsilon-3\delta, \epsilon-2\delta)\}.
\end{array}
\right.
\]

Here $t(x)$ is the  $t$-coordinate of  $x$ and 
$ \psi: (-\epsilon+2\delta , -\epsilon +3\delta) \cup 
(\epsilon-3\delta, \epsilon -2\delta) \rightarrow {\bf R}$ 
is an increasing $C^{\infty}$- function defined by  requiring that

\noindent (i) on $(\epsilon-3\delta, \epsilon-2\delta)$, $\psi$ 
connects smoothly the two functions $\psi_1(t)=t, t\leq\epsilon-3\delta,$ 
and $\psi_2(t)=(\epsilon-2\delta)+\lambda(t-(\epsilon-2\delta)), t\geq 
\epsilon-2\delta;$

\noindent (ii) $\psi$ does similar thing on $(-\epsilon+2\delta,
-\epsilon+3\delta).$

Clearly $H$ so defined has the same level sets as ${\tilde H}$ has in 
$S\times (-\epsilon+2\delta, \epsilon-2\delta).$

\QED

Our assumption now becomes

\medskip

\noindent $\bullet$ {\bf Main Assumption (II)}:

$H_\lambda=\lambda\cdot H$ has no nontrivial
closed orbits of period one, for $0\leq \lambda\leq \lambda_0+{1\over 2}.$
\medskip

Later on we will make some $C^\infty$-small generic perturbation of $H$. Since
the perturbation can be made arbitrarily small, we will assume that the main
assumption (II) also holds for those perturbed $H$.

We now state the consequence of the assumption on the compactness of the
moduli space of cuspidal $(J_\lambda, H_\lambda)$-maps, where $0<\lambda<\lambda_0+1/2.$

\begin{lemma}

Fix any two critical points $c_-$ and $c_+$ of $H_\lambda$, let $\{f_i\}$
be a sequence of $(J_{\lambda_i}, H_{\lambda_i})$-maps of class $A$ connecting $c_-$ and
$c_+$, with $\lambda_i\in [\epsilon, \lambda_0+1/2]$ for some small
$\epsilon >0.$ After reparametrization of the domain 
of $f_i$  and taking subsequence, we have that 
$\{f_i\}$ weakly $C^\infty$-converges to a cuspidal 
$(J_{\lambda_\infty}, H_{\lambda_\infty})$-map
$f_\infty$ of same class $A$ connecting $c_-$ and $c_+.$

\end{lemma}

The domain $\Sigma$ of a cuspidal $(J_\lambda, H_\lambda)$-map $f$ is a union
$\Sigma=\cup_{i=1}^{N_P} P_i\cup_{j=1}^{N_B} B_j$ of its principal components
$P_i$ and bubble components $B_j$. Each $P_i\cong {\bf R}^1\times S^1$ and
the collection $\{P_i\}$ form a chain. Each bubble component $B_j\cong S^2$ is
attached to some $P_i$ or some other $B_k$ at some of its singular points. All
components of $\Sigma$ form a tree. 

\begin{definition}

A continuous map $f:\Sigma\rightarrow V$ is said to be cuspidal $(J_\lambda,
H_\lambda)$-map of class $A$ connecting $c_-$ and $c_+$, if there exist
$N_P+1$ critical points $c_1,\cdots, c_{N_P+1}$ of $H_\lambda, $ with 
$c_1=c_-, $ $c_{N_P+1}=c_+$ such that

\noindent (i) ${\bar\ptl}_{J_\lambda, H_\lambda} f^P_i=0, $ $\lim_{s\rightarrow
-\infty} f^P_i(s,\theta)=c_i, $ $\lim_{s\rightarrow\infty} f^P_i(s,\theta)=c_{i+1}$, 
$i=1,\cdots, N_P.$

\noindent (ii) ${\bar\ptl}_{J_\lambda} f^B_j=0$.

\noindent (iii) $\sum_i[f^P_i]+\sum_j [f^B_j]=A.$

\end{definition}
An analogy of this lemma, in which the Morse function $H_\lambda$ is replaced
by some generic time-dependent Hamiltonian function is proved in 
[F] section 3.
The proof there can be easily adapted to our case as long as we can make 
sure that (i) $H_\lambda$ has no non-trivial periodical orbits of period one for
$0<\lambda <\lambda_0+1/2;$ (ii) all critical points $c$ of $H_\lambda$, 
when considered as a trivial periodical orbit of the time-independent
Hamiltonian function $H_\lambda$ is non-degenerate in the sense of Floer homology.
Now (i) follows from our main assumption and (ii) can be achieved by a
small $C^\infty$-perturbation of $H$. Note that (ii) implies that any
$(J, H)$-map convergent  to $c$ along its ends will converge to $c$
exponentially.

\section{Moduli Space of Stable Maps}

In this section, we will define the various moduli spaces of stable maps
needed to define the Morse theoretical version of GW-invariant.

\subsection{ Stable Curves}

Stable curves will appear as the domains of stable maps, which are to be
defined below. From this Section up to Section \ref{gw} we will only consider semi-stable
  ( connected ) curves of genus zero. Geometrically such a curve $\Sigma$ is a union
of its components $\Sigma_l\cong S^2$ with only double points as its singularities,
and its components form a tree ($H_1(\Sigma)=0$). 

We now define semi-stable ${\cal F}$-curves and ${\cal G}$-curves:

\begin{definition}
An n-pointed semi-stable ${\cal F}$-curves $(\Sigma, l, x)$
is a semi-stable curve $\Sigma$ with $n$ (ordered)
marked points $x=\{ x_1,\cdots, x_n\}$ in $\Sigma$ away from its singular points
such that the components of $\Sigma$ can be divided into principal components
$P_i, i=1,\cdots, N_P$, and bubble components $B_j, j=1,\cdots, N_B.$ The
principal components form a chain in such a way that each $P_i$ has two
distinguished points $z_i$ and $z_{i+1}$, $i=1, \cdots, N_P$ such
that $P_i$ and $P_{i+1}$ join together at $z_{i+1}.$ 
$l$ is the collection of marked lines $l_i$ on $P_i$ connecting its ``ends" 
$z_i$ and $z_{i+1}$. 

\end{definition}

Using the marked line $l_i$ , we may identify $(P_i\setminus\{z_i, z_{i+1}\}; l_i)$ with 
$({\bf R}\times S^1; \{\theta =0\}).$

An n-pointed semi-stable ${\cal G}$-curve $(\Sigma, x)$ can be
obtained from the corresponding ${\cal F}$-curve by simply forgetting all
marked lines $l_i$'s.


Two semi-stable ${\cal F}$-curves $(\Sigma^1; l^1, x^1)$ 
and $(\Sigma^2; l^2,x^2)$ are said to be equivalent if there
is a homomorphism $\phi:\Sigma_1\rightarrow \Sigma_2$ which preserves 
marked points and lines such that the restriction of $\phi$ to any component
of $\Sigma_1$ is a biholomorphic map. We will use $\la \Sigma, l, x\ra $ to denote
the resulting equivalence class of $(\Sigma, l, x).$ Similarly we can define 
equivalence class for semi-stable ${\cal G}$-curves by simply forgetting those
marked lines in the definition of the equivalence of ${\cal F}$-curves, 
and we will use $\la \Sigma,x\ra $ 
to denote the equivalence class of a semi-stable ${\cal G}$-curve 
$(\Sigma, x).$

\begin{definition}

$${\cal FM}_{0,n}=\{ \la \Sigma, l, x\ra  | (\Sigma, l, x) \mbox{ is a semi-stable } {\cal F}\mbox{-curve}\},$$
$${\cal GM}_{0,n}=\{ \la \Sigma, x\ra  | (\Sigma, x) \mbox{ is a semi-stable }{\cal G}\mbox{-curve}\}.$$ 
\end{definition}

There is an obvious forgetting map:
$${\cal FM}_{0,n}\rightarrow {\cal GM}_{0,n}$$ sending 
$\la \Sigma,l, x\ra $ to $\la \Sigma, x\ra .$

From now on, for simplicity, we will call a semi-stable ${\cal F}$-curve or \
a semi-stable ${\cal G}$-curve an ${\cal F}$-curve or ${\cal G}$-curve respectively.

Given an ${\cal F}$-curve $(\Sigma, l, x)$ or a ${\cal G}$-curve $(\Sigma, x)$, there 
is an obvious way to
add minimal number of markings $y_i$ to  an unstable principal component $P_i$ 
and $y_j^k , 1\leq k \leq 2$, to an unstable bubble component $B_j$ to stabilize 
$\Sigma$. 
We will use  $y$ to denote the set of the added markings  and 
$(\Sigma, l, x; y)$ and $(\Sigma, x; y)$ to denote the resulting stabilized 
${\cal F}$-curve
and ${\cal G}$-curve and call them stable ${\cal F}$-curve and stable 
${\cal G}$-curve
respectively. Here stability means that each of its components contains at least
three singular points or marked points in  $x$ or $y$. There is an obvious
forgetting map here from the set of stable ${\cal F}$-curves or ${\cal G}$-curves
to the set of semi-stable ones, sending $(\Sigma, l, x; y)$ to $(\Sigma, l, x)$ or
$(\Sigma, x; y)$ to $(\Sigma, x)$ respectively.

From now on we will use various simplified notations, depending on the context, to denote above curves. For instance we may write $(\Sigma, y)$ for $(\Sigma, l, x; y),$ if no confusion arises.

\noindent $\bullet$ {\bf Local Deformation of $(\Sigma, l, x; y)$}

\medskip

Given a stable ${\cal F}$-curve $(\Sigma; l, x; y)$, 
with double points  
$$d^{P_i}_m\in P_i, m=1, \cdots, M^{P_i} \quad\mbox{ and }\quad 
d_l^{B_j}\in B_j, l=1,\cdots,L^{B_j},$$
let $\alpha_m^{P_i}$ and
$\alpha_l^{B_j}$ be the complex coordinates of the corresponding 
double points of 
$d'^{P_i}_m $ and  ${d'_l}^{B_j}$ of a nearby curve $\Sigma'$ of 
same topological
type.
  Let $\alpha=\{\alpha^{P_i}_m, \cdots, \alpha^{B_j}_l\}$,
$m\leq M^{P_i}-3+r^{P_{i}}, $ and $l\leq L^{B_j}-3+r^{B_{j}}$, where $r^{P_{i}}$and  $r^{B_j}$  are  the number of 
elements in $x$ and $y$ in $ P_{i}$ and  $ B_j$ respectively. 
Let $\theta$ be the collection of all angular
coordinates $\theta_i$ of the third from last double or marked  point 
 of the principal components $P'_i$.
Now $u=(\alpha, \theta)$ gives rise to the 
universal local coordinate of nearby stable \Fc.
 We will use 
$ \Sigma_u$ to denote the nearby curve with coordinate $u$. 

For each double 
point in $\Sigma'=\Sigma_u$, say, $d'_1\in P'_1 $ and $d_2'\in B_2'$ with
$d_1'=d_2'$ in $\Sigma'$, we associate a complex gluing parameter
$t_1=t_2\in D_{\delta}=\{ z\, |\, |z|<\delta\,\}.$ The corresponding
gluing here is the following: cut off the two discs of radius $|t_1|=|t_2|$
of $P_1'$ and $B_2'$ centered at $d_1'$ and $d_2'$ respectively and glue
them back along the boundary circles through a rotation of
$\arg\,\, t_i.$ Let $t$ be the collection of all such gluing parameters. 
Similarly
for each $z_i, i=2,\cdots, N_P-1,$ we associate a gluing parameter 
$\tau_i\in I_\delta=\{ r\,|\, r\in {\bf R}^{+}, \, r<\delta\,\}.$
There is also a similar but simpler gluing process for each $\tau_i$. 
Let $\tau=(\tau_i)$ and $v=(t,\tau).$ Then $(\Sigma_{(u,v)}
=(\Sigma_{(\alpha, \theta, t,\tau)}, l))$ is the local ``universal "
deformation of $(\Sigma, l)$ as an \Fc. 

The ``universal" deformation for
\Gc  can be defined similarly. Since in this case there is no such 
marked lines $l$ appearing, there is no such parameter $\theta$ and
associate to each $z_i$ is a complex parameter ${\tilde{\tau_i}}
=(\tau_i,\theta_i)$ instead of $\tau_i.$ Let ${\tilde t}$ be the
collection of all complex gluing parameters associated with double points and
``ends" of $\Sigma$, the $\Sigma_{(\alpha, {\tilde t})}$ is the
``universal" deformation of $\Sigma$ as a \Gc.

\medskip

\noindent $\bullet$ {\bf Fixed Markings}

\medskip

Recall that in order to define $GW$-invariants, one needs to specify a cycle
$C$ in $H_{*}({\bar {\cal M}}_{g, n+2}.$
Chosing such a cycle will impose restrictions to the possible domains in the bubbling process of the Gromov-Floer compactification of stable maps. For simplicity, we only describe in detail the case $C=\{pt\}.$  The general case can be treated similarly.

Since the main issue here only involves how to fix
marked points, we can treat both
{\Fc}s and {\Gc}s equally
We only formulate the ``fixed marking" process for {\Gc}s.

Let $(S^2; -\infty, +\infty; {\tilde x}_1,\cdots, {\tilde x}_n) $ be a 
fixed a model, where
$-\infty$ and $+\infty$ are the two ``ends" if we identify $S^2\setminus \{-\infty,
+\infty\}$ with ${\bf R}^1\times S^1.$ 
We want to define the notion of a semi-stable curve with
``fixed" marked points. 
markings, $x_1, \cdots, x_n$ ( modeled on $(S^2; {\tilde x})$ ) if 
(i) there exists a principal component $P_i$ and $n$ many of its double
points or marked points, $d_1, \cdots, d_n$ such that 
$(P_i; z_i, z_{i+1}; d_1, \cdots,  d_n)\cong 
(S^2; -\infty, +\infty; {\tilde x}_1,\cdots, {\tilde x}_n); $ (ii) each
marked point $x_i$ of $(\Sigma, x)$ lies on the branch $B(d_i)$ 
consisting of all bubble components with `` root" $d_i$, if $x_i\not = d_i$. 


Let $(\Sigma,x,y)$ be the minimal stabilization of $(\Sigma, x)$, the
next lemma explains why the above two conditions are the  desired ones.

\begin{lemma} \label {fixmarking}

There is a gluing procedure such that for any gluing parameter 
${\tilde t}$ with non zero components,  
$$(\Sigma_{(\alpha, {\tilde t})}; z_1,  z_{N_P+1}, x_1, \cdots, x_n)
\cong (S^2; -\infty, +\infty; {\tilde x}_1,\cdots, {\tilde x}_n),$$
after forgetting those markings $y$ of $(\Sigma_{(\alpha, {\tilde t})}, x, y).$

\end{lemma}

Note here the parameter $\alpha$ is subject to the restriction imposed by
the ``fixed" marking condition.

\proof

Away from those components $ B_{i,j}$ in $B(d_i)$, 
the gluing procedure is the same as
before. Let $x_{i}\in B_{i,l}.$ Because of the tree structure of the 
components of ${\Sigma}$, there is a unique chain of bubble components
$B_{i, j}, j=1, \cdots, l$ of $B(d_i)$ with each $B_{i,j}$ having two particular double
points $d^{j}$ and $d^{j+1}$, connecting $d_i$ and $x_i$. Here $d_i=d^1$ and
$x_i=d^{l+1}$. Now there is a unique identification of each 
$B_{i,j}-\{d^j, d^{j+1}\}\cong {\bf R}^1\times S^1 $ up to translations and 
rotations of ${\bf R}^1\times S^1$. Use the cylindrical coordinate here ( or
the corresponding polar coordinate ) to do the gluing associated with the 
double points $d^j.$ It is easy to see that after forgetting markings other
than $x$, $\Sigma_{\alpha, {\tilde t}}$ has the desired property. 

\QED

\subsection{Stable Maps}

$\bullet$ Given a homology class $A\in H_2(V, {\bf Z})$, a stable $(J, H)$-map
$f$ from an \Fc $(\Sigma, l)$ to $V$ of class $A$ connecting critical point $c_{-}$
and $c_{+}$ of $H$
is a map defined on $(\Sigma, l)\setminus \bigcup^{N_P+1}_{i=1}\{z_i\}$ such
that:

\noindent (i) on each principal component $P_i$, 
$$\frac{\ptl f_i^P}{\ptl s}+J(f^P_i)\frac{\ptl f^P_i}{\ptl\theta}
-\nabla H(f^P_i)=0,$$
where $(s, \theta)\in {\bf R}^1\times S^1$ is the cylindrical coordinate of
$P_i$ and $f^P_i=f|_{P_i-\{z_i, z_{i+1}\}};$

\noindent (ii) there exist $c_i, c=1, \cdots, N_P+1$ with $c_1=c_{-},
c_{N_P+1}=c_{+}$ such that 
$$\lim _{s\rightarrow -\infty} f_i^P(s, \theta)=c_i\,\, \mbox{and }\,
\lim_{s\rightarrow +\infty} f^P_i(s, \theta)=c_{i+1};$$

\noindent (iii) on each bubble component $B_j$, ${\bar{\ptl}}_{J} f^B_j=0$;

\noindent (iv) $\sum_{i} [ f_i^P]+\sum_j [f^B_j]=A$;

\noindent (v) each constant component is stable in the sense
that it has at least three double or marked points.

\medskip

Two such maps $f_1$ and $f_2$ with ${\cal F}$-curves as their domains are said to be 
equivalent if there is an identification 
$$\phi :(\Sigma^1, l^1, x^1)\rightarrow (\Sigma^2, l^2, x^2)$$ 
such that $f_2=f_1\circ\phi.$ Similarly, we can define equivalent relation
for stable $(J, H)$-maps with {\Gc}s as domains. We will use $\la f\ra $ to denote
the resulting equivalence class of $f$. 

We also need the notion of stable $L_k^p$-maps, which can be defined by 
simply requiring that $f$ is a $L^p_k$-map, $k-\frac{2}{p}>1$ satisfying
requirements (ii) , (iv) and (v).

Each stable map $f$ determines an intersection pattern $D_f$ which encodes
the following information:

\noindent (i) the topological type of the domain $\Sigma=\Sigma_f$;

\noindent (ii) the homology classes $[f_i^P], [f^B_j]\in H_2(V, {\bf Z})$;

\noindent (iii) the critical points $c_i, i=1, \cdots N_P+1$.

Note that the topological type of $\Sigma$ is determined by its intersection pattern 
$I=I_\Sigma$, which can be thought as a pairwise correspondence of
the double points of $\Sigma$ lifted to the smooth resolution of $\Sigma.$

Given a stable map $f$, we define its energy
$$E(f)=\sum_i E(f^P_i)+\sum_j\int _{S^2}(f^B_j)^*\omega,$$ where
$E(f^P_i)=\int\int_{{\tiny{\bf R}}^1\times S^1}|\frac{\ptl f^P_i}{\ptl s}|^2.$

Note that if $f$ is a $(J,H)$-map of class $A$ connecting $c_{-}$ and 
$c_{+}$, then
$$E(f)=\omega (A)+f(c_{+})-f(c_{-}).$$

\begin{lemma}\label{lowbound}

For a generic choice of $(J, H)$, there exists a $\delta=\delta(J,H)>0$,
such that for any non-constant stable $(J,H)$-map $f$, $E(f)>\delta.$

\end{lemma}
\proof

Assume that there exists a sequence of $(J,H)$-maps $\{f_i\}$ of class $A$
such that $\lim E(f_i)=0$. By  choosing a suitable subsequence we may assume 
that each $f_i$ has only one principal component and connects two fixed
critical points $c_{-}$ and $c_{+}$. It follows from Gromov-Floer compactness
theorem for cuspidal maps that a subsequence of $\{f_i\}$, still denoted
by $\{f_i\}$, is $C^0-$convergent to a constant map. Hence $c_-=c_+$, $[f_i]=0$,  
for large $i$. If $f_i$ is not a constant, there exists a unique simple
$(J, \frac{1}{m} H)$-map ${\tilde f}_i$ such that $f_i={\tilde f}_i\circ\pi_m$,
where $\pi_m : {\bf R}^1\times S^1\rightarrow {\bf R}^1\times S^1$ is given
by $\pi_m(s,\theta)=(ms, m\theta).$ Here ${\tilde f}_i$ being simple means
that it can not be factorized through further for any $m>1.$ 

Consider the moduli space 
\begin{eqnarray*}
& & {\cal M}^0(c_-, c_+; J, H, A)  \\
& & =\{g\,|\, g:{\bf R}^1\times S^1\rightarrow V \mbox{ is a $(J, H)$-map}, 
[g]=A, g\mbox{ is simple}\,\}. 
\end{eqnarray*}
\noindent Then  for a generic choice of $(J, H)$, 
$$\dim {\cal M}^0(c_-, c_+; J, \frac{1}{m}H, \frac{1}{m}A)=Ind (c_+)
-Ind (c_-)+2c_1(A)/m, $$ which is zero in the case that $c_-=c_+$ and $A=0. $
Clearly $${\tilde f}_i\in {\cal M}^0(c_-, c_+; J, \frac{1}{m}H, 0).$$ 
Since ${\tilde f}_i$ is not a constant, it follows from [FHS] that for a
generic choice of $(J, H)$, ${\tilde f}_i$ has a two dimensional symmetries, 
which implies that $$\dim {\cal M}^0(c_-, c_+; J, \frac{1}{m}H, 0) \geq 2.$$
This is a contradiction.

\QED

An intersection pattern $D$ is said to be effective if $D=D_f$ with $f$
being a stable $(J, H)$-map. Let $e>0$ and define 
$${\cal D}^e=\{ D\,|\, D\mbox{ is  effective, } E(D)\leq e\}, $$
where the energy $E(D)=E(D_f)=E(f).$

\begin{lemma}\label{fip}

${\cal D}^e$ is finite for any $e>0$.

\end{lemma}

\proof

It follows form Gromov-Floer compactness theorem for cuspidal maps that there are
at most finitely  many possible homology classes which can be represented
by some $(J, H)$-map $f$ with $D_f\in {\cal D}^e$. Therefore it is
sufficient to prove that there are only finitely many possible topological
types of $\Sigma_f$ for such  $f$. To this end, we observe that $f$ has
at most $[\frac{e}{\delta}]+1$ non-trivial component. This implies that the
stabilized curve $(\Sigma,y)$ obtained by adding minimal number of markings 
to $ \Sigma$ has at most $2([\frac{e}{\delta}]+1)+n$ markings. This in
turn bounds the number of double points, and hence bounds the number of
components of $\Sigma_f.$

\QED

There is a partial order relation in ${\cal D}^e$ defined as follows
: $D_1=D_{f_1}\leq D_2=D_{f_2}$ if (i) $ \Sigma_{f_2}$ can be
obtained from $\Sigma_{f_1}$ topologically by the gluing 
construction described in Section 2.1; (ii) the homological classes 
represented by the components of $f_1$ and $f_2$ are compatible with the gluing 
construction. ( See the next subsection for the definition of the gluing 
of stable maps.) 

\subsection{\bf Moduli Spaces of Stable Maps}

Now we can define various moduli spaces of stable maps. 

Let ${\cal F}{\cal M}(c_-, c_+; J, H, A) $ be the moduli space
of equivalence classes of stable $(J,H)$-maps   of class $A$ 
connecting $c_-$ and $c_+$ with ${\cal F}$-curves as domains. 

Similarly we can define the moduli space ${\cal GM}(c_-, c_+; J, H, A) $  of
the equivalence classes of stable $(J, H)$-maps of class $A$  
connecting $c_-$ and $c_+$ with ${\cal G}$-curves
as domains.     

Let ${\cal F}{\cal B}^e(c_-, c_+; A)$ be the moduli space of equivalence 
classes of stable
$L^p_k$-maps  of class $A$ connecting $c_-$ and $c_+$ with ${\cal F}$-curves 
as domains, the energy of whose elements is less than $e$.


Since the energy $E(f)$ is bounded for any element in 
${\cal F}{\cal M}(c_-, c_+; J, H, A), $  
${\cal F}{\cal M}(c_-, c_+; J, H, A) \subset 
{\cal F}{\cal B}^e(c_-, c_+; A)$ when $e$ is large enough. We will
choose such an $e$ once for all and omit the superscript $e$ for the 
moduli space of $L^p_k$-maps. 

Similarly we define ${\cal G}{\cal B}(c_-, c_+; A)$. 

We can also restrict to some particular intersection pattern $D\in{\cal D}^e$ 
and define the corresponding moduli spaces. We denote them by 
$${\cal F}{\cal M}^D(c_-, c_+; J, H, A) \quad \mbox{ and }  
\quad {\cal G}{\cal M}^D(c_-, c_+; J, H, A) \quad \mbox{ etc.}$$ 
From now on,  we will omit $c_-$ and $c_+$ in our notations of above
moduli spaces when no confusion arises.


\medskip

\noindent $\bullet$ {\bf Weak topology on ${\cal FM}(J, H, A)$ and 
${\cal GM}(J, H, A)$  }

\medskip

There are two different but equivalent topology on the moduli spaces of 
stable $(J, H)$-maps, the weak $C^{\infty}$-topology and  strong $L^p_k$-
topology.

We start with defining the weak $C^{\infty}$-topology. We will only deal with
${\cal FM}$ and leave the corresponding statements for ${\cal GM}$ to readers. 

\noindent$\bullet$ {Definition of Weakly Convergence}

Given a sequence $\{\la f_i\ra \}_{i=1}^{\infty}$ of equivalence classes of 
stable $(J, H)$-maps with ${\cal F}$-curves as domains, we say that $\{\la f_i\ra \}$ 
is weakly $C^{\infty}$-convergent to a stable $(J, H)$-map $\la f_{\infty}\ra $ if
there are $f_i\in \la f_i\ra , f_{\infty}\in \la f_{\infty}\ra $ such that 
the following conditions hold.

\noindent (i) After
stabilized by adding minimal number of markings, the stabilized domains
$\Sigma_i=\Sigma_{f_i}$ is  convergent to $\Sigma_\infty=\Sigma_{f_\infty}$
in the sense that when $i$ is large enough, there exist identifications
of stable ${\cal F}$-curves, $\phi_i :\Sigma_{(u_i, v_i)}\rightarrow \Sigma_i$ and
$\phi_\infty :\Sigma_{(0,0)}\rightarrow \Sigma_\infty$, such that $(u_i, v_i)
\rightarrow (0, 0)$ as $i\rightarrow \infty.$ Here
 $\Sigma_{(u_i, v_i)}$ is the 
local deformation of $\Sigma_{(0,0)}$ defined before in Section 2.1.

\noindent (ii) Given any compact subset $K\subset \Sigma_{(0,0)}\setminus \{\mbox
{singular points}\}$, there is an obvious embedding 
${\imath}^K_i: K
\rightarrow\Sigma_{(u_i,v_i)}$ through the gluing construction,
 when $i$ is large.
Define $f^K_i=f_i\circ \phi_i\circ {\imath}^K_i: K\rightarrow V$ and 
$f_{\infty}^K
=(f_\infty\circ\phi_\infty)|_K$. We require that $\{f_i^K\}_{i=1}
^\infty$ is $C^{\infty}$-convergent to $f_\infty^K$ for any $K$ as above.

\noindent (iii) $\lim _i E(f_i)=E(f_\infty).$ 

We will call the induced topology on ${\cal FM}(J, H, A)$ and ${\cal GM}(J, H, A)$
the weak $C^\infty$-topology.

\begin{theorem}

${\cal FM}(J, H, A)$ and ${\cal GM}(J, H, A)$ are compact and Hausdorff
with respect to the weak $C^\infty$-topology.

\end{theorem}

The compactness part of this theorem for cuspidal maps is known as
Gromov-Floer compactness theorem. The analysis  there can be adapted here to
prove the corresponding part of our theorem up to some suitable modification.
The Hausdorffness is not true for the moduli space of cuspidal maps, but only
holds for the moduli space of stable $(J, H)$-maps. The complete proof of this
statement is in [LiuT1], Sec 4. We refer our readers to the proof there.

To define the strong $L^p_k$-topology we mentioned before, we need to work
with stable $L^p_k$-maps.

\noindent$\bullet$ {\bf Strong $L^p_k$-topology and local uniformizer}

\medskip

\noindent$\bullet$ Local deformation of stable $(J, H)$-maps.

We start with defining the local deformation of a stable $(J, H)$-map. Again
we only deal with stable maps with ${\cal F}$-curves as domains. Given a stable map
$\la f\ra $, let $f\in \la f\ra $ be a representative with \Fc  $(\Sigma, l)$ as its
domain. Let $(\Sigma_{(u, v)}, l)$ be the local universal deformation. We
define $F_{(u,0)}:\Sigma_{(u, 0)}\rightarrow V$ and $f_{(u, v)}:\Sigma_{(u,v)}
\rightarrow V$ as follows.

Choose a  family of homomorphisms $\phi_{(u,0)}$ of $\Sigma_{(u,0)}$
to $\Sigma_{(0,0)}$ such that the  restriction of $\phi_{(u,0)}$ to each
component of $\Sigma_{(u,0)}$ is a diffeomorphism and it maps all double points
on the components of $\Sigma_{(u,0)}$ to the corresponding double points
of $\Sigma_{(0,0)}$. Moreover, $\phi_{(u,0)}$ is identity on each 
component of $\Sigma_{(u,0)}$ outside a prescribed small neighborhood of its
double points. When $|u|$ is small enough, such a  $\phi_{(u, 0)}$ exists. Note that
$\phi_{(u,0)}$ is not holomorphic. 

We define that $f_{(u,0)}=f\circ \phi_{(u,0)}$. 

Now $f_{(u,v)}$ is obtained from $f_{(u,0)}$ by the following gluing procedure 
with gluing parameter $v$.

It is sufficient to consider the following two simplest cases:

\noindent (i) $f_{(u,0)}=f_1\cup f_2$ with $f_1$ being a principal component and $f_2$
being a bubble component. Let $d_1=d_2$ be their double points, associated
with a complex gluing parameter $t.$ 

\noindent (ii) $f_{(u,0)}=f_1\cup f_2$ with both of them being principal components
jointed at their double point $z$ ( one of their ``ends"). Associate with $z$
a positive real gluing parameter $\tau$. 

The case of gluing two bubble components is the same as case (i) above and
the general case can be reduced to above cases. 

For case (i), let $D_1$ and $D_2$ be the small discs of $\Sigma_1$ and 
$\Sigma_2$ centered at $d_1$ and $d_2$ respectively. Let $(s_i, \theta_i),
i=1,2$ be their cylindrical coordinates given by $w_i=e^{-(s_i+i\theta_i)}.$
Then $\Sigma_{(u,t)}$ is obtained from $\Sigma_{(u,0)}=\Sigma_1\cup
\Sigma_2$ by cutting off $\{(s_i, \theta_i)\, | \, s_{i} > -\log |t|\}
\hookrightarrow D_i$ and gluing back to remaining part of $\Sigma$ 
along the boundaries through a rotation of $arg\, t.$ Choose a cut-off 
function 
\[
\beta(s)=\left\{
\begin{array}{cl}
1 & s<-\log |t|-2\\
0 & s> -\log |t|-1.
\end{array}\right.
\]

We define 

\[
f_{(u,t)}(w)=\left\{
\begin{array}{cl}
f_i(w), & w\in \Sigma_i\setminus \{(s_i, \theta_i)|s_i<-\log |t|-2\}\\
Exp_{f(d)}\beta(s_i)\cdot\xi_i(w), & w\in \{(s_i, \theta_i)|s_i>-\log |t|-2\},  
\end{array}\right.
\]
where $\xi_i(w)$ is defined by $f_i(w)=Exp_{f(d)}\xi_i(w)$ when
$|w|$ is small.

Case (ii) can be treated in a similar way. We leave it to the readers.

\noindent $\bullet$ Local uniformizer and $L^p_k$-topology of ${\cal FB}(A)$.

From now on, we will assume that $\dim V\geq 4.$ Let $(J, H)$ be a generic
pair.  Under the assumption , it is proved in [FHS]  

\begin{theorem}\label{simpt}

Given $f\in {\cal FM}^{D}(c_-, c_+, J, H, A)$, we have either

\noindent(i) $f$ is $\theta$-independent and hence a gradient 
line of $\nabla H$ or

\noindent(ii) there exists an integer $m> 1$, such that $f={\tilde  f}\circ \pi_m$, 
where $\pi_m:{\bf R}^1\times S^1\rightarrow {\bf R}^1\times S^1$ is  given by
$(s,\theta)\rightarrow (ms, m\theta)$ and ${\tilde  f}$ is simple in the 
sense that there is no further factorization.  Moreover if $f$ is already
simple, there exists at least one point $(s_0, \theta_0)\in {\bf R}^1\times
S^1$ such that $f(s,\theta)\not =f(s_0, \theta_0)$ if $(s, \theta)\not =
(s_0, \theta_0)$ and Rank$(df_{(s_0, \theta_0)})=2.$ We will call such point 
$(s_0, \theta_0)$ an injective point. Note that ${\tilde  f}\in
{\cal FM}^{D}(c_-, c_+, J, \frac{1}{m}H, \frac{1}{m}A)$. 
\end{theorem}

As pointed out in [FHS], if the theorem holds for $(J, H)$, so does it for
$(J, \frac{1}{m}H).$ 

Given $\la f\ra \in {\cal FM}(J, H, A)$, choose a representative $f\in \la f\ra .$ If an unstable principal 
component $f_i^P$ is $\theta$-dependent, it covers a simple map ${\tilde 
f}_i^P$. We may assume that there is an  injective point $(s_i, \theta_i)$ of
${\tilde f}^P_i$ lying on the middle circle $\{s_i=0\}.$ 

For simplicity, we may assume that $(s_i, \theta_i)=(0,0)\in l_i$  and use
$y_i=(s_i, \theta_i)$ to stabilize $f_i^P.$ For any unstable bubble component
$f_j^B$, it follows from [M] that, similar to the theorem above,
$f_j^B={\tilde  f}_j^B\circ\pi_j$ such that ${\tilde  f}^B_j $ has
only injective points away from finite points of $B_j$ and $\pi_j: B_j
\rightarrow S_2$ is  a finite branch covering. Choose $y_j^k, 1\leq k\leq 2$
to stabilize $B_j$ in such a way that $\pi_j(y_j^k)$ is an injective point.

Let ${\tilde {\bf H}}_i$ be the local hypersurface of codimension two at 
$f^P_i(y_i)$ such that $f_i^P$ is transversal to ${\tilde {\bf H}}_i$ at 
$y_i$ when $f_i^P$ is not
$\theta$-independent.
We then choose a local 
hypersurface
${\bf H}_i$ of codimension one such that ${\tilde {\bf H}}_i\hookrightarrow
{\bf H}_i$ and $f^P_i|_{l_i}$  is transversal to ${\bf H}_i$ at $y_i$. When
$f_i^P$ is $\theta$-independent, simply choose ${\bf H}_i$ of codimension 
one such that ${\bf H}_i$ is transversal to $f_i^P$ at $y_i=(0,0)$. Similarly
for each unstable bubble $f^B_j$, choose   hypersurface 
${\tilde{\bf H}}^k_j,
1\leq k\leq 2,$ of codimension 2 such that $f^B_j$ is transversal to ${\tilde{\bf H}}^k_j$ at
$y_j^k$. Let ${\bf H}$ be the collection of all those   hypersurfaces 
${\bf H}_i$'s and ${\tilde {\bf H}}^k_j$'s.
Consider the local deformation $f_{(u,v)}$ of $f$ with $\|(u,v)\|<\delta$
for some fixed small $\delta >0.$ Choose  an $\epsilon> 0$, we define a
local uniformizer of ${\cal FB}(A)$ near $\la f\ra \in {\cal FM}(J, H, A)$, 
\begin{eqnarray*}
& &{\cal F}{\tilde U}_{\epsilon}(f; {\bf H})  = \\
& &\{ g=g_{(u,v)}\, |\, \|g_{(u,v)}-f_{(u,v)}\|_{k,p}<\epsilon,
g_{(u,v)}(y_i)\in {\bf H}_i, \,\,  g_{(u,v)}(y^k_j)\in {\tilde {\bf H}}_j^k \},
\end{eqnarray*}
where $g_{(u,v)}:\Sigma_{(u,v)}\rightarrow V$ and $y_i, y_j^k\in 
\Sigma_{(u,v)}$ through gluing. Here the metric on $\Sigma_{(u,v)}$ to define
the $L^p_k$-norm is induced from that of $\Sigma_{(0,0)}$ through gluing. 

Before we state any  properties of 
${\cal F}{\tilde U}_{\epsilon}(f; {\bf H})$, we define the $L^p_k$-topology
on ${\cal FB}(A)$ and ${\cal GB}(A)$ by using a similar construction as above.
We only treat the case ${\cal FB}(A)$ as before. Given $\la f\ra \in {\cal FB}(A)$, 
choose a representative $f\in \la f\ra $ and define
$${\cal F}{\tilde U}_\epsilon (f)=\{ g=g_{(u,v)} \,| \|(u,v)\|<\epsilon,\,
\|g_{(u,v)}-f_{(u,v)}\|<\epsilon\}.$$ 
The process of forgetting those markings $y_i, y_j^k\in\Sigma_{(u,v)}$ induces
a natural projection map 
$$\pi_{\cal F}=\pi_{\cal F}(f):
{\cal F}{\tilde U}_\epsilon(f)\rightarrow {\cal F}U_\epsilon(f)=\pi_{\cal F}
({\cal F}{\tilde U}_{\epsilon}(f))\hookrightarrow{\cal FB}(A).$$

Let ${\cal F}U=\{ {\cal F}U_\epsilon(f)\,|\,f\in \la f\ra ,\, \la f\ra \in
 {\cal FB}(A)\}.$
One can directly check that 
\begin{lemma}\label{top}
${\cal F}U$ form a topological basis on ${\cal FB}(A)$. 
\end{lemma}

We will call the induced topology the (strong) $L^p_k$-topology on 
${\cal FB}(A)$. In particular we get an induced strong  $L^p_k$-topology
on ${\cal FM}(J, H, A)$ as a subspace of ${\cal FB}(A)$. It is proved in
[LiuT1], Section 4 

\begin{theorem}\label{eqn}

The two topologies on ${\cal FM}(J, H, A)$  (  ${\cal GM}(J, H, A)$ ) are
equivalent. In particular, ${\cal FM}(J, H, A)$  ( ${\cal GM}(J, H, A)$ )
is also compact with respect to $L^p_k$-topology. 

\end{theorem}

One of the corollary of this equivalence is 

\begin{cor}\label{Haus}

There \,\, exists \,\,an \,\,open \,\,neighborhood \,\,$W$ \,\,of 
${\cal FM}(J, H, A)$  \,\,( ${\cal GM}(J, H, A)$ )  in ${\cal FB}(A)$
 ( ${\cal GB}(A)$) such that  
$W$ is Hausdorff with respect to the $L^p_k$-topology.  
\end{cor}

We leave its proof to our readers as it will not be used in the rest of this
paper. Now we come back to ${\cal F}U_\epsilon(f;{\bf H}).$

\begin{definition}
$$\Gamma_f=\{\phi | \phi :\Sigma_f\rightarrow \Sigma_f \,\mbox{is an automorphism},
\, f\phi=f\}.$$
Here each $\phi$ is a self identification of $\Sigma_f$ as an
n-pointed \Fc.
\end{definition}

Since each constant component of $f$ is stable, it follows from the existence
of injective points for simple maps that $\Gamma_f$ is finite.

The next lemma explains why ${\cal F}U_\epsilon(f; {\bf H})$ forms a 
local uniformizer of ${\cal FB}(A).$

\begin{lemma}\label{uni}

When $\epsilon$ is small enough, there exists a continuous right action
of $\Gamma_f$ on ${\cal F}{\tilde U}_\epsilon(f; {\bf H})$, which is smooth on each
open strata of ${\cal F}{\tilde U}_\epsilon(f; {\bf H})$. The natural projection
$\pi_{\cal F}: {\cal F}{\tilde U}_\epsilon(f; {\bf H})\rightarrow {\cal FB}(A)$
commutes with $\Gamma_f$. The induced quotient map ${\bar \pi}_{\cal F}
: {\cal F}{\tilde U}_\epsilon(f; {\bf H})/{\Gamma}_f\rightarrow {\cal FB}(A)$ 
gives rise to a homomorphism of ${\cal F}{\tilde U}_\epsilon(f; {\bf H})/{\Gamma}_f$ and 
an open neighborhood of $\la f\ra $ in ${\cal FB}(A).$

\end{lemma}
\begin{remark}
Here the topology on  ${\cal F}{\tilde U}_\epsilon(f; {\bf H})$ is the $L^p_k$-topology,
which can be defined similar to what we did for ${\cal FB}(A)$. 
The smooth structure for each strata ${\cal F}{\tilde U}^D_\epsilon(f; {\bf H})$ is
the obvious one induced from the corresponding Banach manifold of product
of mapping spaces.
\end{remark}

\proof

Define 
$${\tilde \Gamma}_f=\{ \phi \,\|\, \phi:\Sigma_f\rightarrow \Sigma_f \mbox{
is an automorphism. } f\circ\phi\in {\cal F}U_\epsilon(f, {\bf H})\}.$$
We prove first that when $\epsilon>0$ is small enough, ${\tilde \Gamma}_f=
\Gamma_f.$ In fact let $y_f=\{f^{-1}(f(y_i)), f^{-1}(f(y_j^k))\}$ be the
collection of the inverse images of the images of those markings added for 
stabilizing $\Sigma_f$. $y_f$ is a finite set and ${\tilde \Gamma}_f$ is a subgroup
of the permutation group $Sym(y_f).$ Note that $f^{-1}(f(y_i))=y_i$. 
Hence $\phi(y_i)=y_i, \phi\in {\tilde \Gamma}_f. $ But the elements of
${\tilde \Gamma}_f$ may permute different bubble components which lie
on a same principal component. 
Let $m=\min_{\phi\in {Sym}(y_f)}\{\|f-f\circ\phi\|>0\}.$ It is easy to see that when
$0<\epsilon <<m$ we have ${\Gamma}_f={\tilde \Gamma}_f.$

Given $\phi\in \Gamma_f$ and $g\in {\cal F}U_\epsilon(f; {\bf H})$ with 
$g=g_{(u,v)} :\Sigma_{(u,v)}\rightarrow V.$ We want to define the right
action $g*\phi$.  When $\epsilon$ is small enough, $g_{(u,v)}$ and
$f_{(u,v)}$ are $C^1$-close to each other. This implies that near 
$\phi(y_i)=y_i$ and $\phi(y_j^k)$ of $\Sigma_{(u,v)}$ there exist points
$y_i(\phi, g)$ and $y_j^k(\phi, g)$ uniquely determined by $g$ such that 
$g(y_i(\phi,g ))\in {\bf H}_i$ and $g(y_j^k(\phi, g))\in
{\bf H}_j^k.$ Note that here $\phi(y_i)$ and $\phi(y_j^k)$ 
come from the corresponding points in ${ \Sigma}=\phi( \Sigma)$
through gluing.

Now consider stable \Fc  ${\Sigma}=\phi(\Sigma)$ equipped with 
markings $(\phi(y_i),$ $\phi(y_j^k),$ $x)$. There exists a marking preserving
identification 
$$\psi_{(u'', v'')}:(\phi(\Sigma)_{(u'',v'')}; \phi(y_i), \phi(y_j^k), x)
\rightarrow ({\Sigma}_{(u,v)}; y_i(g,\phi), g_j^k(g, \phi), x)$$
for some gluing parameter $(u'', v'').$ Clearly, there is also an identification
induced by $\phi, $
$$\phi_{(u',v')}:({\Sigma}_{(u',v')}, y_i, y_j^k, x)\rightarrow
(\phi(\Sigma)_{(u'',v'')}; \phi(y_i), \phi(y_j^k), x)$$ for
some $(u', v')$. Now we define
$$g*\phi=g\circ\psi_{(u'',v'')}\circ\phi_{(u',v')}:{\Sigma}_{(u'v,v')}
\rightarrow V.$$

Now we prove that for $g_1, g_2\in{\cal F}U_\epsilon(f, {\bf H}),$
$\la g_1\ra =\la g_2\ra \Longleftrightarrow  \exists \phi\in {\Gamma}_f$ such that
$g_1=g_2*\phi.$

 We only need to prove the $\Longrightarrow$ part.

 Suppose
$\la g_1\ra =\la g_2\ra $ with $g_i : {\Sigma}_{(u_i, v_i)}\rightarrow V,$ 
$ i=1, 2$. Then
there exists an identification of ${\cal F}$-curves ${\tilde \phi}:{\Sigma}_{(u_1, v_1)}
\rightarrow {\Sigma}_{(u_2, v_2)}$, which preserves the fixed marked
points $x$ and marked lines but may not preserve those $y$'s, such that
$g_1=g_2*{\tilde\phi}.$ Let $(y_f)_{(u_i, v_i)}\hookrightarrow
{\Sigma}_{(u_i, v_i)}, i=1,2$ be the finite subset of 
${\Sigma}_{(u_i, v_i)}$ corresponding to $y_f$ of ${\Sigma}_f$
through gluing. Consider ${\tilde\phi}(y_{(u_1,v_1)})\hookrightarrow
{\Sigma}_{(u_2, v_2)}$. Then for each element in 
${\tilde\phi}(y_{(u_1,v_1)})$, there is a unique element in 
$(y_f)_{(u_2,v_2)}  $ such that the elements of ${\tilde\phi}(y)$ lie in a
small disc centered at the corresponding element of $(y_f)_{(u_2, v_2)}$. 
This induces an injective map from ${\tilde\phi}(y)$ to $(y_f)_{(u_2, v_2)}$,
hence an injective map ${\tilde\phi}_y: y\rightarrow y_f.$ Now both
$g_1$ and $g_2$ are in ${\cal F}{\tilde U}_\epsilon(f, {\bf H}).$ When 
$\epsilon<<m$, this can happen only if ${\tilde\phi}_y$ is induced from some
$\phi\in {\Gamma}_f.$ Having obtained such a $\phi$, it is easy to see that
$g_1=g_2*\phi.$

Similarly, we can prove that for any $g\in{\cal F}{\tilde U}_\delta(f), $ when 
$\delta<< \epsilon, $ there exists some ${\tilde g}\in {\cal F}U_\epsilon
(f, {\bf H})$ such that $\la g\ra =\la {\tilde g}\ra .$ 

\QED

\noindent$\bullet$ {\bf  Orbifold bundles} 

\medskip

\noindent$\bullet$ Orbifold structure of ${\cal FB}(A)$ 
near ${\cal FM}(J, H, A).$

Let ${\cal F}U_\epsilon(f, {\bf H})=\pi_{\cal F}({\cal F}{\tilde U}_\epsilon
(f, {\bf H}))$, then ${\cal F}U_\epsilon (f, {\bf H})$ is an open neighborhood
of $\la f\ra \in {\cal FB}(A).$ Consider the open covering 
$${\cal FM}(J, H, A)\hookrightarrow \bigcup_{\la f\ra \in {\cal FM}(J, H, A)}{\cal F}
U_\epsilon(f, {\bf H}).$$ The compactness of ${\cal FM}(J, H, A) $ with
respect to the ( strong ) $L^p_k$-topology implies that there exist finite
$f_i, i=1, \cdots, m,$ such that 
$${\cal FM}(J, H, A)\hookrightarrow \bigcup_{i=1}^m {\cal F}U_{\epsilon_i}
(f_i, {\bf H}_i).$$

Now we use $U_i$ to denote  ${\cal F}U_{\epsilon_i}(f_i, {\bf H}_i)$ and
${\tilde  U}_i$ to denote its uniformizer. Let $U=\sum_{i=1}^m U_i.$

\begin{theorem}\label{orb}

$U$ has a stratified orbifold structure with respect to the local uniformizers.
\end{theorem}

\proof

The proof is a routine verification of the definition of orbifolds. We only
indicate the main step and leave the details to  our readers. 

We only need   to prove that if $\la g\ra \in  U_1\cap U_2$ with
$U_i$ being uniformized by ${\tilde U}_i$ with automorphism group $\Gamma_i$,
$i=1, 2,$ then there exists an open neighborhood $U$ of $\la g\ra $, such that
(i) $U\hookrightarrow U_1\cap U_2; $ (ii) $U$ is uniformized by ${\tilde U}$ with
automorphism group $\Gamma$ such that there exist two injection homomorphisms
$\imath_i:\Gamma\rightarrow\Gamma_i$ and two $(\Gamma, \Gamma_i)$-equivariant
embeddings $\lambda_i:{\tilde U}\rightarrow {\tilde U}_i$, $i=1,2.$ Here
the equivariant condition means that for any $\phi\in\Gamma, h\in {\tilde U},$ $ 
\lambda_i(h*\phi)=\lambda_i(h)*\imath_i(\phi), i=1, 2$.

 In our case, let 
$$g^i=g^i_{(u_i, v_i)}\in \la g\ra $$ be the representatives in ${\tilde U}_i
={\cal F}{\tilde U}_\epsilon (f_i, {\bf H}_i), i=1, 2. $ Then $g^i(y^i)\in
{\bf H}_i, i=1,2.$ Here we have used $y^i$ to denote the collection of all
those $y's$ in ${\Sigma}_{(u_i, v_i)}$ through the gluing, which are
originally in ${\Sigma}^i$ and are used to stabilize the unstable components
of $\Sigma^i.$

Now define 
$${\tilde V}_i=\{ g=g_{(u,v)}\,\|g_{(u, v)}-g^i_{(u, v)}\| <\delta_i,
g_{(u,v)}\in {\tilde U}_i|$$  and ${\tilde \Gamma}_i=\{\phi, \, |\,
\phi\in \Gamma_i, g^i*\phi=g^i\}.$ 

It is easy to see that when $\delta_i<<\epsilon_i, i=1,2, $
$(\pi_i)_{\cal F}:{\tilde V}_i/{{\tilde \Gamma}_i}\rightarrow {\cal FB}(A)$ is
an embedding onto some open neighborhood $V_i\hookrightarrow U_1\cap U_2.$
We may assume that $V_1=V_2=W.$ Therefore we get two uniformizer
$({\tilde V}_i, {\tilde \Gamma}_i)$ of $W$. Clearly the inclusion map
$({\tilde V}_i, {\tilde \Gamma}_i)\rightarrow ({\tilde U}_i, \Gamma_i)$
gives rise to an injective $({\tilde \Gamma}_i, \Gamma_i)$-equivariant map.
The theorem is valid if we can prove that $({\tilde V}_1, {\tilde \Gamma}_1)$ and
$({\tilde V}_2, {\tilde \Gamma}_2)$ are equivalent as uniformizers. Now choose
a forgetting marking process for $\Sigma^i_{(u_i, v_i)}$ deleting out
all extra markings in $y^i$ needed to make $\Sigma^i_{(u_i,v_i)}$ stable. 
Let ${\tilde y}^i$  be the remaining markings in
$\Sigma^i_{(u_i, v_i)}$ and ${\tilde{\bf H}}^i$  be the corresponding collection
of local hypersurfaces. Form ${\cal F}{\tilde U}_{\delta_i}({\tilde g}^i, {\tilde
{\bf H}}^i)$ and denote them by ${\tilde W}_i, $ where ${\tilde g}^i$ is same as
 $g^i$ as a map but its domain has on extra marking anymore. One can directly check tha there is an equivariant embeding of 
${\tilde W}_i$ into an open  subset of 
${\tilde V}_i$. Now since ${\tilde g}^i$ has no extra
markings after the forgetting marking process and $\la {\tilde g}_1\ra =\la {\tilde g}_2\ra $, 
we can easily construct an equivariant  equivalence of ${\tilde W}_1$ 
and ${\tilde W}_2$. 

\QED

\noindent $\bullet$ Local $T^{N_P}$-action on 
${\cal F}{\tilde U}_\epsilon^D(f, {\bf H})$.

\medskip


Let $D$ be an intersection pattern with $N_P$ principal components. We now 
define a local  $(S^1)^{N_P}$-action on ${\cal F}{\tilde U}_\epsilon^D(f; {\bf H}).$
We will call such an action a (local) toric $T^{N_P}$-action. Given $g\in 
{\cal F}{\tilde U}_\epsilon^D(f; {\bf H})$, let 
$g=\bigcup_{i=1}^{N_P}\{g^P_i\cup_k g^B_{i,k}\}$, where $g^B_{i,k}$ are the 
bubble components lying above $g_i^P.$ For any 
$$\phi=(\phi_1, \cdots, \phi_{N_P})\in (S^1)^{N_P},$$ with $|\phi|$ small, we 
will define $g*\phi$ by defining the action $\phi_i$ on 
$g^P_i\cup_k g^B_{i,k}, i=1, \cdots, N_P.$ 
 For this purpose we only need to know
how to define the action for $\phi\in S^1$, with $|\phi|$ small, on a stable map
$g\in {\cal F}{\tilde U}^D(f; {\bf H})$ in the following two simplest
cases:

\noindent (i) $g=g^P\cup g^B.$

Let $\Sigma=\Sigma_g=(P, l_P, d_P)\cup (B, y_1, y_2, d_B), $ where $y_1, y_2$
are marked points added  for stabilizing $B$ and $d_P=d_B$ is the double point
of $\Sigma$. Note that $g^B(y_i)\in {\bf H}_i,i=1,2.$ Choose an identification
$(P; l_P)\cong ({\bf R}^1\times S^1, \{\theta=0\})$ so that, under such an identification, 
$R_\phi$ of rotation of $\phi$-angle of $P$ is well-defined. Now the domain
$$\Sigma^\phi=\Sigma_{g*\phi}=(P; l_P,R^{-1}_\phi(d_P))\bigcup 
(B; y_1, y_2, d_B) $$ with double point $R^{-1}_\phi(d_P)=d_B.$

Let $I_\phi:\Sigma^\phi\rightarrow\Sigma $ be given by $R_\phi:P\rightarrow P$
and $Id:B\rightarrow B.$ Note that $I_\phi$ does not preserve the marked lines
$l_P$. We define $g*\phi=g\circ I_\phi.$

\noindent (ii) $g=g^P$ and the domain of $g$ has only one unstable principal component $\Sigma
=(P; l_p, y), y\in l_P$ and $g(y)\in {\bf H}.$ Consider 
$g|_{R_\phi(l_P)}: R_\phi(l_P)   \rightarrow V.$ When $|\phi|,\epsilon $ are small enough.
There is a unique $y_\phi\in I_\delta(R_\phi(y))$ for some given $\delta >0$
such that $g(y_\phi)\in {\bf H}.$ Let $T_\phi$ be the $s$-translation sending
$R_\phi(y)$ to $y_\phi.$ We define $g*\phi:(P; l_P, y)\rightarrow V$ given
by $g*\phi=g\circ T_\phi\circ R_\phi.$

\medskip

We summarize the properties of the $T^{N_P}$-action in the next lemma.

\begin{lemma}\label{act}

For any intersection pattern  $D$ of $N_P$ principal components, there exists a local $T^{N_P}$-action 
on  ${\cal F}{\tilde U}^D(f; {\bf H}), $ which is smooth in the variable 
of ${\cal F}{\tilde U}^D(f; {\bf H})$ and $C^l$-smooth in the variable of  
$T^{N_P},$ where $l=[k-\frac{2}{p}]$. Let 
${\cal F}{\tilde U}^{\bar D}(f; {\bf H})$  be the union of all 
${\cal F}{\tilde U}^{D_{1}}(f; {\bf H})$ with $D_{1} \leq D.$ Then the toric
$T^{N_{p}}$-action has a continuous extension to 
${\cal F}{\tilde U}^{ {\bar D}}(f; {\bf H}).$  The action is free , when all
$f^P_{i}, i=1,\cdots, N_{P},$ is $\theta$-dependent. Moreover, the action of
 $\Gamma_{f}$ commutes with  the local toric action.

\end{lemma}

\proof

We only need to prove the statement concerning  free action.
It follows from our assumption that  when $|u|$ is small,
each principal component of $f_{(u,0)}$ contain at least one point such that
$f_{(u,0)}$ is a local embedding near that point.When $\epsilon $ is small enough,
so does for any $$g\in  {\cal F}{\tilde U}^D(f; {\bf H}).  $$  The desired
conclusion follows from this.
 
\QED

It follows from this theorem that there is always an $S^1$-action on
${\cal F}{\tilde U}^{ {\bar D}}(f; {\bf H})$ via the diagonal map
from $S^1$ to $T^{N_{P}}.$ 
There is an obvious way to extend the $S^1$-action to
${\cal F}{\tilde U}_{\epsilon}(f; {\bf H})$, as we did for extending
the $\Gamma_{f}$-action. Since we will describe a similar process in next 
section, we  refer readers to there for this.

\medskip

\noindent$\bullet$ Orbifold bundle $({\cal L}, W)$.

\medskip

This is an infinite dimensional bundle $\cal L$ over ${\cal FB}(A)$ defined
as follows.

For any $\la f\ra \in {\cal FB}(A)$, we define 
$${\cal L}_{\la f\ra }=\bigcup_{f\in \la f\ra }{\cal L}^p_{k-1}( \wedge^{0,1} (f^*TV))/\sim ,$$
where the equivalence relation $\sim$ is defined via pull-back of the sections
induced from the identification of the domains.

Over $W\hookrightarrow {\cal FB}(A)$, $ {\cal L}$ has an orbifold bundle
structure. In fact, over each uniformizer 
${\tilde U}_i={\cal F}{\tilde U}_{\epsilon_i}(f_i; {\bf H}_i),$ there is a
bundle ${\tilde L}_i\rightarrow {\tilde U}_i$, which form a uniformizer
of ${\cal L}|_{U_i}$ with covering group $\Gamma_i$. For any $g\in 
{\tilde U}_i$, we define 
$$({\tilde{\cal L}}_i)_{(g)}=L^p_{k-1}(\wedge^{0,1}(g^* TV)).$$
The action of $\Gamma_i$ on ${\tilde U}_i$ lifts to ${\tilde {\cal L}}_i$ via pull-back. 
Similarly the local $T^{N_{P}}$-action on ${\tilde U}^D_i$ also lifts to
${\tilde{\cal L}}^D_i$ in the same way. 

The topology and smooth structure of $({\cal L}, W)$ can be described as
follows. 

Fix a gluing parameter $(u,v)$, we use ${\cal F}{\tilde U}^{(u,v)}_\epsilon
(f; {\bf H})$ to denote 
$$\{ g=g_{(u,v)}\,|\, \|g_{(u,v)}-f_{(u,v)}\|<\epsilon, \,\, g(y)\in {\bf H}\}.$$ 
Let ${\tilde {\cal L}}^{(u,v)}$ be the restriction of ${\tilde {\cal L}}$ to 
${\cal F}{\tilde U}^{(u,v)}_\epsilon(f; {\bf H}).$ Then   
${\tilde {\cal L}}^{(u,v)}$ has a local trivialization over 
${\cal F}{\tilde U}^{(u,v)}_\epsilon(f; {\bf H})$ induced from a
$J$-invariant parallel transformation of $(V, J)$ (see for example, [M] and
[F1]). This bundle structure for ${\tilde {\cal L}}^{(u,v)}$   gives rise to a
smooth structure for ${\tilde {\cal L}}^{(u,v)}$.

To define topology for ${\cal L}_W$, it is sufficient to define it for 
${\tilde {\cal L}}\rightarrow {\tilde U}$ as we did for ${\cal FB}(A)$.
Given $\xi\in ({\tilde{\cal L}})_g, g\in {\tilde U}, $ we define 
$\xi_{(u,v)}\in ({\tilde{\cal L}})_{g_{(u,v)}}$ for small $\|(u,v)\|$ as
follows. 

Without loss of generality, we may assume that $g$ has only two components 
$(\Sigma_i, d_i), i=1,2, $ with only one double point $d=d_1=d_2.$ Then
we define $\xi_{(u,0)}=\xi_{(0,0)}=\xi.$ Let $D_\delta(d_i)=\{w_i\,| |w|<\delta\}$
be the $\delta$-disc of $\Sigma_i$ centered at $d_i$ with complex coordinate 
$w_i$. We use $(s_i, \theta_i)$ to denote the corresponding cylindrical
coordinate. Choose a cut-off function
\[
\beta(s)=\left\{
\begin{array}{cl}
1 & s<-1\\
0 & s> 1.
\end{array}
\right.
\]

Note that over $D_\delta(d)$, $g_{(u,v)}(s,\theta)=g_{(u,0)}(s,\theta)$ if
$s<-\log |v|-2.$ We define 
$$\xi_{(u,v)}(s,\theta)=\beta(s_1+\log |v|)\cdot\xi^1_{(u,0)}(s_1,\theta_1)
+(1-\beta(s_1+\log|v|))\xi^2_{(u,0)}(s_2, \theta_2).$$

We now define an $\epsilon$-neighborhood of $\xi$ in ${\tilde {\cal L}}$ by
first define 
$${\tilde U}^{(u,v)}_\epsilon (\xi)=\{\eta=\eta_{(u,v )}\,|\,\eta_{(u,v)}
\in {\cal L}^p_{k-1}(g^*_{(u,v)}, TV), \, \|\eta_{(u,v)}-\xi_{(u,v)}\|_{k-1, p}
<\epsilon\}$$
for each fixed $(u,v).$ Then for each $(u,v)$, using the parallel transformation
to move ${\tilde U}^{(u,v)}_\epsilon(\xi)$ to the fiber of ${\cal L}^{(u,v)}$
over $h_{(u,v)}$,  with $\|h_{(u,v)}-g_{(u,v)}\|<\epsilon.$ We use 
${\tilde U}_\epsilon (\xi)$ to denote the collection of all images of 
${\tilde U}^{(u,v)}_\epsilon(\xi)$ under the parallel transformation. The 
collection of all ${\tilde U}_\epsilon(\xi)$ form a base of a topology on
${\tilde {\cal L}}\rightarrow {\tilde U}.$

Given the Hamiltonian function $H$, we can define a section $s_H^i:
{\tilde U}_i\rightarrow {\tilde{\cal L}}_i$ for the bundle 
${\tilde{\cal L}}_i$ as follows. For $f=f^P_j\cup f_k^B, $ we define 
$s_H^i(f)|_{B_k}\equiv 0$ and 
$$s^i_H(f)|_{P_j}=\nabla H\circ f^P_jds -J\nabla H\circ f^P_j d\theta.$$
${\bar{\ptl}}_J$-operator also induces an obvious section on 
${\tilde{\cal L}}_i\rightarrow{\tilde U}_i$ by sending $g\in {\tilde U}_i$ to
${\bar{\ptl}}_Jg=dg+J\circ dg\circ i.$ We use  
${\bar{\ptl}}^i_{J,H}$ to denote ${\bar{\ptl}}_J+s^i_H.$

We summarize what we have achieved this far in the following theorem. 
\begin{theorem}

There is a (stratified) orbifold bundle ${\cal L}$ over an open neighborhood
$W$ of ${\cal FM}(J, H, A)$ in ${\cal FB}(A).$ The orbifold structure on 
${\cal L}^D$ is compatible with the local $T^{N_P}$-action. The Hamiltonian
function $H$ induces a $\Gamma_i$-equivariant continuous section 
${\bar\ptl}^i_{J, H}$ of $({\tilde {\cal L}}_i, {\tilde U}_i)$, which
is smooth over each ${\tilde U}_i^{(u,v)}.$ Moreover when restricted to 
$({\tilde{\cal L}}^D_i, {\tilde U}^D_i)$, ${\bar{\ptl}}^i_{J,H}$ is
local $T^{N_P}$-equivariant.
\end{theorem}

\section{$T^{N_P}$-equivariant obstruction sheaf (local theory)}

In this section, we will construct an obstruction sheaf ${\tilde {\cal R}}_f$
over  a $T^{N_P}$-invariant uniformizer  
${\cal F}{\tilde U}_\epsilon(f^e; {\bf H}$  of $W$, where 
${\cal F}{\tilde U}_\epsilon(f^e; {\bf H}$ is the `completion' of 
${\cal F}{\tilde U}_\epsilon(f; {\bf H}$ so that the local $T^{N_{P}}$ action can be extended into a global one. 
The space ${\tilde K}_f=\Gamma({\tilde{\cal R}}_f)$ of the sections of
$ {\cal R}_f $
is a subspace
of $\Gamma({\tilde{\cal L}}, {\cal F}{\tilde U}_\epsilon (f^e; {\bf H})). $ 
Our
goal in this section is to prove Theorem \ref{ext}, which claims that each 
element of ${\tilde K}_f$, when restricted to 
${\cal F}{\tilde U}_\epsilon^D(f^e; {\bf H})$ is $(T)^{N_D}$-invariant, 
where $N_D$ is the number of principal components of $D$. 

Most of this section will be devoted to  construct directly the section space
${\tilde K}_f$. Our results here can be expressed as a solution to an abstract
extension problem. 

For any $f\in\la f\ra \in{\cal FM}(J, H, A)$, let  $({\tilde {\cal L}},
{\cal F}{\tilde U}_\epsilon(f;{\bf H})) $  be a local uniformizer with a local
$T^{N_{P}}$ -action on  each $({\tilde{\cal L}}^D, {\cal F}{\tilde U}^D_{\epsilon}(f; {\bf H})).$
Fix  a finite dimensional subspace $K\hookrightarrow {\cal L}_f
=\{\xi\,|\,\xi\in\wedge^{0,1}(f^* T(V))\}$ with the following properties

\noindent (i) there exists a $\delta>0$ such that for any element $\eta\in K,
\eta_{D_\delta}=\eta|_{D_\delta}=0$, where $D_\delta$ is the union of all 
$\delta$-discs centered at  double points of $\Sigma_f$;

\noindent (ii) for any $\eta\in K$, $\eta|_{P_i}=0$ if $f^P_i$ is a $\theta$-independent principal
component.

The main question that we want to answer  in this section is whether it is possible to extend each
element $\eta\in K$ into a local section ${\tilde \eta}$ of the bundle 
${\tilde{\cal L}}\rightarrow{\cal F}{\tilde U}^D_\epsilon(f; {\bf H})$ near
$f$ in such a way that

\noindent (a) ${\tilde \eta}$ is (locally) $T^{N_P}$-equivariant;

\noindent (b) ${\tilde\eta}$ is smooth on ${\cal F}{\tilde U}^{(u,v)}(f;{\bf H})$ for
any fixed $(u,v)$; 

\noindent  (c) ${\tilde \eta}$ is continuous. 

One can also formulate the same question for ${\cal F}{\tilde U}_\epsilon(f^e; {\bf H}$ with $T^{N_{P}}$ action.

There is a very simple way to extend elements in $K\rightarrow{\cal L}_f$ 
locally. The vanishing property $\eta_{D_\delta}=0$ for any $\eta\in K$ 
implies that $\eta$ can be extended over the local deformation $f_{(u,v)}$,
for small $\|(u,v)\|$, in an obvious way. We then use parallel transformation
to extend it to ${\cal F}{\tilde U}^{(u,v)}(f;{\bf H})$, for fixed $(u,v).$ 
Let $\eta^e$ be the resulting extension of $\eta. $ We will call it the canonical
extension of $\eta.$ In general $\eta^e$ may not be $T^{N_P}$-equivariant. One
may try to use the usual averaging process to make $\eta^e$ into a 
$T^{N_P}$-equivariant section. However there are two obvious difficulties that
make it impossible to directly use this usual averaging process. First of all
our $T^{N_P}$-action is only locally defined. Secondly, even in the case
that the usual averaging process is applicable to $\eta^e$, the resulting
$T^{P_{N}}$-equivariant section may not be equal to $\eta$ at $f$. 
In fact it may happen that it is even not close  to $\eta$ at $f$. In this
case , the transversality argument needed for gluing may fail. 

To see the nature of difficulties better,   we give a different 
description of ${\cal F}{\tilde U}_\epsilon(f;{\bf H}). $ 

We define a new uniformizer ${\cal F}{\tilde U}_\epsilon(f; {\tilde {\bf H}}).$
If there is no unstable $\theta$-dependent principal component of $f$, 
${\cal F}{\tilde U}_\epsilon(f; {\tilde {\bf H}})$  is just    
${\cal F}{\tilde U}_\epsilon(f:{\bf H}). $     Otherwise, we may assume that 
$f^P_i, i=1,\cdots, N^1_P$ are the $\theta$-independent unstable principal
components and $f_m^P, m=N^1_P+1, \cdots, N_P^2$ are the $\theta$-dependent
unstable principal components. Then for each such $\theta$-dependent
principal component $f_m^P$, there exists a local hypersurface 
${\tilde{\bf H}}_m$ at $f^P_m((0,0))$ of codimension two transversal to 
$f^P_m$ at $(0,0)$. Now for each such unstable component $P_m$, we will also
allow $y_m$, the marking used  for stabilizing $P_m,$  to vary in $\{s=0\}$. This will introduce a new real parameter
$\theta_m$ for the local deformation of $f$ and $\Sigma_f.$  It describes
the $\theta$-coordinate for $y_m.$ Let ${\tilde u}$ be the collection of 
the parameters in $u$ together with $\theta_m$'s. We now define the local
deformation $f_{\tilde u}=f\circ R_{\tilde u}$, where $R_{\tilde u}$ is
the rotation which  brings  $\theta_m$ to $(0,0)$ on each $P_m$. 

Let $D=D_f$. We set
\begin{eqnarray*}
& &{\cal F}{\tilde U}^D_\epsilon (f, {\tilde{\bf H}})=\\
& &\{ g=g_{\tilde u}
,| \|g_{\tilde u}-f_{\tilde u}\|_{k,p}<\epsilon, g^P_i(y_i)\in {\bf H}_i, 
g^P_m(y_m)\in {\tilde {\bf H}}_m, g^B_j(y^k_j)\in {\tilde {\bf H}}^k_j \}.
\end{eqnarray*}

We can  define ${\cal F}{\tilde U}_\epsilon (f, {\tilde{\bf H}})$  
similarly by using the local deformation $f_{({\tilde u}, v)}$, where
$f_{({\tilde u}, v)}$ is obtained from $f_{\tilde u}$ by the   gluing
with gluing parameter $v$. 

We now define the corresponding $T^{N_P}$-action. We start with the action
on ${\cal F}{\tilde U}^D_\epsilon (f, {\tilde{\bf H}}).$  Clearly, we only
need to define $g_m^P*\phi$, for $g\in 
{\cal F}{\tilde U}^D_\epsilon (f, {\tilde{\bf H}})$ and $\phi\in S^1. $ 
The domain of $g_m^P*\phi$ is $(P_m, l_m, R^{-1}_\phi(y_m))$, and we define
$g_m*\phi=g^P_m\circ R_\phi.$

We now extend the induced $S^1$-action to
${\cal F}{\tilde U}_\epsilon (f, {\tilde{\bf H}})$. Given $g=g_{({\tilde u}, v)}$
and $\phi\in S^1, $ we define the $\phi$-rotation of the principal components,
still denoted as $R_\phi$, 
$$R_\phi:(\Sigma_{(\tilde u, v)}; y_i, R^{-1}_\phi(y_m), y^k_j)\rightarrow
(\Sigma_{(\tilde u, v)}; R_\phi(y_i), y_m, y^k_j),$$
which is a rotation of $\arg \phi$ on each principal component with respect
to marked line there and is identity on each bubble component. Consider 
$g\circ R_\phi$. There are unique $y_i'\in I_\delta(y_i), y_m'\in D_\delta
(R^{-1}_\phi(y_m))$ and ${y_j^k}'\in D_{\delta}(y^k_j)$ such that
$g\circ R_\phi(y_i')\in {\bf H}_i$,\,\, 
$g\circ R_\phi(y_m')\in{\tilde{\bf H}}_m$ 
\,\,and \,\, $g\circ R_\phi(y_j^k)\in {\tilde {\bf H}}^k_j.$ Then there is a unique new
parameter $({\tilde u}', v')$ such that 
$$(\Sigma_{({\tilde u}', v')}; y_i, R^{-1}_\phi(y_m), y^k_j)\cong (\Sigma_{(\tilde u, 
v)}, y'_i, y_m', {y^k_j}')$$
under the identification map $\psi: \Sigma_{({\tilde u}', v')}\rightarrow
\Sigma_{(\tilde u, v)}.$ We define $g*\phi= g\circ R_\phi\circ\psi.$

We remark that one can use a similar construction to extend the local 
$S^1$-action on ${\cal F}{\tilde U}^D_\epsilon (f, {{\bf H}})$
defined in previous section to 
${\cal F}{\tilde U}_\epsilon (f, {{\bf H}})$. 

Now observe that the action above is well-defined for all 
$\phi\in T^{N_P}$,
not just locally for $|\phi|$ small. 
This suggests us to enlarge ${\cal F}{\tilde U}_\epsilon 
(f, {{\bf H}})$ so that $T^{N_P}$-action can be defined globally. 
For this purpose, we fix a large integer $M$. 
 For each element 
$I=(i_1, \cdots, i_{N_P})\in ({\bf Z}/M{\bf Z})^{N_P}$, set 
$$\phi_I=(\phi_{i_1}, \cdots, \phi_{i_{N_P}})$$ where $\phi_{i_k}=\frac{2\pi i_k}
{M}.$ We define $f_I$ and ${\cal F}{\tilde U}_\epsilon (f_I, {{\bf H}})$
as follows. $f_I$ is same as $f$ when restricted to any stable principal
components
or bubble components. When restricted to unstable principal component $P_m$,
the domain of $f_I$ is $(P_m, l_m, R^{-1}_{\phi_{i_m}}(y_m)).$ Note that
here the marked point $R^{-1}_{\phi_{i_m}}(y_m)$ used for stabilizing $P_m$ for
$f_I$ lies on $R^{-1}_{\phi_{i_m}}(l_m).$ We define $f_I$ to be $f\circ
 R_{\phi_{i_m}} $ on $P_m.$ Now we can define ${\cal F}{\tilde U}_\epsilon
(f_I, {\bf H})$ and local $T^{N_P}$-action on it by the very same formula as we
did for ${\cal F}{\tilde U}_\epsilon (f, {{\bf H}})$. 

Now form the disjoint union 
$$\coprod_{I}{\cal F}{\tilde U}_\epsilon (f_I, {{\bf H}}), I\in 
({\bf Z}/M{\bf Z})^{N_P}.$$   
We introduce an equivalence relation between the elements in    
${\cal F}{\tilde U}_\epsilon (f_I, {{\bf H}})$ and elements in its 
neighbors. A typical neighbor has a form 
${\cal F}{\tilde U}_\epsilon (f_J, {{\bf H}})$ with all $j_k=i_k, k\not =l$
and $j_l=i_l+1$ for some $1\leq l\leq N_P.$ Without loss of generality, we
may assume that $l=1$, $I=(i, {\hat I}), J=(i+1, {\hat I}). $
Given such $I$ and $J$ there is a coordinate change $T_{I,J}$ from 
${\cal F}{\tilde U}_\epsilon (f_I, {{\bf H}})$ to 
${\cal F}{\tilde U}_\epsilon (f_J, {{\bf H}})$ as follows. 
For any $g\in {\cal F}{\tilde U}^D_\epsilon (f_I, {{\bf H}})$,
$D=D_f$, we define $T_{I,J}(g)=g$ along all stable principal components, 
bubble components and those unstable principal components $g_l^P, l>1.$
Now consider the unstable principal component $g^P_1.$ Since 
$g^P_1(R^{-1}_{\phi_i}(y_1))\in {\bf H}_1$, there exists a unique
$y'_{i+1}$ lying on $R^{-1}_{\phi_{i+1}}(l_1)$ and near $R^{-1}_{\phi_{i+1}}(y_1)$
such that $g^P_1(y'_{i+1})\in {\bf H}_1, $ when $g$ is close enough to the
local deformation of $f_I$ and $M$ is large enough. Let $T_r$ be the 
$s$-translation for $P_1$ that brings $R^{-1}_{\phi_{i+1}}(y_1)$ to
$y'_{i+1}$. We define $T_{I,J}$ along $g^P_1$ to be $g^P_1\circ T_r.$ One can 
easily extend $T_{I,J}$ from fixed intersection pattern $D$ to the general case
by using a similar process of extending $\Gamma_f$-action before.

Let ${\cal F}{\tilde U}_\epsilon(f^e; {\bf H})$ denote the above disjoint
union quotienting out the equivalence relation introduced by these $T_{I,J}$'s.

\begin{lemma}\label{main}

${\cal F}{\tilde U}_\epsilon(f^e; {\bf H})$ is a ( stratified ) Banach
manifold. The local $T^{N_{P}}$-action on ${\cal F}{\tilde U}_\epsilon(f_I; {\bf H})$ 
are compatible with the ``coordinate changing" maps $T_{I,J},$
$ I, J\in 
({\bf Z}/M{\bf Z})^{N_P}. $ This defines a local $T^{N_P}$-action on 
${\cal F}{\tilde U}_\epsilon(f^e; {\bf H})$ , which can be  extended into
a (global ) $T^{N_P}$-action.

\end{lemma}

\proof

One can directly check that the local $T^{N_P}$-action on each coordinate
chart ${\cal F}{\tilde U}_\epsilon(f_I; {\bf H})$ is preserved under coordinate
changes $T_{I,J}$. The usual process to complete a local action to a global
one for a compact Lie group is applicable here, which yields a well-defined
$T^{N_P}$-action. 

\QED

\begin{lemma}

There is a $T^{N_P}$-equivariant equivalence map 
$$\rho: {\cal F}{\tilde U}_\epsilon(f^e; {\bf H}) \rightarrow 
{\cal F}{\tilde U}_\epsilon(f; {\tilde {\bf H}}) .$$

\end{lemma}

\proof

We only need to define 
$$\rho^D: {\cal F}{\tilde U}^D_\epsilon(f^e; {\bf H})\rightarrow 
{\cal F}{\tilde U}^D_\epsilon(f;{\tilde  {\bf H}})$$ with $D=D_f$. The extension from
$\rho^D$ to $\rho$ is a routine procedure. We omit it here. 

Without loss of generality, we may assume that $f$ contains some unstable
principal component $f^P_m$, $m< N_P$, which is $\theta$-dependent. 
Given $g\in {\cal F}{\tilde U}_\epsilon(f^e; {\bf H})$, we may assume that
$g\in {\cal F}{\tilde U}_\epsilon(f; {\bf H})=
{\cal F}{\tilde U}_\epsilon(f_0; {\bf H}).$ Clearly, we only need to define
$\rho^D$ along 
 those unstable principal component $g^P_m$. There exists a
unique point ${\tilde y'}_m$ near $y_m$ in the domain $P_m$ such that
$g^P_m({\tilde y'}_m)\in {\tilde {\bf H}}_m\hookrightarrow {\bf H}_m.$

We define the domain of the $m$-th unstable principal component of $\rho(g)$ to be 
$(P_m, l_m, {\tilde y}_m),$ where ${\tilde y}_{m}$ has same
 $\theta$-coordinate as ${\tilde y'}_m$ has  and  lies
on th central circle $\{s=0\}.$ Now we define
$$(\rho(g))^P_m=g^P_m\circ Tr^m: (P_m, l_m, {\tilde y}_m)\rightarrow V,$$
where $Tr^m$ is the $s$-translation of $P_m$ that brings ${\tilde y}_m$
to ${\tilde y'}_m.$ One can directly check that the above definition
of $\rho$, which is given by using particular coordinate of 
${\cal F}{\tilde U}_\epsilon(f^e; {\bf H})$ is actually compatible with 
coordinate changing maps $T_{I, J}$. A direct calculation shows that $\rho$
is one-to-one and commutes with the $T^{N_P}$-actions defined on the 
corresponding spaces. 

\QED

Because of this lemma, we may use 
${\cal F}{\tilde U}_\epsilon(f;{\tilde {\bf H}})$ to replace 
${\cal F}{\tilde U}_\epsilon(f^e; {\bf H})$ for our problem of finding a 
$T^{N_P}$-equivariant extension of $K$.  

\begin{theorem}\label{ext}

There exists an extension of $K$ over 
${\cal F}{\tilde U}_\epsilon(f^e; {\bf H})$ , which has the property described at
the beginning of this section. 

\end{theorem}

\proof 

We only need to prove the corresponding statement for 
${\cal F}{\tilde U}_\epsilon(f; {\tilde {\bf H}})$.
The following two cases are to be considered. 

\noindent (i) All unstable principal components $f^P_m$ are $\theta$-dependent;
\noindent (ii) some of the unstable principal components are 
$\theta$-independent. 

\medskip
\noindent Case (i): Recall that in this case there is another kind of local
deformation $f_{(\alpha, {\tilde t})}.$ The domain 
$\Sigma_{(\alpha, {\tilde t})}$
of $f_{(\alpha, {\tilde t})}$ is obtained from $\Sigma=\Sigma_f$, with
the gluing parameter $(\alpha, {\tilde t})$, where $\alpha$ is the first
component of $u=(\alpha, \theta)$ and ${\tilde t}$ is obtained from $v=(t,\tau)$ by
adding a rotational component to $\tau$ (see Section 2.1 ). 
$f_{(\alpha, {\tilde t})}$ is obtained from $f_{(\alpha, 0)}$ through a similar
gluing process to the previous one with the gluing parameter $\tilde t.$ We
define
$${\cal G}{\tilde U}_\epsilon(f;{\tilde  {\bf H}})=\{g=g_{(\alpha, {\tilde t})}\, |\,
|g_{(\alpha, {\tilde t})}-f_{(\alpha, \tilde t)}|<\epsilon, g(y_i)\in 
{\tilde H}_i, g(y^k_j)\in {\tilde H}^k_j\}.$$
There is an obvious projection map
$${ p}_f: {\cal F}{\tilde U}_\epsilon(f; {\tilde{\bf H}})\rightarrow 
{\cal G}{\tilde U}_\epsilon(f;{\tilde  {\bf H}}),$$
given by simply forgetting the marked lines
of the elements of ${\cal F}{\tilde U}_\epsilon(f; {\tilde {\bf H}})$. Let 
${\cal L}^{\cal F}$ and ${\cal L}^{\cal G}$ be the corresponding bundles
over the above two spaces. Clearly, $({ p}_f)^*({\cal L}^{\cal G})={\cal L}
^{\cal F}.$ We may identify $K\hookrightarrow \Gamma({\cal L}^{\cal F}_f)$ with
a subspace of $\Gamma({\cal L}^{\cal G}_f).$ Let $K'$ be the
corresponding subspace. Suppose that $\dim K=k$ and $\{e_1, \cdots, e_k\}$ is
a basis for $K$. Let $e_i', i=1, \cdots, k$ be the corresponding elements in
$K'$ with $p^*_f(e_i')=e_i.$

Since each $e_i'$ vanishes near all double points of $f$, we can easily 
get a canonical extension $(e_i')^c$ of a section of the bundle ${\cal L}^{\cal G}$ 
over ${\cal G}{\tilde U}_\epsilon(f;{\tilde  {\bf H}})$ as we did for stable maps
with ${\cal F}$-curves as domains described before in this section. We define
${\tilde e}_i=p^*_f((e'_i)^c).$ 

Now for any $\phi\in T^{N_P(g)}$ and 
$g\in {\cal F}{\tilde U}_\epsilon(f; {\tilde{\bf H}})$, 
$p_f(g*\phi)=p_f(g).$\label{comm}
This implies that ${\tilde e}_i, i=1, \cdots, k$, has the required invariant property.

\medskip

\noindent Case (ii): 

We may assume in this case that $f^P_i, i=1, \cdots, N_P^1$ be the 
$\theta$-independent unstable principal components and $f^P_m, m=N_P^1+1, \cdots, 
N_P^2$ be the other unstable principal components. Note that in this case
$N_P^1<N_P.$

We now define a local deformation of $f$, which is  a mixture of those
deformations as stable $ {\cal F}$-maps and stable ${\cal G}$-maps. The idea
is to deform the $\theta$-independent part of unstable principal components
as stable $ {\cal F}$-maps and the rest as stable ${\cal G}$-maps. When the topological type  is fixed,  the domain
of such a deformation is described by the parameter $\alpha$ appeared in the 
corresponding deformation as stable ${\cal G}$-maps, since in this case all those markings $y_m$ of $P_m$, $m\leq N_P^1, $
are fixed. The gluing parameters that control the topological type of the
deformation can be described as follows. 

On \Fc part and \Gc part of $\Sigma_f$, we use the usual gluing parameter respectively. 
That is we associate the gluing parameter $\tau=(\tau_2, \cdots, \tau_{N_P^1})$ to
the ends $z_2, \cdots, z_{N_P^1}$ and ${\tilde t}$ to those double points
and ``ends" in \Gc part. Now the key point is to associate the ``ends"
$z_{N_P^1+1}$, the double point that divides $\Sigma_f$ into the two parts, 
with a complex gluing parameter $w.$ We will use ${\tilde t}'$ to denote 
$(\tau, w, \tilde t)$. Set $\alpha'=\alpha. $ We have the deformation 
$ \Sigma_{(\alpha', {\tilde t}')}$ of $\Sigma$. The deformation of $f_{(\alpha', 0)}$
and $f_{(\alpha', {\tilde t}')}$ can be defined in a similar way as before.

We now define
$$
{\cal H}{\tilde U}_\epsilon(f; {\tilde{\bf H}})=\left\{      
g_{(\alpha', {\tilde t}')} \left\vert\begin{array}{cl}
|g_{(\alpha', {\tilde t}')}-f_{(\alpha',{\tilde t}')}|<\epsilon, &
g_{(\alpha', {\tilde t}')} (y^k_j)\in {\tilde {\bf H}}_j^k ,\\
g_{(\alpha', {\tilde t}')}(y_i)\in {\bf H}_i, &   
g_{(\alpha', {\tilde t}')} (y_m)\in {\tilde {\bf H}}_m, 
\end{array}
\right.
\right\},
$$
where $y_i\in P_i$ of a $\theta$-independent unstabl  principal component, 
$y_m\in P_m$ of a  $\theta$-dependent unstabl  principal componen and 
$y^k_j \in B^k_j$  of an unstable bubble component.

Fix an intersection pattern $D.$ We have the following two cases:
 
\noindent(a) the gluing parameter $w$  determined by $D$ is zero.

In this case, we may assume that each element $g$ in the strata has
$N_P(D_1)$ those principal components, whose domains are  obtained
from the domains of $\theta$-independent unstable principal components
of f through gluing. Let $N_P(D_2)$  be  the number of the other
principal components of g.
We decompose  ${\cal F}{\tilde U}^D_\epsilon(f; {\tilde{\bf H}})$ as 
the product,  
$${\cal F}{\tilde U}^{D_1}_\epsilon(f; {\tilde{\bf H}})
\times {\cal F}{\tilde U}^{D_2}_\epsilon(f; {\tilde{\bf H}}),$$ where the first
factor contains those $g^P_i, 1\leq i\leq N_P(D_1)$ and the second contains
all the other components of $g$.  Similarly, we can decompose 
${\cal H}{\tilde U}^D_\epsilon(f; {\tilde{\bf H}})$ as 
the product,  ${\cal H}{\tilde U}^{D_1}_\epsilon(f; {\tilde{\bf H}})
\times {\cal H}{\tilde U}^{D_2}_\epsilon(f; {\tilde{\bf H}})$. Let
$p_i: {\cal F}{\tilde U}^D_\epsilon(f; {\tilde{\bf H}})\rightarrow
{\cal F}{\tilde U}^{D_i}_\epsilon(f; {\tilde{\bf H}}), i=1,2,  $ be the 
projection. Then we define $$p^D_f(g)=(p_1(g), p^{D_2}_fp_2(g)),$$ where
$$p^{D_2}_f:{\cal F}{\tilde U}^{D_2}_\epsilon(f; {\tilde{\bf H}})\rightarrow 
{\cal H}{\tilde U}^{D_2}_\epsilon(f; {\tilde{\bf H}})$$  is defined as in
case (i). Clearly $p^D_f$ commutes with the $T^{N_P(D_2)}$-action on the
second factor. 

\noindent (b) $w\not = 0.$ In this case, we need to add one more factor
${\cal F}{\tilde U}^{D_3}_\epsilon(f; {\tilde{\bf H}})$ to the above 
decomposition of ${\cal F}{\tilde U}^D_\epsilon(f; {\tilde{\bf H}})$ , which
corresponds to the principal component of $g$ ``passing through" $w$
together with all bubbles lying on this component. Let 
${\cal H}{\tilde U}^{D_i}_\epsilon(f; {\tilde{\bf H}}), i=1,2,3, $ be the 
corresponding decomposition of 
${\cal H}{\tilde U}^D_\epsilon(f; {\tilde{\bf H}})$ . We define 
$p_i, i=1,2,3$ similarly. Now one can directly verify that
${\cal F}{\tilde U}^{D_3}_\epsilon(f; {\tilde{\bf H}})$ and
${\cal H}{\tilde U}^{D_3}_\epsilon(f; {\tilde{\bf H}})$ are homeomorphic to each
other( in the case $w\not= 0$) and diffeomorphic when restricted to subspace with fixed deformation
type. Let $$p^{D_{3}}_f :{\cal F}{\tilde U}^{D_3}_\epsilon(f; {\tilde{\bf H}})
\rightarrow {\cal H}{\tilde U}^{D_3}_\epsilon(f; {\tilde{\bf H}})$$ be
the corresponding homeomorphism.
We define $$p^D_f(g)=(p_1(g), p^{D_2}_f\cdot p_2(g), p^{D_{3}}_f(p_3(g))).$$ Again
$p^D_f$ commutes with the $T^{N_P(D_2)}$-action, where $N_P(D_2)$ is the
number of principal components in 
${\cal F}{\tilde U}^{D_2}_\epsilon(f; {\tilde{\bf H}})$. One can prove that all 
these $p^D_f$'s, pasted together, define a continuous map
$$p_f:{\cal F}{\tilde U}_\epsilon(f; {\tilde{\bf H}})\rightarrow
{\cal H}{\tilde U}_\epsilon(f; {\tilde{\bf H}}),$$ which is smooth when 
restricted to each subspace of fixed deformation type. 

Now as in case (i), we can easily get a canonical extension ${\tilde e}'_i, 
i=1, \cdots, k, $ of $e'_i$ where ${\tilde e}'_i$ is a section of the
bundle ${\cal L}^{\cal H}$ over ${\cal H}{\tilde U}_\epsilon(f;{\tilde {\bf H}}).$ 
We then consider ${\tilde b}_i=p^*_f({\tilde e}'_i)$ of the corresponding
section of ${\cal L}^{\cal H}$ over 
${\cal F}{\tilde U}_\epsilon(f; {\tilde{\bf H}})$. 
Unlike the case (i), in this
case, ${\tilde b}_i$ does not have the required invariant property. In fact, given
$g\in {\cal F}{\tilde U}_\epsilon(f; {\tilde{\bf H}})$, if $D=D_g$
is in the case (a) above, then ${\tilde b}_i$ is already $T^{N_P(D)}$-invariant
and we simply define ${\tilde e}^D_i={\tilde b}^D_i$ over
${\cal F}{\tilde U}^D_\epsilon(f; {\tilde{\bf H}})$, since  in this case, 
 $p^D_f$ commutes with the  $T^{N_P(D_2)}$-action acting on 
${\cal F}{\tilde U}^{D_2}_\epsilon(f; {\tilde{\bf H}})$ and  $b_i$ can  be 
chosen to be zero along those elements of 
${\cal F}{\tilde U}^{D_1}_\epsilon(f; {\tilde{\bf H}})$ by our assumption
on $K$. In the case that $D$ is in case (b), ${\tilde b}_i^D$ so defined
only has $T^{N_P(D_1)+N_P(D_2)}$-invariant. There is an extra $S^1$-action
coming from the ``rotations" of elements of 
${\cal F}{\tilde U}^{D_3}_\epsilon(f; {\tilde{\bf H}})$. Given $g\in 
{\cal F}{\tilde U}^{D}_\epsilon(f; {\tilde{\bf H}})$ and $\phi\in S^1$, we 
use $$A^D_\phi:{\cal F}{\tilde U}^{D}_\epsilon(f; {\tilde{\bf H}})
\rightarrow {\cal F}{\tilde U}^{D}_{\epsilon}(f; {\tilde{\bf H}}) $$ 
to denote this action for fixed $\phi\in S^1$. We  define
$${\tilde e}^D_i=\frac{1}{2\pi}\int_{S^1} {A^D_\phi}^*({\tilde b}^D_i) 
d\phi.$$
It follows from the definition that ${\tilde e}^D_i$ is $T^{N_P(D)}$-invariant. 
One can directly verify that ${\tilde e}^D_i$ is compatible to each other 
when $D$ varies and hence gives rise to a well-defined continuous
extension ${\tilde e}_i$ of $e_i$ with all the required properties of the theorem.

We note that in the case (ii) above, we have assumed that the first 
$N^1_P$ principal components of $f$ are unstable and $\theta$-independent. 
This assumption simplifies the way of choosing gluing parameter ${\tilde t}'$
of $f_{(\alpha', {\tilde t}')}$ and other related constructions. The general
case can be treated in a similar manner, but with more complicated notations. 

Now let ${\tilde K}=span \{ {\tilde e}_1, \cdots, {\tilde e}_k\}.$

\QED

\begin{lemma}

When $p$ is even, there exists a $T^{N_P}$-invariant continuous cut-off
function ${\tilde \beta}$ on ${\cal F}{\tilde U}_\epsilon(f^e; {\bf H})$
with $0\leq {\tilde\beta}(g)\leq 1$  such that  ${\tilde \beta} (g)=1$ 
 in a neighborhood of the $T^{N_P}$-orbit
of $f$ and ${\tilde\beta}=0$ near the boundary of 
${\cal F}{\tilde U}_\epsilon(f^e; {\bf H})$. Moreover, ${\tilde\beta}$ is smooth
on ${\cal F}{\tilde U}^{(u,v)}_\epsilon(f^e; {\bf H})$.

\end{lemma}

\proof

Consider 
$$p_f:{\cal F}{\tilde U}_\epsilon(f; {\tilde {\bf H}})\rightarrow 
{\cal H}{\tilde U}_\epsilon(f;{\tilde  {\bf H}})$$ in Theorem \ref{ext}. Suppose
that  we can construct a continuous cut-off function
$\beta_1:{\cal H}{\tilde U}_\epsilon(f; {\bf H})\rightarrow [0,1]$ such    
that $\beta_1=1$ in a neighborhood of $p_f(f)$ and $\beta=0$ near the boundary
of ${\cal H}{\tilde U}_\epsilon(f;{\tilde  {\bf H}})$. Repeating the process of finding
${\tilde e}_i$ with the required smoothness in previous Theorem, we set
${\tilde\beta}_1=p_f(\beta_1)$ and ${\tilde \beta}^D=\frac{1}{2\pi}\int
(A^D_\phi)^*({\tilde\beta}^D_1) d\phi.$ As before, ${\tilde \beta}^D$
will be pasted together to get a $T^{N_P}$-invariant cut-off function
$\tilde\beta$ with the desired property. To construct $\beta_1$, we note 
that when $p$ is even, the function $\rho^g_{(u,v)}(h)=\|h-g\|^p_{k,p}$ is
a smooth function on ${\cal F}{\tilde U}^{(u,v)}_\epsilon(f; {\bf H})$. Using    
this we can easily construct the desired $\beta_1.$

\QED

\begin{cor}

In Theorem \ref{ext}, we may choose ${\tilde K}$ in such a way that all its 
elements vanish near the boundary of ${\cal F}{\tilde U}_\epsilon(f^e; {\bf H}).$    

\end{cor}

\section{Transversality and Gluing}

In this section we will establish the transversality of perturbed
${\bar\ptl}_{J, H, \nu}$-operation over a $T^{N_P}$-invariant
uniformizer ${\cal F}{\tilde U}_\epsilon(f^e; {\bf H})$, where the perturbation    
term $\nu$ is a generic element of ${\tilde K}$. Therefore ${\bar\ptl}_{J,H,}^{\nu}$
is a $T^{N_P}$-equivariant transversal section of ${\cal L}$. Its zero set
${\cal F}{\tilde{\cal M}}^{\nu}(f^e)$ is a cornered smooth manifold  of
dimension $Ind (c_+)-Ind (c_-)+2c_1(A)+1$, with a $T^{N_P}$-action acting on it.

Our method in this section is an adaption of the method in [LiuT1] ( see also
[L1]). In this section often we will only quote results in [LiuT1] and indicate
necessary changes to incorporate the $T^{N_P}$-action here. We refer reader
for the detailed proof in [LiuT1].

\subsection{Transversality}

We start with giving a local coordinate charts for 
${\cal F}{\tilde U}_\epsilon(f; {\bf H})$  near $f$   and local 
trivialization  of ${\cal L}$ over those coordinate charts. 

Let $y_i\in P_i$ and $y^k_j\in B_j, 1\leq k\leq 2$, be the marked points
added to $\Sigma_f$ for stabilizing it, and ${\bf H}_i$ and ${\tilde {\bf H}}^k_j$ be
the local hypersurfaces at $f^P_i(y_i)$ and $f^B_j(y^k_j)$ used before for
slicing.

For each ${\bf H}_i$ and ${\tilde {\bf H}}^k_j$, let $h_i=T_{f(y_i)}{\bf H}_i$, 
$h^k_j=T_{f(y^k_j)}{\tilde {\bf H}}_j^k$. We may assume that both ${\bf H}_i$ and
 ${\tilde {\bf H}}_j^k$ are totally geodesic so that they are the local images of $h_i$ and $h^k_j$ under the exponential map. 

We define 
$$L^p_k(f ^*TV,h)=\{\xi, \, |  \, \xi\in L^p_k(f^* TV), 
\, \xi (y_i)\in h_i, \xi(y^k_j)\in h^k_j  \},$$ 
where the values of $\xi$ from different components are 
the same at double points. 

Similarly, we define
 $L_k^p(f_{(u, v)}^*TV, h).$

 Let
$${\tilde V}^{(u,v)}_\epsilon=\{\xi\,| \xi \in L^p_k(f_{(u,v)}^* TV; h); 
\, \|\xi\|_{k,p}<\epsilon\},$$ and 
${\tilde V}_\epsilon=\bigcup_{\|(u,v)\|<\delta} {\tilde V}^{(u,v)}_\epsilon $, 
which is a ``bundle" over $\Lambda_\delta=\{(u,v)\,|\, \|(u,v)\|<\delta\,\}.$

The coordinate chart
$$ Exp^{(u,v)}_f:{\tilde V}^{(u,v)}_\epsilon\rightarrow
{\cal F}{\tilde U}^{(u,v)}_\epsilon(f; {\bf H})$$ 
for ${\cal F}{\tilde U}_\epsilon
^{(u,v)}(f;{\bf H})$ is given by
$\xi\rightarrow Exp_{f_{(u,v)}}\xi.$

The ``coordinate chart" for ${\cal F}{\tilde U}_\epsilon(f; H)$ is given
by
$$Exp_f=\bigcup_{(u,v)\in \Lambda_\delta} Exp^{(u,v)}_f:{\tilde V}_\epsilon=\bigcup_{(u,v)\in \Lambda_\delta}
{\tilde V}^{(u,v)}_\epsilon\rightarrow {\cal F}{\tilde U}_\epsilon (f; {\bf H})
=\bigcup_{(u,v)\in \Lambda_\delta}
{\cal F}{\tilde U}^{(u,v)}_\epsilon(f; {\bf H}).$$

We already defined trivialization of ${\cal L}^{(u,v)}$ over 
${\cal F}{\tilde U}^{(u,v)}_\epsilon(f; {\bf H})$ by using the parallel 
transformation. Let
$$\psi^{(u,v)}:{\cal F}{\tilde U}^{(u,v)}_\epsilon(f; {\bf H})\times
L^p_{k-1}(\wedge^{0,1}(f^*_{(u,v)}TV))\rightarrow
{\cal L}^{(u,v)}$$
denote the trivialization here. Then 
$$\gamma^{(u,v)}=\psi^{(u,v)}\circ(Exp^{(u,v)}_f\times
Id):{\tilde V}^{(u,v)}_\epsilon\times
L^p_{k-1}(\wedge^{0,1}(f^*_{(u,v)}TV))\rightarrow {\tilde{\cal L}}^{(u,v)}$$ gives
rise to a trivialization of ${\cal L}^{(u,v)}$ in terms of above coordinate
chart. 

Let $\gamma=\cup_{(u,v)\in\Lambda_\delta}\gamma^{(u,v)}.$ Then 
$$\gamma:\bigcup V^{(u,v)}_\epsilon\times 
L^p_{k-1}(\wedge^{0,1}(f^*_{(u,v)}TV))
\rightarrow{\tilde{\cal L}}$$ is a local ``trivialization" of 
${\tilde{\cal L}}$.

Now under these local coordinate chart and local trivialization, the 
${\bar{\ptl}}_{J,H}$-section of ${\tilde{\cal L}}$ becomes:
$$F^1_{(u,v)}=\pi_2\circ(\gamma^{(u,v)})^{-1}\circ{\bar{\ptl}}_{J,H}\circ 
Exp^{(u,v)}_f: V^{(u,v)}_\epsilon
\rightarrow  {\cal L}^p_{k-1}(\wedge^{0,1}(f^*_{(u,v)}TV)). $$

Let $F^1=\bigcup_{(u,v)\in \Lambda_\delta} F^1_{(u,v)}$. 
Note that $F^1_{(u,v)}$ is smooth. 
Let $$L^1_{(u,v)}=(DF^1_{(u,v)})_{f_{(u,v)}}: 
L^p_k(f^*_{(u,v)}TV; h)\rightarrow L^p_{k-1}(\wedge^{0,1}(f^*_{(u,v)}TV)). $$

\begin{lemma}
Under our assumption that all critical points $c_i$ of $H$ are non-degenerate in the sense of Floer homology, $L^1_{(u,v)}$ is a Fredholm operator.

\end{lemma}

In general we don't expect that $L^1_{(u,v)}$ is surjective, 
even for a generic choice of $(J, H)$. Failure of the transversality by only 
perturbing the parameter $(J,H)$ has been considered
as a major difficulty in Floer homology and quantum cohomology. 
As we mentioned in the introduction of this paper, this difficulty 
had been overcome through the work [FO], [LiT] and [LiuT1]. 
Following the method we developed in [LiuT1], we define 
$K=K_f=coker L^1_{(0,0)}$, then 
$$L^1_{(0,0)}\oplus  E: L^p_f(f^*TV; h)\oplus K\rightarrow  
L^p_{k-1}(f^*TV))$$ is surjective, where $E$ is the inclusion. 
We can actually choose a modified $K$, 
still denoted as $K$, such that 
$L^1_{(0,0)}\oplus  E$ is surjective and that $K$ has the property described 
in the beginning of last section. As in there we extend $K$ to ${\tilde K}$, 
whose elements are
$T^{N_P}$-invariant section of ${\cal L}$ over 
${\cal F}{\tilde U}_\epsilon(f^e; {\bf H}).$ 
Consider ${\tilde K}$ as a ``bundle " over  
${\cal F}{\tilde U}_\epsilon(f^e; {\bf H})$ with
bundle projection $\pi_{\tilde K}$. Now we define 
${\bar{\ptl}}_{J,H}^{\tilde K}: {\tilde K}\rightarrow {\tilde {\cal L}}$ 
given by
$$\nu\rightarrow {\bar{\ptl}}_{J,H}(\pi_{\tilde K}(\nu))+\nu.$$
In terms of above local coordinate chart and local trivialization,  each element ${\tilde \nu \in {\tilde K}}$ gives
rise to a  map :
$$V_\epsilon\rightarrow  \bigcup_{(u,v)\in \Lambda_\delta}{\cal L}^p_{k-1}
(\wedge^{0,1}(f^*_{(u,v)}TV))$$ and
${\bar{\ptl}}^{\tilde K}_{J,H}$ becomes a  
function $F=\cup_{(u,v)\in \Lambda_\delta} F_{(u,v)}, $ with 
$$F_{(u,v)}=F^1_{(u,v)}\oplus{\tilde E}: V^{(u,v)}_\epsilon
\oplus{\tilde K}\rightarrow
{\cal L}^p_{k-1}(\wedge^{0,1}(f^*_{(u,v)}TV))$$ given by 
$(\xi,\nu)\rightarrow F^1_{(u,v)}+\nu(\xi).$
Clearly, $(DF_{(u,v)})_{(0,0)}=L^1_{(u,v)}+E_{(u,v)}$ where 
$E_{(u,v)}:K\rightarrow
L^p_{k-1}(\wedge^{0,1}(f^*_{(u,v)}TV))$ is given by 
$e\rightarrow{\tilde e}|_{f_{(u,v)}}.$ Let 
$L_{(u,v)}=L^1_{(u,v)}+E_{(u,v)}$. We know that $L_{(0,0)}$ is surjective. 
We want to prove that when $\delta$ is small enough, 
for any $(u,v)\in \Lambda_\delta$, $L_{(u,v)}$ is also surjective. 
In fact in order
to do the gluing, we need somewhat  more. We need to prove that when 
$\|(u,v)\|$ is 
small enough, $L_{(u,v)}$ has a uniformly right inverse with respect to some 
suitable exponential weighted norm on the domain and range of $L_{(u, v)}$. 
We will only consider the following simplest case for defining these norms, 
since this case already contains essential points of general case. 

We assume that $f=f^P\cup f^B$ of two components with $(P, d_1)\cup(B, d_2)$ 
of double points $d_1=d_2.$ Let $z_1, z_2$ be the ends of $P$ and $y_1, y_2$ be
the marked points of $B$. Identify $D_{\epsilon_i}(d_i)$ with ${\bf R}^1\times
S^1=\{(s_i,\theta_i)\}, i=1,2 $ with $d_i$  corresponding to $s_i=+\infty. $ 
 Those local deformations of $f$ coming from only moving double point $d$ along the central circle of $P$ do not play any role in the following definitions. For the reason of simplicity, 
we omit them here. Therefore the local deformation $f$ can be described by a single complex parameter $t\in D_\delta.$ Given $\xi\in L^p_k(f^*_tTV; h)$, we define ${\tilde \xi}^0=\int _{c_{t}}\xi|_{c_t}d\theta,$ where $c_t=\{s_i=-\log |t|\}$ is the central circle 
 of $\Sigma_{f_t}.$

Note that when $s_i>-\log |t|-1$, both $f^P(s_1, \theta_1)$ and $f^B(s_2, \theta_2)$ are 
just $f(d)$. Hence the above definition of ${\tilde \xi}^0$ makes sense and ${\tilde \xi}^0\in T_{f(d)}V.$ We may think ${\tilde \xi}^0$ as a vector field along $f_t$ with $s_i>-\log |t_{0}|$ for some fixed $t_{0}.$ Multiple ${\tilde \xi}^0$ with a fixed cut-off function in $\Sigma_t$, we extend ${\tilde\xi}^0$ to an element $\xi^0\in L^p_k(f^*_tTV, h).$ Let $\xi^1=\xi-\xi^0$. 

Now for  $\xi\in L^p_k(f^*_tTV, h)$, $\eta\in L^p_{k-1}(\wedge^{0,1}(f^*_tTV))$, we define 
$\|\xi\|_{k,p;\mu}=\|e^{\mu s}\xi^1\|_{k,p}+|\xi^0|$ and
$\|\eta\|_{k-1, p; \mu}=\|e^{\mu s}\eta\|_{k-1, p},$ where
$0< \mu< 2\pi$ is fixed and $|\xi^0|=|{\tilde \xi}^0|$ is the Euclidean norm
of $\xi^0\in T_{f(d)}V.$
For any element $(\xi, e)\in L^p_k(f^*_{(u,v)}TV, h)\oplus K$, we define
$\|(\xi,e)\|_{k,p;\mu}=\|\xi\|_{k,p;\mu}+|e|.$

It is proved in [LiuT1] that
\begin{pro}

When $\|(u,v)\|$ is small enough, under these $\mu$-exponential weighted 
norms
$$L_{(u,v)}:L^p_k(f^*_{(u,v)}TV, h)\oplus K\rightarrow L^p_{k-1}(\wedge^{0,1}
(f^*_{(u,v)}TV))$$ has a uniform right inverse $G_{(u,v)}$ in the sense that 
there exists a constant $C_1=C_1(f),$ only depends on $f=f_{(0,0)}$ such
that for any $\eta\in L^p_{k-1}(\wedge^{0,1}(f^*_{(u,v)}TV))$
$$\| G_{(u,v)}(\eta)\|_{k,p;\mu}\leq C_{1}\|\eta\|_{k-1,p;\mu}.$$

\end{pro}
\begin{cor}

$L_{(u,v)}$ is surjective when $\|(u,v)\|$ is small.

\end{cor}

We will use $L^p_{k,\mu}(f^*_{(u,v)}TV,h)$ and 
$L^p_{k-1, \mu}(\wedge^{0,1}(f^*_{(u,v)}TV))$ to 
denote the corresponding spaces equipped with the $\mu$-exponential weighted
norms. 

\subsection{Gluing} 

Now a direct computation shows that the local deformation  $f_{(u, v)}$ is an asymptotic solution of ${\bar\partial}_{J,H}g=0$ when $f$ is a stable 
$(J,H)$-map. 
More precisely, we have 
\begin{lemma} \label{3.7} 
$$\lim_{(u, v)\rightarrow 0}\|{\bar\partial}_{J,H}f_{(u, v)}\|_{k-1,p;\mu}=0.$$ 
\end{lemma} 
 
To do gluing, we also need an estimate on the second order term $Q_{(u, v)}$ 
in the Taylor expansion of $F_{(u, v)}:$ 
$$V_{\epsilon}^{(u, v)}\subset L^p_{k,\mu}(f^{*}_{(u, v)}TV,h)\rightarrow L^p_{k-1,\mu}(\wedge 
^{0,1}(f^{*}_{(u, v)}TV)),$$ 
where $Q_{(u, v)}$ is defined by 
$$F_{(u, v)}(\xi)=F_{(u, v)}(0)+L_{(u, v)}(\xi)+Q_{(u, v)}(\xi).$$ 
\begin{lemma} \label{3.8} 
There exists a constant $C_2=C_2(f)$ only depending on $f$ such that for any 
$\xi_{(u, v)},\,\eta_{(u, v)}\in L^p_{k,\mu}(f^{*}_{(u, v)}TV,h),$ 
$$(i)\quad \|Q(\xi_{(u, v)})\|_{k-1,p;\mu}\leq C_2\|\xi_{(u, v)}\|_{\infty} 
\|\xi\|_{k,p;\mu};$$ 
\begin{eqnarray*} 
(ii)&\|Q(\xi_{(u, v)})-Q(\eta_{(u, v)})\|_{k-1,p;\mu}&\\ 
&\leq C_{2}(\|\xi_{(u, v)}\|_{ 
k,p;\mu}+\|\eta_{(u, v)}\|_{k,p;\mu})\|\xi_{(u, v)}-\eta_{(u, v)}\|_{k,p;\mu}.& 
\end{eqnarray*} 
\end{lemma} 
\proof 
 
The corresponding statement was proved in [F1] when $k=1$, and $1-\frac{2}{p}> 
0.$ The general case here follows from that by a direct induction argument. 
\QED 
\begin{lemma}{(Picard method)} 
Assume that a smooth map $F: X\rightarrow Y$ from Banach spaces $(X,\|\cdot\|)$ 
to $Y$ has a Taylor expansion 
$$F(\xi)=F(0)+DF(0)\xi +Q(\xi)$$ 
such that $DF(0)$ has a finite dimensional kernel and a right inverse 
$G$ satisfying 
$$\|GQ(\xi)-GQ(\eta)\|\leq C(\|\xi\|+\|\eta\|)\|\xi-\eta\|$$ 
for some constant $C$. Let $\delta_1=\frac{1}{8C}.$ If $\|G\circ F(0)\| 
\leq\frac{\delta_1}{2}$, then the zero set of $f$ in $B_{\delta_1}= 
\{\xi,\,|\,\|\xi\|<\delta_1\}$ is a smooth manifold of dimension equal 
to the dimension of $ker DF(0).$ In fact, if 
$$K_{\delta_1}=\{\xi\,|\xi\in ker DF(0),\,\|\xi\|<\delta_1\}$$ 
and $K^{\perp}=G(Y),$ then there exists a smooth function 
$$\phi:K_{\delta_1}\rightarrow K^{\perp}$$ 
such that $F(\xi+\phi(\xi))=0$ and all zeros of $f$ in $B_{\delta_1}$ 
are of the form $\xi+\phi(\xi).$ 
\end{lemma} 
The proof of this Lemma is an elementary application of Banach's fixed point 
theorem (see [F1]).  Now we apply the Picard method above to our case with 
$$X=V^{(u,v)}_{\epsilon,\mu}\oplus {\tilde K}_\epsilon\hookrightarrow
L^p_{k,\mu}(f^*_{(u,v)}TV,h)\oplus{\tilde K},$$ 
$$Y=L^p_{k-1, \mu}(\wedge^{0,1}(f^*_{(u,v)}TV)),\, \mbox{and} \, F=F_{(u,v)}.$$
We have 

\begin{theorem}\label{gluing}

When $\epsilon$ is small enough, the solution set 
${\cal FM}^{K,(u,v)}_{\epsilon,\mu}(f)$ of the equation $F_{(u,v)}=0$ in
$V^{(u,v)}_{\epsilon,\mu}\times {\tilde K}_\epsilon$ is a smooth manifold of
dimension $Ind (c_+)-Ind (c_-)+2c_1(A)+1+ r-2n_\alpha-2n_t-n_\theta-n_\tau,$ 
where $r=\dim K$ and $n_\alpha, n_t, n_\theta$ and $ n_\tau$ are the numbers of zero
components in $\alpha, \theta, t\, \mbox{and}\, \tau$ respectively. Here 
$u=(\alpha, \theta), v=(t,\tau).$ Moreover $DF_{(u,v)}$ is surjective along the
zero set ${\cal FM}^{K,(u,v)}_{\epsilon,\mu}(f)$. 

\end{theorem}

Note that the last statement follows from the fact that $F_{(u,v)}$ induces
a local differentiable embedding:
$$F_{(u,v)}\oplus\pi_N: V^{(u,v)}_{\epsilon, \mu}\oplus K_\epsilon\rightarrow
L^p_{k-1, \mu}(\wedge^{0,1}(f^*_{(u,v)}TV))\oplus N^{(u,v)}, $$ 
where $N^{(u,v)}=N$ is the  kernel of $L_{(u,v)}$ and $\pi_N$ is  the projection with respect to the orthogonal decomposition
$L^p_{k,\mu}(f^*_{(u,v)}TV,h)\oplus K=N\oplus N^{\perp}.$ Here the $L^2$-inner
product in the above decomposition is defined in terms of the ``standard metric"
on $\Sigma_{(u,v)}$ induced from $\Sigma_{(0,0)}$ through gluing.

The exponential weight $\mu$-norm and the usual $(k,p)$-norm on 
$L^p_{k,\mu}(f^*_{(u,v)}TV,h)$ are not uniformly equivalent with respect
to $(u,v)$. However, it is proved in [LiuT1] that

\begin{theorem}\label{surj} When $\epsilon'<<\epsilon, $ the solution set
$\{ F_{(u,v)}(\xi)=0\}$ in $V^{(u,v)}_{\epsilon'}\times {\tilde K}_{\epsilon'}$ is 
contained in ${\cal F}{\tilde{\cal M}}^{K,(u,v)}_{\epsilon, \mu}(f).$ 

\end{theorem}
 In other words, as far as the solution set are concerned,
the exponential weight $\mu$-norm and the usual $(k,p)$-norm are equivalent. 
 
Because of this, we can reformulate Theorem\ref{gluing}.

\begin{theorem}

The solution set 
${\cal F}{\tilde{\cal M}}^{K,(u,v)}_{\epsilon}(f)$  in $V^{(u,v)}_\epsilon
\times{\tilde K}_\epsilon $ is a smooth manifold of $\dim Ind (c_+)-Ind (c_-)
+2c_1(A)+1+r-2(n_\alpha+n_t)-(n_\theta+n_\tau).$

\end{theorem}

Let $N^{(u,v)}_{\epsilon}$ be the $\epsilon$-ball in $N^{(u,v)}$
centered at origin. Let $N_\epsilon=\cup_{(u,v)\in\Lambda_\delta} N_\epsilon
^{(u,v)}\cong \Lambda_\epsilon\times N^{(0,0)}_\epsilon.$ There is a diffeomorphism
of $T^{(u,v)}:N_\epsilon^{(u,v)}\rightarrow   
{\cal F}{\tilde{\cal M}}^{K(u,v)}_{\epsilon}(f)$ as described in Picard method. 

Let $$T=\bigcup_{(u,v)\in\Lambda_\delta}T^{(u,v)}: N_\epsilon\rightarrow
{\cal F}{\tilde{\cal M}}^{K}_{\epsilon}(f)=\bigcup_{(u,v)\in \Lambda_\delta}
{\cal F}{\tilde{\cal M}}^{K(u,v)}_{\epsilon}(f)$$ be the induced ( continuous )
identification. We give ${\cal F}{\tilde{\cal M}}^{K}_{\epsilon}(f)$
the (cornered) smooth structure of $N_\epsilon$ induced from $T$. 

\begin{theorem}\label{smooth}
${\cal F}{\tilde{\cal M}}^{K}_{\epsilon}(f)$ is a cornered smooth manifold
of dimension $Ind (c_+)-ind (c_-)+2c_1(A)+r+1$.  The smooth structure of
${\cal F}{\tilde{\cal M}}^{K}_{\epsilon}(f)$ is induced from 
${\cal F}{\tilde{\cal M}}^{D_f, K}_{\epsilon}(f)\times \Lambda_\delta$ under
the gluing map $T$. 

\end{theorem}

We now come to a $T^{N_P}$-invariant version of above theorem.

 Let 
${\cal F}{\tilde U}_\epsilon(f^e,{\bf H})$ be a $T^{N_P}$-invariant
uniformizer near $f$. Then 
$${\cal F}{\tilde U}_\epsilon(f^e,{\bf H})=\bigcup_I{\cal F}{\tilde U}_
{\epsilon_I}(f_I,{\bf H}), I\in ({\bf Z}/M{\bf Z})^{N_P(D_f)},$$
where $M$ is a  fixed large  integer. Note that the  $T^{N_P}$-orbit of 
$f$ is contained in ${\cal F}{\tilde U}_\epsilon(f^e,{\bf H})$. In
particular, for each $I$, there exists a $\phi_I\in T^{N_P}$ such that
$f_I=f*\phi_I.$ Now assume that we have chosen a $T^{N_P}$-equivariant
extension ${\tilde K}$ of $K$ over 
${\cal F}{\tilde U}_\epsilon(f^e,{\bf H})$. By exploring the naturality
of all relevent construction, it is easy to see that if $K$-perturbed
operator ${\bar\ptl}_{J,H}^{\tilde K}$ is a transversal section in 
${\cal F}{\tilde U}_{\epsilon_0}(f,{\bf H})$ for some sufficiently small
$\epsilon_0$, so is it on all ${\cal F}{\tilde U}_{\epsilon_I}(f_I,{\bf H})$,
$I\in ({\bf Z}/M{\bf Z})^{N_P}$ for $\epsilon_I$ small enough. Because of the 
compactness of the $T^N_{P}$-orbit of $f,$ we may choose
$M$ to be large enough so that 
$\bigcup_I{\cal F}{\tilde U}_{\epsilon_I}(f_I,{\bf H})$ already covers the
orbit of $f$. We will still use 
${\cal F}{\tilde U}_\epsilon(f^e,{\bf H})$ to denote this union. Now we have

\begin{theorem}

When $M$ is large enough and $\epsilon$ is small enough, the perturbed
${\bar\ptl}_{J,H}^{\tilde  K}$-operator is a transversal $T^{N_P}$-equivariant 
section of ${\tilde{\cal L}}$ over ${\cal F}{\tilde U}_\epsilon(f^e,{\bf H})
\times {\tilde K}_\epsilon.$ The solution set 
${\cal F}{\tilde{\cal M}}^{\tilde K}_\epsilon(f^e)$ of ${\bar\ptl}^{\tilde K}_{J,H}
\xi=0$ is a cornered smooth manifold of dimension
$Ind (c_+)-Ind (c_-)+2c_1(A)+r+1.$ The cornered smooth structure of 
${\cal F}{\tilde{\cal M}}^{\tilde K}_\epsilon(f^e)$ is obtained from the
corresponding (cornered) smooth structure of 
${\cal F}{\tilde{\cal M}}^{\tilde K, D_f}_\epsilon(f^e)\times \Lambda_\delta$
under the gluing map $T$.

\end{theorem}
 For any $\nu\in {\tilde K}_\epsilon$, let 
 ${\cal F}{\tilde{\cal M}}^{\nu}_\epsilon(f^e)=\pi^{-1}_2(\nu)$, where
 $\pi_2:{\cal F}{\tilde U}_\epsilon(f^e, {\bf H})\times{\tilde K}_\epsilon
 \rightarrow {\tilde K}_\epsilon$ is the projection. We have
 
\begin{theorem}

For a generic choice of $\nu\in {\tilde K}_\epsilon$, 
${\cal F}{\tilde{\cal M}}^{\nu}_\epsilon(f^e)$ is a 
(cornered) smooth manifold of dimension $Ind (c_+)-Ind (c_-)+2c_1(A)+1$.
There is a continuous $T^{N_P}$-action on 
${\cal F}{\tilde{\cal M}}^{\nu}_\epsilon(f^e)$, which is smooth on
each of its strata. 

\end{theorem}

\section{$T^{N_P}$-invariant Virtual Moduli Cycles}

In this section, we will globalize these $T^{N_P}$-invariant local moduli
spaces ${\cal F}{\tilde{\cal M}}_\epsilon^{\nu_f}(f^e)$ to get a 
$T^{N_P}$-invariant virtual moduli cycles.

Note that we may cover the moduli space ${\cal FM}(J,H,A)$ by using
$$\bigcup_{\la f\ra}{\cal F}U_{\epsilon_f}(f^e; {\bf H}_f), \la f\ra\in
{\cal FM}(J, H, A), $$ where 
$${\cal F}U_{\epsilon_f}(f^e,{\bf H}_f)=\pi_f({\cal F}{\tilde U}_{\epsilon_f}
(f^e, {\bf H}_f))$$ is the  $\pi_f$-image of the $T^{N_P}$-invariant uniformizer
${\cal F}{\tilde U}_{\epsilon_f}(f^e,{\bf H}).$ We may assume that there 
exist finite many $f_i$'s, say, $i=1, \cdots, q$, such that (i) $\bigcup_{1\leq i
\leq q}{\cal F}U_{\epsilon_i}(f^e_i; {\bf H}_i)$ already cover ${\cal FM}(J,H,A)$;
(ii) perturbed ${\bar\ptl}^{{\tilde K}_i}_{J,H}$-operator is transversal to zero section  on
${\cal F}{\tilde U}_{\epsilon_i}(f^e_i;{\bf H}_i)\times {\tilde K}_i$ and
(iii) ${\tilde K}_i$ vanishes near the boundary of ${\cal F}{\tilde U}_{\epsilon_i}
(f^e_i, {\bf H}_i).$

We now use $W_i$ to denote ${\cal F}U_{\epsilon_i}(f_i; {\bf H}_i)$ with 
a $T^{N_P}$-invariant uniformizer ${\tilde W}_i$ and cover group $\Gamma_i$.
Let ${\tilde {\cal L}}_i$ be the corresponding bundle over ${\tilde W}_i$. Note
that the $\Gamma_i$-action commutes with the action of $T^{N_P}$. This
implies that the $T^{N_P}$-action descends to $W_i\hookrightarrow {\cal FB}(A).$

In order to globalize these perturbed moduli spaces ${\cal F}{\tilde {\cal M}}^{\nu_i}
(f^e_i)$, we need to know how these perturbation terms $\nu_i$ change from
${\tilde W}_i$ to ${\tilde W}_j$.  The idea now is to use a fiber product
construction of the covering $({\tilde W}_i,\pi_i)$ as a replacement of 
``intersections" of ${\tilde W}_i$'s. 

For this purpose, let ${\tilde{\cal N}}$ be the nerve of the  covering $W=\{W_i,
1\leq i\leq q\}.$ We will use elements of ${\tilde {\cal N}}$ as  indices,
In other words. we define the following multi-indices set
$${\cal N}=\{I=(i_1, \cdots, i_n)\,|i_1<i_2<\cdots<i_n, W_{i_1}\cap
W_{i_2}\cap\cdots W_{i_n}\not =0\}.$$

We define the length $l(I)=n,$ for $I=(i_1, \cdots, i_n).$ 
Note that ${\cal N}$ has an obvious partial order induced by inclusion. For any
$I=(i_1,\cdots, i_n)\in {\cal N}$, we will use ${W}_I$ to denote
$W_{i_1}\cap W_{i_2}\cdots \cap W_{i_n}$, and ${{\cal L}_I}=
{\cal L}|_{W_I}.$

 There are $n$  uniformizing systems 
$$({\widetilde{\cal L}}_{i_1,\cdots,{\widehat{i_k}},\cdots,i_n}, 
{\widetilde W}_{i_1,\cdots, 
\widehat{i_k},\cdots, i_n}; \pi_{i_1,\cdots, \widehat{i_k}, \cdots, 
i_n})$$ 
of $$({ {{\cal L}}}_{ I}, {  W}_{I }), 
 $$ 
with covering group $\Gamma_{i_k},$ induced from 
$$\pi _{i_k}: 
({\widetilde{\cal L}}_{i_k}, {\widetilde W}_{i_k})\rightarrow 
({ {{\cal L}}}_{i_k},{ W}_{i_k}),$$ where $${\widetilde W}_{i_1,\cdots,\widehat 
{i_k},\cdots,i_n}=(\pi_{i_k})^{-1}({  W}_{I })$$ and 
$${\widetilde{\cal L}}_{i_1,\cdots, \widehat{i_k},\cdots, i_n}= 
{\widetilde{\cal L}}|_{{\widetilde W}_{i_1,\cdots, \widehat{i_k},\cdots, i_n}}.$$

We want to construct the pull-back of these  morphisms, which is  denoted by 
 $$\pi_{ I}: 
  ({\widetilde {\cal L}}_{ I}^{\Gamma_{ I}}, 
  {\widetilde W}_{ I}^{\Gamma_{ I}})\rightarrow 
  ({ {{\cal L}}}_{    I}, {  W}_{ I})$$ with covering group 
  $$\Gamma_{ I}=\Gamma_{i_1}\times \cdots\times\Gamma_{i_n}.$$ 
  We  define first  
\begin{eqnarray*} 
{\widetilde W}^{\Gamma_I}_I  &=&
\left \{u\,|\,u\in\prod_{k=1}^n {\widetilde W}_{i_k}, 
 \begin{array}{ccc} 
        \pi_{i_k}(u_k)& \in & {W}_{I},\\ 
                          \pi_{i_k}(u_k)& = & \pi_{i_l}(u_l) 
                                                \end{array}\right\}. 
\end{eqnarray*} 
  We define  $\pi _{ I}$  to be  the composition of $\prod_{k=1} 
^n \pi_{i_k}$ restricting to ${\widetilde W}_I^{\Gamma_I}$ with $\triangle^{-1}_n$ of 
the inverse of n-fold diagonal. If $J=(j_1,\cdots, j_m)\subseteq I=(i_1,\cdots, 
i_n)$, there exists an obvious projection map 
$$\pi^I_J: {\widetilde W}^{\Gamma_I}_I 
\rightarrow {\widetilde W}^{\Gamma_{J}}_J$$ 
 induced from the 
corresponding projection $\prod_{i_k\in I}{\widetilde W}_{i_k}$ to 
$\prod_{j_l\in J}{\widetilde W}_{j_l}$ 
such that $\pi_J\circ \pi^I_J= \iota^I_J\circ\pi_I$ when restricted to the inverse image of $\pi^I_J$, where 
$\iota^I_J$ is the inclusion ${ W}_I\hookrightarrow{  
W}_J.$ 
 
All the above constructions can be directly extended to bundle case and we get a 
system of bundles 
$\{p_I:{\widetilde{\cal L}}_I^{\Gamma_I}\rightarrow {\widetilde W}_I^{\Gamma_I}\}$, 
$I\in {\cal N}.$ 
 
Note that for any fixed $I$ with $l(I)>1,$ ${\widetilde W}_I^{\Gamma_I}$ is not a 
( stratified) smooth manifold in general but rather a ( stratified) smooth 
variety, 
i.e., locally it is a finite union of ( stratified) smooth manifold. In fact 
for $u\in {\widetilde W}^{\Gamma_I}_I$ with $u=(u_1,\cdots, u_n)$, ${\bar u}=\pi_{i_k}(u_k)$, 
we can choose  an open neighborhood ${ U}$ of ${\bar u}$ in ${  
W}_I$ and consider   the inverse image 
${\widetilde U}_k=\pi^{-1}_{i_k}({U})$ of $u_k$ in ${\widetilde W}_{i_k}.$ 
When $U$ is small enough, there exist $(n-1)$ 
equivalence maps $\lambda_k:{\widetilde U}_1\rightarrow {\widetilde U}_k,$ 
$k=2, \cdots, n.$ 
Composing with the actions of automorphism group   $\Gamma_{u_k}$ of 
 ${\widetilde U}_k,$   we get $\prod_{i=2}^n|\Gamma_{u_i}|$ equivalence 
maps: 
$$\phi_k\lambda_k:{\widetilde U}_1\rightarrow {\widetilde U}_k, 
\quad k=2, \cdots,n,\quad \phi_k\in \Gamma_{u_k}.$$ 
Clearly $u=(u_1, u_2, \cdots,u_n)\in\prod_{k=1}^n {\widetilde U}_k$ is contained in 
${\widetilde W}^{\Gamma_I}_I$ if and 
only if $u_k =\phi_k\lambda_k (u_1)$ for some $\phi_k\in \Gamma_{u_k}$, 
 $k>1.$  Thus, in general we can identify a neighborhood 
 ${\widetilde U}$ of $u$ in ${\widetilde W}^{\Gamma_{I}}_{I}$ with an 
 union of  $\prod_{i\not =j}^n|\Gamma_{u_i}|$  copies of 
 ${\widetilde U}_{j}.$ 

Note that  the action of $\Gamma_{i}$ commutes with $T^{N_{P}}$-action.  That implies that all our constructions here are $T^{N_{P}}$-equivariant.  

We summarize up the above discussion in the following lemma. 
\begin{lemma} 
 
There exists a pull-back 
$$\pi_I : 
({\widetilde{\cal L}}^{\Gamma_I}_I,{\widetilde W}_I^{\Gamma_I}) 
\rightarrow ({ {{\cal L}}}_I, {  W}_I)$$ of the $n$  uniformizing systems 
$$\pi_{i_1,\cdots,\widehat{i_k},\cdots, i_n}  : 
({\widetilde{\cal L}}_{i_1,\cdots, \widehat{i_k},\cdots, i_n}, 
{\widetilde W}_{i_1,\cdots,\widehat{i_k}, 
\cdots,i_n})\rightarrow  ({ {{\cal L}}}_I, { W}_I)$$ 
in the category of (stratified) smooth varieties  with the automorphism 
group $\Gamma_I.$ 
 
For 
any $J\subset I$, there exists a projection 
$$\pi_J^I: 
({\widetilde{\cal L}}^{\Gamma_I}_I, {\widetilde W}_I^{\Gamma_I}) 
\rightarrow ({\cal L}_J^{\Gamma_J}, { W}_J^{\Gamma_J}),$$ 
whose generic fiber contains $\frac {|\Gamma_{I}|} {|\Gamma_{J}|}$ points, 
where $\frac {|\Gamma_{I}|} {|\Gamma_{J}|} =\prod_{i_{k}\in I\setminus J} |\Gamma_{i_k}|$ 
. It satisfies  the relation that $\pi_J\circ\pi_J^I={\iota}^I_J\circ \pi_I$ 
for each $I\in {\cal N},$ when restricted to the inverse image of $\pi ^I_J.$
 Moreover, all constructions can be done in a $T^{N_{P}}$-equivariant 
manner. 
\end{lemma} 

\begin{lemma} 
There exists an open covering $\{{  V}_{I }\},$   $ I\in {\cal N}$ of 
${{{\cal  FM}}}(J,H; A)$ 
such that 
 
(i)$V_ I\subset W_I$, for all 
 $I\in {\cal N};$ 
 
(ii) $Cl({  V}_{I_1})\cap Cl({  V}_{I_2})\not =\emptyset$ only if  
$I_1 < I_2,$  or $I_2 < I_1.$

 Moreover, all $V_{I}$ can be chosen to be $T^{N_{P}}$-invariant. 
\end{lemma} 
\proof 
 
We may assume that there exist open sets ${ W}_i^1\subset\subset 
{  W}_i$, $i=1, \cdots, q$ such that $\{{  W}_i^1, i=1,\cdots, 
q\}$ already forms a covering of ${{{\cal FM}}}(J,H, A)$. 
For each fixed $i$ we can find pairs of open sets $  {W}_i^j 
\subset\subset {  U}_i^j$, $j=1,\cdots,q-1$ such that 
$${ W}^1_i\subset\subset {  U}_i^1\subset\subset 
{ W}_i^2\subset\subset { U}_i^2\cdots\subset\subset 
{  W}_i^q={  W}_i.$$ 
Now define 
$${  V}_{i_1,\cdots,i_n}={  W}_{i_1}^n\cap{ W}^n_{i_2} 
\cdots\cap {  W}_{i_n}^n\setminus (\cup_{J\in {\cal N}_{n+1}} 
Cl({  U}_{j_1}^n)\cap Cl({  U}^n_{j_2})\cdots\cap Cl({  
U}^n_{j_{n+1}})),$$ 
where $J=(j_1,\cdots, j_{n+1}).$ 
 
Clearly the family $\{{  V}_{i_1,\cdots,i_n}, (i_1,\cdots, i_n)\in{\cal 
N}\}$ so constructed satisfies the conditions in the lemma. 

To get a $T^{N_{P}}$-invariant construction, we only need to run through above construction for
${\cal  GM}(J,H; A)$ 
and define $V_{I}$ as the lifting of the corresponding sets constructed by using  stable
${\cal G}$-maps.
\QED 
 
Now we define 
$${\widetilde V}_I=(\pi_I)^{-1}({ V}_I),\quad {\widetilde E}_I=(\pi_I)^{-1} 
({ {{\cal L}}}_{I}).$$ 
Then the bundle $({\widetilde E}_I, {\widetilde V}_I)$ are still a pair of 
(stratified ) smooth varieties and for 
any $J\subset I$ the projection $\pi^I_J$ still can be defined when 
restricted to $(\pi^I_J)^{-1}({\widetilde E}_J, {\widetilde V}_J)\cap ({\widetilde E}_I, 
{\widetilde V}_I).$ 
Since locally ${\widetilde V}_I$ is a finite union of its ( stratified smooth) components, 
we will say 
a continuous section $S_I:{\widetilde V}_I\rightarrow {\widetilde E}_I$ 
to be smooth 
if locally, $S_I$ 
restricted to any of those  components is  (stratified )smooth. For a smooth section 
$S_I$, we say that  $S_I$ is  
transversal to zero section if locally, $S_I$ restricted to any 
 of the smooth components of ${\widetilde V}_I$ is transversal to zero section . 
 
 Now let $({\widetilde E},{\widetilde V})$ be the collection 
 $\{({\widetilde E}_I,{\widetilde V}_I),\pi^I_J; J\subset I\in {\cal N}\}.$ 
 We can define a global 
 section $S=\{S_I;I\in {\cal N}\}$ of such a system by requiring the 
 obvious compatibility condition: 
 $$(\pi^I_J)^{*}S_J=S_I|_{{\pi^I_J}^{-1}({\widetilde V}_J)}.$$ 
 $S$ is said to be transversal to zero section if each $S_I$ is. 
 
 Now the section ${\overline{\partial}}_{J,H}:{ W}\rightarrow 
 { {{\cal L}}}$ gives rise to a global section of the bundle system 
 $({\widetilde E},{\widetilde V})$ in an obvious way. Our goal now is to perturb 
 ${\overline{\partial}} 
 _{J,H}$ to get a global transversal section. To this end, we need to know 
 how an element $\nu_i\in K_{i}$ can be interpreted as a global section 
 of $({\widetilde E},{\widetilde V})$ first.  

 \begin{lemma} 
Each ${\tilde \nu_i}\in {\widetilde K}_i$ gives rise to a global section, 
denoted by same notation ${\tilde{\nu}}_i= 
\{({\tilde{\nu}}_i)_{I}; I\in {\cal N}\}$, 
of the system $({\widetilde E}, {\widetilde V}),$
which is $T^{N_{P}}$-equivariant. 
\end{lemma} 
\proof 

By multiplying with some $\Gamma_i$-equivariant cut-off 
 function 
 $\beta_i$, we  may assume that the support of each element ${\tilde \nu_i}$ is 
 contained 
 in ${\widetilde W}_i^1=\pi_i^{-1}({ W}_i^1)$ and that $\{{ U}_i^0;i=1,\cdots, 
 q\}$ already forms a covering of ${{{\cal FM}}}(J,H,A)$, 
 where ${\widetilde U}_i^0=\{u\,|\,u\in {\widetilde W}_i,\beta_i(u)>0\}$ and 
 ${ U}_i^0= 
 \pi_i({\widetilde U}_i^0).$   
 Now since each $\nu_i$ vanishes near the boundary of ${\widetilde W}_i$, we may consider 
 it as a global multi-valued section ${\bar\nu}_i$ of ${ {{\cal L}}} 
 \rightarrow { W}$ supported in ${ U}_i^0\subset\subset 
 { W}_i^1.$ 
Let $I\in {\cal N}$ with $i\not\in I$ and consider ${ V}_I.$ 
 
Recall that if $I=\{i_1,\cdots, i_n\}$ then 
$${ V}_I={  W}^n_{i_1}\cap { W}_{i_2}^n\cdots 
\cap { W}_{i_n}^n\setminus \cup_{J\in{\cal N}_{n+1}} 
 Cl({  U}^n_{j_1})\cdots \cap Cl({ U}^n_{j_{n+1}})$$ 
 with $J=(j_1,\cdots, j_{n+1}).$ Since $i\not\in I,$ 
 \begin{eqnarray*} 
 { V}_I&\subseteq &{ W}_{i_1}^n\cap\cdots\cap{  
 W}^n_{i_n}\setminus  
 { W}_{i_1}^n\cap\cdots\cap{  
 W}^n_{i_n}\cap Cl({ W}_i^1)\\ 
 & \subseteq & { W}\setminus Cl({ W}_i^1). 
 \end{eqnarray*} 
Therefore, the intersection $Cl({  U}^0_i)\cap Cl({ V}_I)=\emptyset.$ 
Hence ${\tilde{\nu}}_i|_{{ V}_I}\equiv 0$ for any $I\in {\cal N}$ 
with $i\not\in I.$ We define $({\tilde{\nu}}_i)_I\equiv 0$ if $i\not\in I. $ 
 
Now assume that $i\in I.$ 
 
When $l(I)=1 $ hence $I=\{i\}$, ${\widetilde V}_I={\widetilde V}_i$ and 
$({\tilde {\nu}}_i)_I$ 
is just ${\tilde \nu}_i:{\widetilde W}_i\rightarrow {\widetilde{\cal L}}_i$ restricted 
to ${\widetilde V}_i.$ 
 
 If we denote 
$\{i\}$ by $I_i$, then for any $I$  with $n=l(I)>1$, we have 
$$\pi^I_{I_i}:({\widetilde{\cal L}}^{\Gamma_I}_I,{\widetilde W}_I^{\Gamma_I}) 
\rightarrow 
({\widetilde {\cal L}}_{I_{i}},{\widetilde W}_{I_{i}})\subset 
({\widetilde{\cal L}}_i,{\widetilde W}_i).$$ 
Therefore, $(\pi^I_{I_i})^{\ast}({\tilde{\nu}}_i)_{I_i}$ gives rise to a section 
of ${\widetilde E}_I\rightarrow {\widetilde V}_I,$ denoted by $({\tilde \nu}_i)_I. 
$ Clearly the 
section $({\tilde\nu}_i)_I$, $I\in {\cal N}$ so constructed are compatible 
to each other and yields a well-defined global section ${\tilde \nu}_i 
=\{({\tilde \nu}_i)_I, I\in {\cal N}\}$ of the system $({\widetilde E},{\widetilde V}).$ 

Finally, we note that  $\beta_{i}$ can be chosen to be 
$T^{N_{P}}$-equivariant, 
which implies that  ${\tilde \nu}_{i}$ so constructed is also $T^{N_{P}}$-equivariant.

\QED 
 
Let  
 $K=\oplus_{i=1}^q K_i.$   Consider the system
$$({\widetilde E} \times  {\widetilde K}_{\delta}, 
{\widetilde V}\times {\widetilde K}_{\delta})=\{({\widetilde E}_I
\times {\widetilde K}_{\delta}, {\widetilde V}_I\times 
{\widetilde K}_{\delta}); I\in {\cal N}\}$$ 
of bundles, where ${\widetilde K}_{\delta}$ is a $\delta$-neighborhood 
of zero of ${\widetilde K }$ under the identification of  ${\widetilde K}$ 
and $K$.  We now defined global section 
${\overline{\partial}}_{J,H}^{\tilde K}$  
given by $$({\overline{\partial}}_{J,H}^{{\tilde K}})(u_I, {\tilde \nu})=
{\overline{\partial}}_{J,H}u_{I}+ 
{\tilde{\nu}}_I(u_I)$$ for any 
$(u_I, {\tilde \nu})\in {\widetilde V}_I\times 
{\widetilde K}_{\delta}.$ 

\begin{theorem}  
                      
${\overline{\partial}}_{J,H}^{{\tilde K}}$ is a smooth section of 
 $({\widetilde E}\times {\widetilde K}_{\delta}, {\widetilde V}\times 
{\widetilde K}_{\delta})$, which is transversal to zero section. It follows that when $\delta$ 
is small enough for a generic choice of the perturbation term ${\tilde \nu}\in {\widetilde K}_{\delta}$ 
the section ${\overline{\partial}}_{J,H}^{\nu}:{\widetilde V} 
\rightarrow {\widetilde E}$ is 
transversal to zero section and that the family of perturbed moduli spaces 
$${\widetilde{\cal M}}^{\nu}=
\{{\widetilde{\cal M}}_I^{\nu}=({\overline{\partial}}_{J,H}^{\nu_I} )^{-1}(0); 
I\in {\cal N}\}$$ is compatible in the sense that 
$$\pi^I_J({\widetilde{\cal M}}^{\nu}_I)=
{\widetilde{\cal M}}_J^{\nu}\cap (Im \pi_J^I),\quad J\subset I.$$ 
 
\end{theorem} 
 \proof 
 
${\overline{\partial}}_{J,H}^{{\tilde K}}$ is obviously  (stratified) smooth. 
Since  
$$({\overline{\partial}}_{J,H}^{{\tilde K}_i})|_{{\widetilde U}^0_i} : 
{\widetilde U}^0_i \times 
({\widetilde K}_{i})_{\delta}\rightarrow \pi^{\ast}_1({\widetilde{\cal L}} 
_i|_{{\widetilde U}_i^0})$$ is  transversal to zero section,  so is  
${\overline{\partial}}_{J,H}^{{\tilde K}} $ on $ {\widetilde U}^0_i \times 
{\widetilde K}_{\delta}, $ for any $1\leq i\leq q.$

Now $\{{ U}^0_i; i=1, \cdots, m\}$ already forms a covering of 
${{{\cal FM}}}(J,H; A).$  Outside  the inverse images of $\pi_{i}$ of this covering, ${\overline{\partial}}_{J,H}^{\tilde K}=
{\overline{\partial}}_{J,H }.$ This implies that the zero set of  ${\overline{\partial}}_{J,H}^K$  is contained in the inverse image of the covering. However,   
${\overline{\partial}}_{J,H}^K$ is already transversal to zero section  in 
$(\pi^I_{I_i})^{-1}({\widetilde U}^0_i)\hookrightarrow {\widetilde W}_{ I}$ 
for any $I\in {\cal N},$ where $I_i=\{i\}.$ 
 This proves the transversality for 
${\overline{\partial}}_{J,H}^{{\tilde K}}.$ It follows from implicit function theorem 
applied locally to each  smooth component of ${\widetilde V} 
=\{{\widetilde V}_I;I\in{\cal N}\}$ 
that 
$$({\overline{\partial}}_{J,H}^{{\tilde K}})^{-1}(0)
=\{({\overline{\partial}}_{J,H}^{{\tilde K}} )^{-1}_I 
(0);I\in {\cal N}\}\subset {\widetilde V}\times {\widetilde K}_{\delta}$$ is a family of ``cornered" 
(stratified) smooth subvarieties. 
 
Let $$\pi:({\overline{\partial}}_{J,H}^{{\tilde K}})^{-1}(0)\rightarrow {\widetilde K}_{\delta}$$ be 
the restriction of projection of ${\widetilde V}\times {\widetilde K}$ to ${\widetilde K}$. It is easy to see that 
Smale-Sard theorem is still applicable in this case. We conclude that for `` 
generic" choice of ${\tilde \nu}\in {\tilde  K},$ ${\overline{\partial}}_{J,H}^\nu$ is 
  a transversal section of $({\widetilde E},{\widetilde V}).$ 
 
The compatibility of the family of zero set 
$$\{{\widetilde{\cal M}}^{\nu}_I; I\in {\cal N} 
\}=\{ {\overline{\partial}}_{J,H}^{\nu} )^{-1}_I(0)\}$$ follows from the fact that 
${\overline{\partial}}_{J,H}^\nu$   is a global section of $(\widetilde E,\widetilde V).$ 
 
\QED
 
We can give a canonical orientation for each ${\tilde {\cal M}}^{\nu_I}
_I$ (see [F1], [FH]). Now we get a family of ``singular cells" of ${\cal FB}(A)$, 
given by $\pi_I:{\tilde{\cal M}}^\nu_I\rightarrow{\cal FB}(A)$ for
any $I\in {\cal N}$. If $I\subset J$, then $\pi_I$ and $\pi_J$ are related
by the $\frac{|\Gamma_J|}{|\Gamma_I|}$-folded covering $\pi^J_I$ in the 
overlap
$$(\pi^J_I)^{-1}({\tilde {\cal M}}^\nu_I)\hookrightarrow {\tilde{\cal M}}
^\nu_J.$$ This suggests that the family of the rational ``singular" chains
of ${\cal FB}(A)$, defined by $S^\nu_I=\frac{1}{|\Gamma_I}\pi_I$ are 
compatible each other when restricted to those overlaps above. Therefore 
after been identified over those overlaps,  $\{S^\nu_I, \, I\in {\cal N}\}$ form a
well-defined rational relative ``singular"  cycle of ${\cal FB}(A)$. 

We
formally write this as $$C^\nu=\sum_{I\in{\cal N}}S^\nu_I.$$
 Now all 
$T^{N_P}$-actions carry over to the cycle $C^\nu$. In particular there is an $S^1$-action
on the top strata of $C^\nu$. To state our main theorem, we need to indicate
the dependence on the critical points $c_+, c_-$ and homology class $A$ in our
notations. We will write ${\tilde{\cal M}}^\nu_I(c_-, c_+; A)$ and
$S^\nu_I(c_-, c_+; A)$ etc. To define the boundary of $C^\nu(c_-, c_+; A)$, we
also need to indicate where the fixed marking $x=(x_1, \cdots, x_n)$ is located in
our notations. We have

\begin{theorem}\label{vcyc}

For generic choice of $(J,H, \nu )$, fix any two critical points
$c_-, c_+$ of $H$. There is a rational virtual moduli ( relative ) cycle 
$$C^\nu(c_-, c_+; A)=\sum_{I\in{\cal N}}S^\nu_I(c_-, c_+; A)$$ of
${\cal FB}(A),$ whose dimension is equal to
$ Ind (c_+)-Ind (c_-)+2c_1(A)+1.$  The boundary 
$\ptl C^\nu(c_-, c_+; A) $ is
$$\sum_{c, A_-, A_+}(C^\nu(c_-, c; A_-, x)\times
C^\nu(c, c_+; A_+)\cup C^\nu(c_-, c; A_-)\times C^\nu(c, c_+; A_+, x),$$
$A_{-}+A_{+}=A,$
 whose dimension   is  equal to
 $Ind (c_+)-Ind (c_-)+2c_1(A).$
Moreover there are $S^1$ and $T^2$-actions on $C^\nu(c_-, c_+; A)$ and
its boundary respectively.
\end{theorem}

We note that for any $\phi\in S^1$ and $g\in C^\nu(c_-, c_+;A)$, the action
$\phi$ on $g$ is given by a rotation of the domain $(\Sigma_g,l,x)\cong
(S^1\times{\bf R}; {\tilde l}, {\tilde x})$, which changes the relative
position of the fixed marking $x$ with respect to the marked line $l$. This
implies that the $S^1$-action on $C^\nu(c_-, c_+; A)$ is free. 

\section {GW-invariants and Weinstein Conjecture}\label{gw}

  In this section  we will  prove
Theorem 1 for genus zero case. 

Using the moduli cycle we obtained in last section, we define a Morse
theoretic version of GW-invariants under our main assumption that ${\tilde H}$ has
no closed orbits. 

Let $E_{\cal F}:{\cal FB}(A)\rightarrow V^n.$
$F_l:{\cal FB}(A)\rightarrow{\cal GB}(A)$ given by forgetting the marked
 We now compose the moduli cycle 
$C^\nu(c_-, c_+; A, x)$ in ${\cal FB}(c_-, c_+; A)$ with the evaluation
map $E_{\cal F}$ and define a (relative) moduli cycle 
$E_{\cal F}\circ C^\nu(c_-, c_+; A)$ in $V^n.$ We denote it by 
$C^\nu_x=C^\nu_x(c_-,c_+; A).$

Given $\beta_i\in H_*(V; {\bf Q}), i=1, \cdots, n$, for simplicity, we may
assume that each $\beta_i$ can be represented by a smooth manifold of $V$. We still 
use $\beta_i$ to denote this representative.
Then $\beta=\beta_1\times\cdots\beta_n$ is a cycle in $V^n$. Assume that
the codimension of $\beta$ in $V^n$ is 
$$\dim C^\nu_x(c_-, c_+; A)-1.\eqno(**)$$ Now
consider $C^\nu_x$ as a map from the `` domain'' 
$\frac{1}{|\Gamma|}{\tilde {\cal M}}^\nu=\sum_{I\in{\cal N}}
\frac{1}{|\Gamma_I|}{\tilde {\cal M}}^\nu_I$ to $V^n.$ By perturbing
$\beta_i$ slightly, 
we may assume that $C^\nu_x$ is transversal to $\beta$ and 
$(C^\nu_x)^{-1}(\beta)$
is a compact submanifold of 
$\frac{1}{|\Gamma|}{\tilde {\cal M}}^\nu.$ The
dimension assumption above implies that $\dim (C^\nu_x)^{-1}(\beta)=1. $ 
Since $C^\nu_x$ is $S^1$-equivariant, $(C^\nu_x)^{-1}(\beta)$ also carried an $S^1$-action, which
is free as we mentioned above. This implies that $(C^\nu_x)^{-1}(\beta)/S^1$ is
a finite set. We will use $C^\nu(c_-, c_+; A,\beta)$ to denote this set. There
is an induced orientation on $C^\nu(c_-, c_+; A,\beta)$.

\begin{definition}
Let $(J,H)$ be a generic pair with $H$ being $C^\infty$-close to ${\tilde H}.$ 
Given a homology class  $A\in H_2(V; {\bf Z})$ 
  and  a cycle and a cocycle of
 Morse-Witten complex,  
represented by certain linear combinations $c-_, c_+$ of critical points
 of $H$, choose $\beta_i\in H_*(V; {\bf Q}), i=1,\cdots, n, $
such that (**) above holds.
 We define the Morse theoretical version of GW-invariant $\Phi_{A, J, H}(c_{-}, c_{+})$ by specify its value at all such $\beta$'s.
$$\Phi_{A, J, H}(c_-, c_+, \beta_1,\cdots, \beta_n)\equiv \#(C^\nu(c_-, c_+;
J, H, A, \beta))\in {\bf Q},$$
where $\nu\in K_\epsilon$ is a generic element of $K_\epsilon$. Note that on
the righthand side of above equality, the number is counted with sign.
\end{definition}

Now let $H_\lambda=\lambda\cdot H$, and choose a corresponding generic
$J_\lambda.$ We get a $\lambda$-dependent generic pair $(J_\lambda, H_\lambda).$
Note that all $H_\lambda, \lambda >0,$ has same critical point set. We
can define a $\lambda$-dependent $GW$-invariant 
$\Phi_{A, J_{\lambda}, H_{\lambda}}(c_-, c_+).$

 The key point needed
to prove Theorem 1 is the following invariant property of the GW-invariants.

\begin{theorem}\label{invt}
$$\Phi_{A,J_{\lambda},H_{\lambda}}(c_-, c_+,\beta_1,\cdots, \beta_n)$$ is independent of the
choice of $\lambda, \lambda\in [\epsilon, \lambda_0+1/2),$ when 
$(J_{\lambda},H_{\lambda})$  is generic.
\end{theorem}

\proof

Fix $0<\lambda_-<\lambda_+,$ let $\Lambda=[\lambda_-, \lambda_+]\subset
[\epsilon, \lambda_0+1/2)$ be the interval where the parameter $\lambda$
varies. We can run through everything developed in the previous sections to
incorporate the parameter $\lambda.$ Therefore, we will have 
${\cal FB}_\Lambda(A)$, ${\cal F}{\cal M}_\Lambda(c_-, c_+, J_\lambda,
H_\lambda, A)$ etc. Here for instance, ${\cal FB}_\Lambda(A)=\{(f,\lambda)
\,|\,f\in {\cal FB}(A), \lambda\in \Lambda\}.$ We can similarly define the
virtual (relative) moduli cycle 
$C^\nu_\Lambda(c_-, c_+, J_\lambda, H_\lambda; A)$, $C^\nu_{x,\Lambda}
(c_-, c_+; J_\lambda, H_\lambda,A)$ and $C^\nu_\Lambda(c_-, c_+; J_\lambda,
H_\lambda, A,\beta)$. We may assume that $\nu$ has been chosen is such a way
that at two end points $\lambda_-$ and $\lambda_+$ of $\Lambda$, $\nu_{\lambda_-}$
and $\nu_{\lambda_+}$ are also generic so that $C^{\nu_{\lambda_-}}(c_-, c_+, J_{\lambda_-},
H_{\lambda_-}, A, \beta)$ and $C^{\nu_{\lambda_-}}(c_-, c_+, J_{\lambda_+},  
H_{\lambda_+}, A,\beta)$ are well-defined. Now the crucial step is the following

\begin{lemma}\label{bnd}
When the condition (**) on dimension holds, $C^\nu_\Lambda(c_-, c_+; J_\lambda,
H_\lambda, A,\beta)$ is a one dimensional (relative) virtual moduli cycle. 
It has the boundary
\begin{eqnarray*}
& & \ptl C^\nu_\Lambda(c_-,c_+; J_\lambda, H_\lambda, A,\beta) \\
& = & C^{\nu_{\lambda_-}}
(c_-, c_+; J_{\lambda_-}, H_{\lambda_-}, A,\beta)\\
& \cup & C^{\nu_{\lambda_+}}(c_-, c_+, J_{\lambda_+}, 
H_{\lambda_+}, A, \beta)  \\
& \cup & \{\cup_{\lambda,c}M(c_-, c; H_{\lambda})\times 
C^{\nu_{\lambda}}(c,c_+;J_\lambda,
H_\lambda,A, \beta)\}\\
& \cup & \{\cup_{\lambda,c}C^{\nu_\lambda}(c_-, c;J_\lambda, 
 H_\lambda, A,\beta)\times
M(c, c_+; H_\lambda)\},  
\end{eqnarray*}
where in the third term $c$ runs through all critical points of $H_{\lambda}$ 
such that
$Ind (c)-Ind (c_-)=1$ and in the fourth term $Ind(c_+)-Ind(c)=1$. Here 
$M(c_-,c; H_\lambda)$ is the moduli space of unparametrized gradient lines of
$\nabla H_\lambda$ connecting $c_-$ and $c.$ So is $M(c, c_+; H_\lambda)$ in
a similar way.
\end{lemma}

\proof 

The boundary $\ptl C^\nu_\Lambda(c_-, c_+; J_\lambda,H_\lambda, A,\beta)$
certainly contains the four terms listed in the lemma. We need to prove
that there is no other terms. Consider one of the components in the
boundary $\ptl C^\nu_\Lambda(c_-, c_+; J_\lambda,H_\lambda, A,\beta)$. By an
analogy of Theorem \ref{vcyc}. we may assume that it has the form 
$$\cup_{\lambda_i}  C^{\nu_{\lambda_i}}(c_-, c; A_-)\times C^{\nu_{\lambda_i}}   
(c, c_+; A_+, x),$$
where $A_-+A_+=A$ and $\lambda_i\in \Lambda$ are of finitely many. Now the
evaluation map
$$E=\cup_{\lambda_i}E_{\lambda_i}:\cup_{\lambda_i}C^{\nu_{\lambda_i}}(c_-, c; A_-)
\times C^{\nu_{\lambda_i}}(c, c_+; A_+, x)\rightarrow V^n$$
only involves the second components above. The dimension condition (**)
implies that after quotienting out the $S^1$-action of second factor, $E^{-1}(\beta)$
is already zero dimensional and hence is a finite set for our generic choice
$(J, H,\nu).$ Now for any element $g\in E^{-1}(\beta)/S^1$, there is another 
$S^1$-action acting on the first factor. However the isotropy group $I_g
\hookrightarrow S^1$ is either finite or $S^1$. The first case contradicts
to the finiteness of $E^{-1}(\beta)/S^1$. This implies that any element 
$g$ in first factor $C^{\nu_i}(c_-, c, A_-)$ is $S^1$-invariant. This,
in turn, implies that $A_-=0$, $A_+=A$, and $Ind(c)-Ind(c_-)=1.$ Now
$M(c_-, c; H_{\lambda_i})$ is obtained in ${\cal F}{\cal M}(c_-, c;
J_\lambda,H_\lambda, 0)$ as an isolated compact component   when $Ind(c)-Ind(c_-)=1,$
and it is just the fixed points set of the $S^1$-action. It follows from this
and the vanishing property along $M(c_-, c;H_{\lambda_i})$ of
elements in $K$ that $$C^{\nu_{\lambda_{i}}}(c_-,c_+; A_-)=M(c_-, c; H_{\lambda_i}).$$

Note that in the last two terms of the expression for the boundary 
$$\ptl C^\nu_\Lambda(c_-, c_+; J_{\lambda}, H_{\lambda}, A_{\lambda}),$$  there are
only  finite many $\lambda\in 
\Lambda$ involved for the dimensional reason. 

\QED

We now prove that counting
algebraically, the last two terms have no contribution to the boundary
operator. Recall that we have defined $c_-$ and $c_+$ as two cycles of the
Morse-Witten complex with respect to Morse function $H$ and $-H$
respectively. This implies that when $Ind(c)-Ind(c_-)=1$ or 
$Ind(c_+)-Ind(c)=1$, $\#(M(c_-,c; H))=0,$ $\#(M(c,c_+; H))=0$, 
where both numbers are counted algebraically. Therefore for any fixed 
$\lambda_i$ and fixed critical point $c_-, c$ of $H_{\lambda_i}$ with
$Ind(c)-Ind(c_-)=1$, 
\begin{eqnarray*}
\#\{M(c_-,c; H_{\lambda_i})\times C^{\nu_{\lambda_i}}
(c, c_+; J_{\lambda_i},
H_{\lambda_i}, A,\beta)\}=& &\\
\# \{ M(c_-, c; H_{\lambda_i})\}\times \#
\{C^{\nu_{\lambda_i}}(c, c_+; J_{\lambda_i}, H_{\lambda_i}; A,\beta)\}=0.
\end{eqnarray*}
This proves that the boundary operator, in the sense of algebraic topology,
gives 
\begin{eqnarray*}
\ptl C^\nu_\Lambda(c_-,c_+; J_\lambda, H_\lambda, A,\beta)& = &\\
C^{\lambda_+}(c_-, c_+; J_{\lambda_+}, H_{\lambda_+}, A, \beta) & &\\
-C^{\lambda_-}(c_-, c_+; J_{\lambda_-}, H_{\lambda_-}, A,\beta) & &
\end{eqnarray*}
as rational virtual moduli cycles. 

The invariance of $\Phi_{A, J_{\lambda}, H_{\lambda}}(c_{-}, c_{+})$ about $\lambda$ follows.

\QED

\noindent {\bf Proof of     Theorem  1.1 in the case of genus zero: }

Recall that $c_-$ and $c_+$ are the Morse theoretic representations of two
homology classes $\alpha_{-}$ and $\alpha_+$  with $supp (\alpha_-)\subset
V_+$ and $supp(\alpha_+)\subset V_-.$ We have assumed that the usual GW-invariant
$\Psi_{A, n+2}(\alpha_-, \alpha_+)\not=0,$ hence
$\Psi_{A, n+2}(\alpha_-, \alpha_+,$ $\beta_1,\cdots, \beta_n)\not=0$ for certain $\beta_{i}$'s. The usual
GW-invariant for general symplectic manifolds are established in the work
of [FO] and [LiT]. We refer readers to these references for the relevant definition.

As we mentioned in Theorem \ref{fifth} whose proof is in [LiuT3], that when 
$0<\lambda$ is small enough, we have
$$\Psi_A(\alpha_-,\alpha_+,\beta_1,\cdots,\beta_n)=\Phi_{A,J_\lambda,
H_\lambda}(c_-,c_+, \beta_1,\cdots,\beta_n).$$

Theorem \ref{invt} and our assumption imply that $\Phi^{\lambda}_{A,J_\lambda,
H_\lambda}\not=0$ for any $\lambda\in \Lambda.$ However, we will prove in a
moment 
\begin{lemma}\label{empty}
  
   When $\lambda \in \Lambda$ is large enough, the moduli
space ${\cal FM}(c_-, c_+; J_\lambda, H_\lambda, A)$ is empty. 

\end{lemma}

This implies
that $$\Phi_{A,J_\lambda, H_\lambda}(c_-, c_+, \beta_1,\cdots, \beta_n)=0$$
when $\lambda\in \Lambda$ is large enough. 
We get a contradiction. This implies that our
main assumption that $S$ and hence ${\tilde H}$ has no closed orbits is
not correct. This finishes the proof of Theorem 1.

To prove this last lemma, choose an element $f\in {\cal FM}(c_-, c_+; J_\lambda,
H_\lambda, A)$, and calculate its energy. Assume that $f=\cup_{i=1}^{N_P}
f^P_i\cup_{j=1}^{N_B}f^B_j$ with $f^P_i$ connecting critical points $c_i$ and
$c_{i+1}$ of $H_\lambda=\lambda\cdot H.$ Since $\sum_i [f_i^P]+\sum_j[f^B_j]
=A$ , we have 
\begin{eqnarray*}
0\leq E(f) &= &\sum_iE(f^P_i)+\sum_j E(f^B_j)   \\
& = &\sum_i (H_\lambda(c_{i+1})-H_\lambda(c_i)+\omega ([f^P_i]))+\sum_j \omega
([f^B_j])\\
& = & \lambda (H(c_+)-H(c_-))+\omega(A). \quad\quad\quad\llap(Main \,Estimate) 
\end{eqnarray*}
Therefore,  
$$\lambda\leq \frac{\omega(A)}{H(c_-)-H(c_+)}.$$ Now $H$ is $C^0$-closed
to the given ${\tilde H}$, which  implies that 
$$\frac{\omega(A)}{H(c_-)-H(c_+)}\leq\frac{\omega}{{\tilde H}(c_-)-{\tilde H}(c_+)}
+\frac{1}{2}=\frac{\omega(A)}{2\epsilon}+\frac{1}{2}.$$
Note that since $\Psi_A\not =0$, the class $A$ can be represented by some
$J$- holomorphic sphere. We have $\omega (A)>0.$  We conclude that if 
${\cal FM}(c_-, c_+; J_\lambda, H_\lambda, A)$ is not empty, then
$0<\lambda<\frac{\omega(A)}{2\epsilon}+\frac{1}{2}.$ This proves the lemma.

\QED

\section {Higher Genus Case}\label {gw1}
In this section, we will prove Theorem 1 in higher genus case. Our main observation 
here is that the main energy estimate can be carried out in a coordinate-free
manner. In particular, existence of a preferable cylindrical coordinate
of $S^2\setminus \{-\infty, +\infty\}$ in genus zero case does not play
any essential role in the estimate.

The method we developed in the previous sections can be adapted to higher 
genus to define the $GW$-invariants $\Psi_{A, g, n+2}$ (see [FO] and [LT]
for the details and other methods).

 To define higher genus perturbed $GW$-invariants $\Phi _{A, J,H, g, n+2},$
 we need to modify the definition of $(J,H)$-maps of genus zero case. For that purpose,
 we need to find    certain  $1$-forms on $(\Sigma, j)$ as  a replacement 
 of $ds$ and $d\theta.$ 
 Let ${\bar {\cal M}}_{g,n+2}$ be the Deligne-Mumford copmpactification of 
 ${\cal M}_{ g,n+2}$ of stable curves of genus $g$ with $n+2$ (ordered) marked points.
  We will  use $-\infty$ and $+\infty $ 
 to denote the first  two marked points. It is well-known that
 ${\bar {\cal M}}_{ g,n+2}$
is an orbifold. 
 We can  define the orbifold  bundle of closed $1$-forms with poles at double points,
 ${\cal CF}_{n+2}\rightarrow {\bar {\cal M}}_{g,n+
2}$  as follows.

Given  any $ \la \Sigma, j \ra
 \in {\bar {\cal M}}_{g, n+2},$ let $ (\Sigma,j) $  be a representative 
of it. The fiber $ ( {\cal CF}_{n+2})_{\la \Sigma,j \ra}$
$$ =\left\{ \xi| \xi \,\, \mbox{is  a closed 1-form over } \,\Sigma\setminus
 \{\mbox {double points},-\infty, +\infty\}, \,\, Res(\xi, d)=0\right \}$$ modulo the equivalence  relation  induced  by  the action of the automorphism group of $(\Sigma, j),$ where $d$ is a
 double point of $\Sigma$ and the residue $Res (\xi, d)$ is defined to be 
 $\int_{C_{+}}\xi+ \int_{C_{-}}\xi$. Here   $C_{+}$ and $C_{-}$ are the boundaries
 of small discs $D_{+}$and $D_{-}$ of $\Sigma$ centered at $d$, oriented  by
 the induced orientation of $D_{+}$ and $D_{-}.$  Note that here we consider  $-\infty$ and $+\infty$ as a single  double point.


 Given $(\Sigma, j)\in
\la \Sigma, j \ra \in 
 {\bar {\cal M}}_{g, n+2}$, let 
${\widetilde  W}_{\epsilon }(\Sigma)$  be a  local uniformizer
 of $ {\bar {\cal M}}_{g, n+2}$ 
near $\Sigma$ with cover group $\Gamma _{\Sigma}.$ 
 For simplicity, we may  assume that   $ \Sigma$ has  $3g-3+n +2 $  
 double points. Then the local deformation of  $ \Sigma$ can be completely described by  the gluing process as we did for genus zero case,  with
 the gluing parameter $t$ of $3g-3+n +2$ complex components. 
We will use $\Sigma_{t}$ to denote the corresponding deformation. 
Let $D_{l,+}$ and
$D_{l,-}$ be two fixed discs at a double point $d_{l}\in \Sigma.$ After gluing, we get a corresponding annulus  $C_{l, t}\cong S^1\times [0, L_{l}(t)]$ in $\Sigma_{t}.$
 By taking the intersection of the $\Gamma_{\Sigma}$-orbit of $C_{l,t}$,
we may assume that each $C_{l,t}$ is $\Gamma_{\Sigma}$-invariant. In general 
cylindrical coordinate $(s_l, \theta_l)$ along $C_{l,t}$ is only defined modulo the action of $\Gamma_{\Sigma}$. 
 Since the action preserves marked points, the 
 cylindrical coordinates near the two ends $-\infty $ and $+\infty$ are always well -defined for any elements in $ {\bar {\cal M}}_{g,n+ 2}.$

\begin{pro}

There exists a continuous section $Y^*$ of the bundle ${\cal CF}\rightarrow {\bar {\cal M}}_{g,n+2}$ such that the value of $Y^*$ (lifted to an
uniformizer)
at any $\Sigma_{t}
\in  {\widetilde  W}_{\epsilon }(\Sigma)$ has the property that
 $$Y^*( \Sigma_{t})|_{C_{l,t}}= a_{l}d\theta_{l},$$
 $$Y^*( \Sigma_{t})|_{C_{-\infty,t}}= d\theta_{-\infty} ,$$ 
 $$Y^*( \Sigma_{t})|_{C_{+\infty,t}}=d\theta_{+\infty}, $$  
 where $a_{l}$ is a constant. 
\end{pro}

\proof 

\noindent Step I. Construct $Y^*$ on each ${\widetilde W}_\epsilon (\Sigma)$.

We cosider the case that $n=0$ first.

Fix a smooth $\Sigma_{t_0}\in {\widetilde W}_\epsilon (\Sigma)$, let $C_l=C_l(t_0), $
$l=1,\cdots, L,$ be the all simple closed geodesic of $(\Sigma, m(j_0)),$ where 
$m(j_{0})$ is the hyperbolic metric which  corresponds to the complex structure on $\Sigma_{t_{0}}$. The
degeneracy from $\Sigma_{t_0}$ to $\Sigma$ can be described by shrinking 
 $3g-3 +2$ many of $C_l$'s of $\Sigma_{t_0}.$ For simplicity, we  assume that first $3g-3 +2$ 
 $C_l$'s are to be shrunk to double points of $\Sigma$. We may assume that 
 the length of $C_l(t)$ is bounded above for $\Sigma_t\in {\widetilde W}_\epsilon
 (\Sigma)$ and the upper bound has been achieved at $\Sigma_{t_0}$. Then
 each element $\Sigma'\in {\widetilde W}_\epsilon(\Sigma)$ can be obtained from 
 $\Sigma_{t_0}$ by cutting along some of $C_l, l=1,\cdots, L,$ and insert
 some neck $N_i\cong S^1\times [0, M_i]$ through gluing along boundaries. 
 Therefore we only need to construct $Y^*$ at $\Sigma_{t_0}$ with described
 property of theorem. To see this, we note that
by using the local translation invariance of 
 $a_l d\theta_l=Y^*(\Sigma_{t_0})|_{C_{l,t_0}}$ one can easily extend $Y^*$ over
 $N_i$, hence define $Y^*(\Sigma')$ for any $\Sigma'\in {\tilde W}_\epsilon(\Sigma).$

 To construct $Y^*(\Sigma_{t_0})$, we cut off $D_{+\infty}, D_{-\infty}$
of the neighbourhoods of $+\infty$ and $-\infty$  and
 glue their boundary $C_{+\infty}$ and $C_{-\infty}$ back to get a curve ${\tilde\Sigma}$
 of genus $g+1$ with $C_0=C_{+\infty}=C_{-\infty}\hookrightarrow {\tilde \Sigma}.$ By our
 assumption all $C_l\hookrightarrow\Sigma,  1\leq l \leq L,$ is still contained in ${\tilde \Sigma}$
 and $C_0\cap C_l=\emptyset. $ Let $e_1,\cdots, e_{2g+2}$ be the
 generators of $H_1({\tilde \Sigma}).$ 
 Write $C_l=\sum_j a_{lj}e_j, 0\leq l\leq L.$
 Since each $[C_l]\not =0$ in $H_1({\tilde \Sigma}), {\bf Z})$ , for any 
 fixed $l$, there exist some $j,$  $1\leq j\leq 2g+2$ such that $a_{l,j}\not=0$.
 Clearly, there exist some $x=(x_1,\cdots, x_{2g+2})$ such that 
 $\sum_j a_{0,j}x_j=1 $ and $\sum_j a_{l,j}x_j\not =0, 1\leq l\leq L.$
 By de Rham theorem, we can find a ${\tilde Y}^*$ such that 
 $\la {\tilde Y}^*, e_j\ra =x_j$. This implies that 
 $\la {\tilde Y}^*, C_0\ra =1$ and $\la {\tilde Y}^*, C_l\ra \not =0,$
$1\leq l\leq L.$ 
 By adding
 an exact one-form and back to $\Sigma_{t_{0}}$ , we can find 
 a ${\tilde Y}^*(\Sigma_{t_0})$
 with the desired property,  hence a section 
 ${\tilde Y}^*$ on ${\tilde W}_\epsilon(\Sigma).$
This completes the local construction for $n=0$. 

To deal with the general case, we may assume, for simplicity,  that
 $\Sigma_{t_{0}}=\Sigma^1_{t_{0}}\cup \Sigma^2_{t_{0}},$ where
 $\Sigma^1_{t_{0}}$ is a genus
$g$ curve with $k+2$ marked points, $k<n$, without any unstable rational components after forgetting its  last $k$ marked points, and $\Sigma^2_{t_{0}}$ consists of only stable  rational components (bubbles) carrying the rest $n-k$ marked points.
We may assume further that $\Sigma^2_{t_{0}}$ has only one component with two
marked points and one double point $d$. The local deformation of
 $\Sigma_{t_{0}}$ consists of two parts: the deformation of $\Sigma^1_{t_{0}}$ 
in $ {\bar {\cal { M}}}_{k+2}$, with the resulting surface $\Sigma^1_{t^1}$ associated to the   complex gluing parameter $t^1$ of $3g-3+k+2$ components, and a further deformation of  $\Sigma^1_{t^1}\cup \Sigma^2_{t_{0}}$ with a gluing parameter $t^2$  associated to the double point $d.$  It is easy to see that we can define $Y^*$ for all local deformation
$\Sigma^1_{t^1}$ by pulling back the $Y^*$ defined for the case $n=0$  through the  local projetion given by forgetting the last $k$ marked points. We  define $Y^*$ at $\Sigma^1_{t^1}\cup \Sigma^2_{t_{0}}$ by declaring its value to be zero along $\Sigma^2_{t_{0}}.$  Extending $Y^*$ over the final deformation given by  $t^2$ can be easily obtained. We leave it to the reader. 

 \noindent Step II. 

By taking the average of the action $\Gamma_\Sigma$ on ${\tilde Y}^*$, we get a
section $Y^*$ over $W_\epsilon(\Sigma)={\tilde W}_\epsilon(\Sigma)/\Gamma_\Sigma,$
with the desired property.

By using a partition of identity subject to a finite covering of ${\bar
{\cal M}}_{g,2}$ given by $\{{\tilde W}_{\epsilon_i}(\Sigma_i)/\Gamma_{\Sigma_i}\}$ , 
we can paste these local $Y^*$ defined on ${\tilde W}_{\epsilon_i}(\Sigma_i)$ 
together to get a well-defined $Y^*$ with the required property. 

\QED


What we need is slightly more. We need to define $Y^*({ \Sigma})$ for
any semi-stable curve ${\Sigma}$ of genus $g$. Each 
semi-stable curve $\Sigma$ can be
obtained from some stable curve $\Sigma'$ by inserting first some unstable principal
components $(P_l; (z_l)_-, (z_l)_+)$'s with $P_l\setminus \{(z_l)_-, (z_l)_+\}\cong
{\bf R}^1\times S^1$ at a double point or at the two ends of $\Sigma'$, then
adding some bubble components $B_j$'s. Clearly,  $Y^*$ can be extended in an obvious way to include all semi-stable curves in its domain.

 We will use $\Sigma_l$ to denote 
component $\Sigma'$ and write $\Sigma=\cup_{l}\Sigma_l\cup_{i} P_i\cup_{j} B_j.$
We now define $(j, J,H)$-map $f$ with domain $(\Sigma, j)$ of semi-stable curve of  
genus $g$ by using the
following equations:

\noindent (i) on $B_j$, $df^B_j+J(f^B_j)\circ df\circ j=0;$

\noindent (ii) on $\Sigma_l$, $df^\Sigma_l+J(f^\Sigma_l)\circ df\circ j
-\nabla H(f^\Sigma_l) j^*(Y^*)+J(f^\Sigma_l)\nabla H(f^\Sigma_l) Y^*=0;$

\noindent (iii) on $P_i$, same as (ii) for $f_i^P.$ 

As in genus zero case, we impose the obvious asymptotic condition along all
ends of $\Sigma.$ 

By using local convergence described in Section 1, together with the local
translation invariance of $Y^*$ along necks, we can prove the Gromov-Floer
compactness theorem for $(j,J, H)$-maps. We remark that each $(j, J,H)$-map
$f$ is defined on $\Sigma\setminus \{\mbox{double points}\}\cup 
\{-\infty,+\infty\}$ and we consider double points of $\Sigma$ and
$-\infty, +\infty$ as ends of $f$. Along those ends $f^\Sigma_l$ or $f^P_i$ 
are convergent to (successive) critical points of $\nabla H$, if the corresponding
$a_l \not= 0$. Otherwise, they are $(j, J)$-holomorphic along these ends, hence
can be extended smoothly over the double points.

Set
$X^*=j^* Y^*.$ Let $X, Y=j(X)$ be the dual vector fields of $X^*$ and $Y^*$. 
We define energy 
\begin{eqnarray*}
E(f) & = &\sum_i{\int\int}_{P_i} \la df^P_i(X), df^P_i(X)\ra X^*\wedge Y^*\\
 &+ &\sum
{\int\int}_{\Sigma_l}\la df^\Sigma_l(X), df^\Sigma_l(X)\ra  X^*\wedge Y^*
+\sum_j{\int\int}_{B_j} (f^B_j)^* \omega
\end{eqnarray*}

Since $X^*\wedge Y^*$ is compatible with the orientation of $P_i$ or $\Sigma_l$
at its non-zero points, we have $E(f)\geq 0.$ 

We can now recover our main estimate.
\begin{theorem}
 For $\lambda>\lambda_0+{1\over 2}$, there is no $(j, J_\lambda,
 H_\lambda)$-map connecting $c_-$ and $c_+$ of class $A$.
\end{theorem}

\proof

\begin{eqnarray*}
0\leq E(f) & =& {\int\int}_B f^*\omega +{\int\int}_{P\cup_l\Sigma_l} \omega
(df(X), df(Y))X^*\wedge Y^*  \\
& + & \lambda{\int\int}_{P\cup_l\Sigma_l} \la \nabla H,
df(X)\ra X^*\wedge Y^*\\
&= & {\int\int}_\Sigma f^*\omega+\lambda {\int\int}_{P\cup_l\Sigma_l}
(d(H\circ f)(X)X^*+d(H\circ f)(Y^*))\wedge Y^*\\
& = & {\int\int}_\Sigma f^*\omega +  \lambda{\int\int}_{P\cup_l\Sigma_l} d(H\circ
f\cdot Y^*)\\
 & = & \omega (A)+\lambda(H\circ f(+\infty)-H\circ f(-\infty))\\
& +& \lambda\sum_i (H\circ f(d_i) Res (Y^*(\Sigma), f(d_i)))\\
& = & \omega (A)+\lambda(H(c_+)-H(c_-)).
\end{eqnarray*}
This implies that 
$$\lambda\leq \frac{\omega(A)}{H(c_-)-H(c_+)}\leq \lambda_0.$$

\QED

By using $(j, J_{\lambda}, H_{\lambda})$-maps, we can now repeat our construction from Section 3 to Section 7 to 
get a virtual moduli cycle in higher genus case and to define perturbed
GW-invariant $\Phi_{A, J_\lambda, H_\lambda, g, n+2}$ for $\epsilon\leq
\lambda\leq \lambda_0+{1\over 2}.$ What is left is to prove the  higher genus
case of Theorem 1.6.

\begin{theorem}

$$\Phi_{A,J_\lambda, H_\lambda, g, n+2}(c_-, c_+, \beta_1,\cdots, \beta_n)$$
is well-defined, independent of the choices of $\lambda\in [\epsilon,
\lambda_0+{1\over 2}]$. It is equal to zero when $\lambda>\lambda_0.$
\end{theorem}
\proof 

The last statement follows from our estimate above. The first part can be
proved in a similar way as before. It fact, if there is no unstable principal
component $f^P_i$ appearing as connecting orbits, th invariance of $\Phi_{A, 
J_\lambda, H_\lambda, g, n+2}$
with respect to $\lambda$ can be proved in same way as the usual GW-invariant
$\Psi_{A,J_\lambda, g, n+2}.$ For each fixed intersection pattern, those
unstable principal components $f^P_i$ 
of connecting orbits form several strings, each connecting
two critical points of $H$. We may apply the theory for genus zero case to
each of those strings separately. The general case now follows by combining
above two cases.

\QED


\begin{thebibliography}{W}   



%
\bibitem[C]{c}
 W.Chen, 
 On  Weinstein conjecture for Stein Surfaces with Boundary,
{\it   Preprint} August (1997).
          
%
\bibitem[EH]{eh}
 Y. Eliashberg and  H. Hofer,
A Hamiltonian characterization of the three ball, 
{\it Differential Integral Equation} {\bf 7} (1994).

\bibitem[F]{f}
 A. Floer, Symplectic fixed points and holomorphic spheres, 
{\it Comm. Math. Phys. } {bf 120}(1989), pp. ~575-611. 

\bibitem[FH]{fh}
 A. Floer, H. Hofer,
 Coherent orientations for periodical orbit problems in symplectic
geometry.
                 
%
\bibitem[FHS]{fhs}
 A. Floer, H. Hofer and D. Salamon,
 Transversality in the elliptic Morse theory for the symplectic action,
 {\it Duke Math. Journal} {\bf 80}(1995), pp.~251-293.

%

\bibitem[FHV]{fhv}
 A. Floer, H. Hofer and C. Viterbo,
 The Weinstein conjecture in $P\times C^l$,
{\it  Mathematische Zeitschrift.} {\bf 203} (1989), pp. ~355-378.
%
\bibitem[FO]{fo}                        
Fukaya and Ono,
Arnold conjecture and Gromov-Witten invariants,
{\it Preprint} (1996)
%
\bibitem[G]{g}
M. Gromov,
Pseudo holomorphic curves in symplectic manifolds,
{\it Invent. Math.} {\bf 82} (1985), pp. ~307-347.
%
\bibitem[HV1]{hv1}
H. Hofer and C. Viterbo, 
 The Weinstein conjecture in the presence of $J$-holomorphic
 spheres,
{\it  Comm. Pure Appl. Math.} {\bf Vol. XLV} (1992), pp. ~583-622.


\bibitem[HV2]{hv2}
H. Hofer and C. Viterbo, 
 The Weinstein conjecture in cotangent bundles and related results,
{\it Ann. Scuola Norm. Sup. Pisa Cl. Sci.} (4) {\bf 15} (1988).


\bibitem[HZ]{hz}
H. Hofer and E. Zehnder,  Periodic solutions of hypersurfaces and
a result by C. Viterbo,
{\it  Inv. Math.} {\bf 90} (1987), pp. ~1-9.

\bibitem[HWZ]{hwz}
H. Hofer, K, Wysocki and E. Zehnder, A characterization of the tight three-sphere, {\it Duke Math.J.} {\bf 81} (1995).

\bibitem[L1]{l1}     %
G.Liu,                                     
Associativity of  quantum  multiplication, 
{\it  Preprint}  (1994), will appear in {\it Comm. Math. Phys.}





\bibitem[L2]{l2}     %
G.Liu,                                     
The equivalence of quantum cohomology and Floer cohomology, 
{\it  Preprint}  (1995).




%
\bibitem[LT]{lt}
J. Li and G. Tian,
Virtual moduli cycles and GW-invariants  of general symplectic manifolds,
{\it Preprint} (1996).


\bibitem[LiuT1]{liut1}
G. Liu and G. Tian,
 Floer homology and Arnold conjecture,
{\it  Preprint}  (1996).

\bibitem[LiuT2]{liut2}
G. Liu and G. Tian,
Weinstein Conjecture and GW- Invariants (announcement), 
{\it Preprint,} March (1997).

%

\bibitem[LiuT3]{liut3}
G. Liu and G. Tian,
On the equivalence of multiplicative structures in Floer and
quantum cohomology, in preparation.

\bibitem[PSS]{pss}
S. Piunikhin, D. Salamon and M. schwarz, Symplectic
Floer-Donoldson  theory and quantum cohomology,
{\it Preprint} (1994).      

\bibitem[R]{r}
P. Rabinowitz, Periodic solutions of Hamiltonian systems
{\it  Comm. Pure Appl Math. } {\bf 31}  ( 1979), pp. ~336-352.

\bibitem[RT]{rt}  
Y. Ruan and G. Tian, 
Bott-type   symplectic Floer cohomology and its multiplication structures,
{\it preprint } (1994).

\bibitem[V]{v}
 C. Viterbo, A proof of the Weinstein conjecture ${\bf {R}}^{2n}$
{\it  Ann. Inst.Henri Poincare, Analyse nonlineare, } {\bf 4}  ( 1987),
pp. ~337-356.
%
%
\bibitem[W]{w}
A. Weinstein,  
 Periodic orbits for convex Hamiltonian systems,
{\it Proc. Symp. Pure
Ann. Math.   } {\bf  108}  ( 1978). pp. ~507-518.
%
                 


%
\end{thebibliography}
\end{document}